\documentclass[mathpazo]{cicp}

\newcommand{\ud}{\textrm{d}}
\newcommand{\df}[2]{\frac{\partial #1}{\partial #2}}
\newcommand{\dd}[2]{\frac{\ud #1}{\ud #2}}
\newcommand{\con}{{\mathbf u}}
\newcommand{\ent}{{\mathbf v}}
\newcommand{\fx}{{\mathbf f}}
\newcommand{\gx}{{\mathbf g}}
\newcommand{\Dx}{{\mathbf D}}
\newcommand{\jph}{{j+\frac{1}{2}}}
\newcommand{\jmh}{{j-\frac{1}{2}}}
\newcommand{\tp}{\tilde{p}}
\newcommand{\tu}{\tilde{u}}
\newcommand{\avg}[1]{\overline{#1}}
\newcommand{\ept}{\varepsilon^{(2)}}
\newcommand{\epf}{\varepsilon^{(4)}}
\newcommand{\kapt}{\kappa^{(2)}}
\newcommand{\kapf}{\kappa^{(4)}}

\newcommand{\alog}[1]{\hat{#1}}

\begin{document}

\title{Kinetic energy preserving and entropy stable finite volume schemes for compressible Euler and Navier-Stokes equations}
\author[P.~Chandrashekar]{Praveen Chandrashekar}
\address{TIFR Center for Applicable Mathematics, Bangalore 560065, India}
\email{{\tt praveen@math.tifrbng.res.in} (P.~Chandrashekar)}

\begin{abstract}
Centered numerical fluxes can be constructed for compressible Euler equations which preserve kinetic energy in the semi-discrete finite volume scheme. The essential feature is that the momentum flux should be of the form $f^m_\jph = \tp_\jph + \avg{u}_\jph f^\rho_\jph$ where $\avg{u}_\jph = (u_j + u_{j+1})/2$ and $\tp_\jph, f^\rho_\jph$ are {\em any} consistent approximations to the pressure and the mass flux. This scheme thus leaves most terms in the numerical flux unspecified and various authors have used simple averaging. Here we enforce approximate or exact entropy consistency which leads to a unique choice of all the terms in the numerical fluxes. As a consequence novel entropy conservative flux that also preserves kinetic energy for the semi-discrete finite volume scheme has been proposed. These fluxes are centered and some dissipation has to be added if shocks are present or if the mesh is coarse. We construct scalar artificial dissipation terms which are kinetic energy stable and satisfy approximate/exact entropy condition. Secondly, we use entropy-variable based matrix dissipation flux which leads to kinetic energy and entropy stable schemes. These schemes are shown to be free of entropy violating solutions unlike the original Roe scheme. For hypersonic flows a blended scheme is proposed which gives carbuncle free solutions for blunt body flows. Numerical results for Euler and Navier-Stokes equations are presented to demonstrate the performance of the different schemes.\\
\end{abstract}


\maketitle

\section{Introduction}

The numerical solution of compressible Euler and Navier-Stokes (NS) equations by the finite volume method is now a routine task in many industries.  Due to their non-linear hyperbolic nature, solutions of Euler equations can be discontinuous with the presence of shocks or contact discontinuities. Discontinuous solutions must necessarily satisfy the Rankine-Hugoniot jump conditions which are a consequence of conservation. However it is well known that discontinuous solutions that satisfy the jump conditions can be still non-unique and an additional entropy condition has to be imposed in order to select the unique weak solution. In the case of Euler equations, there is a natural entropy condition which comes from the entropy condition in thermodynamics which must also be satisfied by the numerical scheme. Additionally other global balance equations like that for the total kinetic energy must also be consistently approximated by the numerical solutions. The finite volume method requires the computation of the inviscid and viscous fluxes across the boundaries of the finite volumes. The design of these fluxes must incorporate the properties of the Euler/NS equations like entropy condition and kinetic energy preservation. There exists a vast library of numerical flux functions for the Euler equations and some of these like the Godunov scheme and kinetic scheme can be shown to satisfy the entropy condition.  The popular Roe scheme~\cite{Roe1981357} does not satisfy the entropy condition and can give rise to entropy violating shocks near sonic points. Various entropy fixes for Roe scheme have been proposed which involve preventing the numerical dissipation from vanishing at sonic points. Tadmor~\cite{tadmor1987} proposed the idea of entropy conservative numerical fluxes which can then be combined with some dissipation terms using entropy variables to obtain a scheme that respects the entropy condition, i.e., the scheme must produce entropy in accordance with the second law of thermodynamics. This is a more mathematically rigorous approach to construct entropy stable schemes for conservation laws. However some of these entropy conservative numerical fluxes have to be computed with quadrature rules since the integrals involved in the definition of the flux cannot be evaluated explicitly. For the Euler equations, Roe proposed explicit entropy conservative numerical fluxes~\cite{roehyp2006,Ismail:2009:AEE:1552569.1552681} which are augmented by Roe-type dissipation terms using entropy variables. These schemes do not suffer from entropy violating solutions that are observed in the original Roe scheme. However for strong shocks, even the first order schemes can produce oscillations indicating that the amount of numerical dissipation is not sufficient. Roe~\cite{Ismail:2009:AEE:1552569.1552681} proposed modifying the eigenvalues of the dissipation matrix which lead to non-oscillatory solutions. The modification of the eigenvalues is such that the amount of entropy production is of the correct order of magnitude for weak shocks.  The availability of cheap entropy conservative fluxes allows us to use the procedure of~\cite{doi:10.1137/S003614290240069X} to develop high order accurate entropy conservative schemes. Matrix dissipation can be added following the ENO procedure of~\cite{doi:10.1137/110836961} to develop arbitrarily high order accurate entropy stable schemes for the Euler equations on structured grids.

Faithful representation of kinetic energy evolution is another desirable property of a numerical scheme~\cite{Jameson:2008:FKE:1331553.1331559}. This is important for direct numerical simulation of turbulent flows where the kinetic energy balance plays an important role in the evolution of turbulence~\cite{Morinishi2010276,Shoeybi20105944,doi:10.1146/annurev-fluid-122109-160718}. The scheme is also stable in the sense that spurious kinetic energy is not produced by the numerical fluxes. Godunov schemes are found to have wrong order of kinetic energy dissipation and entropy production which leads to excessive damping of flow structures~\cite{Thornber20084853}. One way to construct kinetic energy preserving schemes is to use a skew symmetric form for the non-linear convective terms combined with finite difference operators which satisfy summation by parts property~\cite{Blaisdell1996207,Honein2004531,Kennedy20081676}. The skew symmetric form is also found to yield smaller aliasing errors which makes them more robust~\cite{Blaisdell1996207}. For compressible simulations, it is necessary to also satisfy conservation of mass, momentum and energy and hence it is preferable to solve the Navier-Stokes equations in conservation form using the finite volume method. The essential feature for a numerical flux in a semi-discrete finite volume method to correctly capture the kinetic energy balance is that the momentum flux should be of the form $f^m_\jph = \tp_\jph + \avg{u}_\jph f^\rho_\jph$ where $\avg{u}_\jph = (u_j + u_{j+1})/2$ and $\tp_\jph, f^\rho_\jph$ are {\em any} consistent approximations to the pressure and the mass flux. This scheme thus leaves most terms in the numerical flux unspecified and various authors have used simple averaging.  Subbareddy and Candler~\cite{Subbareddy20091347} have proposed a fully discrete finite volume scheme for the compressible Euler equations which preserves kinetic energy but the resulting scheme is implicit. All of these kinetic energy preserving schemes are however not entropy conservative, while on the other hand, the entropy conservative schemes do not have the kinetic energy preservation property. It is thought that for DNS of compressible flows, a numerical scheme which preserves kinetic energy and satisfies entropy condition is desirable since such schemes would be non-linearly stable. Schemes which satisfy entropy condition are found to lead to stable density fluctuations in compressible isotropic turbulence simulations, while schemes which do not have this property can be unstable with respect to density fluctuations~\cite{Honein2004531,doi:10.1146/annurev-fluid-122109-160718}. Using a canonical splitting of the flux with exponential entropy, Gerritsen and Olsson~\cite{Gerritsen1996245} construct entropy stable schemes which are however not conservative with respect to mass, momentum and energy. In~\cite{Honein2004531}, skew-symmetric form of convective terms is used to enforce better entropy consistency which is then shown to lead to schemes capable of computing high Reynolds number turbulence without artificial dissipation or filtering.

In the present work, we construct explicit centered numerical fluxes for the compressible Euler equations which are entropy conservative and also preserve kinetic energy in the case of the semi-discrete finite volume scheme. All the numerical fluxes presented here can be used in any finite volume scheme on structured or unstructured grids. Two versions of the numerical flux are constructed, one of which is only approximately entropy consistent but has simpler expressions, while the second one is exactly entropy conservative but involves certain logarithmic averages requiring more computations. Due to lack of upwinding, the schemes  are not stable for discontinuous solutions and for NS equations on coarse meshes for which shocks may not be well resolved. They yield stable solutions for Navier-Stokes equations when used on very fine meshes where the physical viscosity is enough to stabilize the scheme. However for Euler equations and for NS equations on coarse meshes, the centered fluxes are unstable and must be augmented with dissipation terms. Firstly, we construct scalar dissipation terms using second and fourth order differences as in the JST scheme~\cite{jst}. The second order dissipation terms which are active near shocks are kinetic energy and entropy stable while the fourth order dissipation is active only in smooth regions of the flow. Secondly, we also use entropy variable based matrix dissipation flux similar to the Roe scheme which can be shown to lead to entropy generation~\cite{tadmor1987}. The eigenvalue modification of Roe~\cite{Ismail:2009:AEE:1552569.1552681} is used to compute strong shocks without oscillations. All the schemes are shown to give entropy consistent solutions in cases where the Roe scheme would give entropy violating shocks. The matrix dissipation scheme can also be made kinetic energy dissipative by a proper choice of the eigenvalues appearing in the dissipation matrix. The entropy stable schemes with matrix dissipation preserve stationary contacts exactly but also suffer from 1-D shock instability and the carbuncle phenomenon. A modification of the eigenvalues in the dissipation flux is suggested which avoids these problems but is still able to accurately compute shear flows like boundary layers. The performance of the schemes is demonstrated on inviscid and viscous test cases involving shocks and contact discontinuities.

The rest of the paper is organized as follows. Section~(\ref{sec:1d}) introduces the 1-D Navier-Stokes equations and finite volume method and section~(\ref{sec:kep}) discusses the kinetic energy preservation property. The entropy condition is introduced and the entropy conservative fluxes are derived in section~(\ref{sec:entropy}) together with some justification for their uniqueness. Scalar artificial dissipation terms which satisfy kinetic energy and entropy stability are derived in section~(\ref{sec:scal}) while matrix dissipation flux is treated in section~(\ref{sec:mat}). Section~(\ref{sec:monotone}) shows numerical results for stationary 1-D shock problem and modifications of the scheme for monotone resolution of shocks, while section~(\ref{sec:hyb}) introduces a hybrid scheme. The 2-D equations and numerical fluxes are discussed in section~(\ref{sec:2d}) and section~(\ref{sec:results}) ends with some numerical results on a range of test problems.

\section{1-D NS equations and finite volume scheme}
\label{sec:1d}
The one dimensional Navier-Stokes equations can be written in vector conservation form as
\begin{equation}
\df{\con}{t} + \df{\fx}{x} = \df{\gx}{x}
\end{equation}
where $\con$ is the set of conserved variables and $\fx$, $\gx$ are the inviscid and viscous fluxes whose expressions are given by
\begin{equation}
\con = \begin{bmatrix}
\rho \\
\rho u \\
E
\end{bmatrix} = \begin{bmatrix}
\rho \\
m \\
E
\end{bmatrix}, \qquad \fx = \begin{bmatrix}
\rho u \\
p + \rho u^2 \\
(E+p)u
\end{bmatrix} = \begin{bmatrix}
m \\
p + um \\
(E+p)u
\end{bmatrix}, \qquad \gx = \begin{bmatrix}
0 \\
\tau \\
u \tau - q
\end{bmatrix}
\end{equation}
In the above equations $\rho$ is the density, $u$ is the velocity, $p$ is the pressure and $E$ is the total energy per unit volume which for a perfect gas is given by $E=p/(\gamma-1) + \rho u^2/2$, where $\gamma$ is the ratio of specific heats which is taken to be constant. Moreover, $\tau$, $q$ are the shear stress and heat flux for which we take the Newtonian and Fourier laws respectively, leading to
\[
\tau = \frac{4}{3}\mu \df{u}{x}, \qquad q = -\kappa \df{T}{x}
\]
where $\mu$, $\kappa$ are the coefficient dynamic viscosity and heat conduction respectively. In the absence of the viscous fluxes $\gx$, the resulting Euler equations form a hyperbolic system of conservation laws. Consider a partition of the one dimensional domain into uniform finite volumes of size $\Delta x$ and the $j$'th cell is the interval $(x_\jmh,x_\jph)$. The semi-discrete finite volume scheme is given by
\begin{equation}
\Delta x \dd{\con_j}{t} + \fx_\jph - \fx_\jmh = \gx_\jph - \gx_\jmh
\label{eq:sdfvm}
\end{equation}
where $\con_j$ is the cell average value in the $j$'th cell and $\fx_\jph$, $\gx_\jph$ are numerical inviscid and viscous fluxes respectively at the interface $x_\jph$. In the numerical computations, the above set of ordinary differential equations will be solved using a strong stability preserving Runge-Kutta scheme~\cite{Shu1988439}.
\section{Kinetic energy preserving scheme}
\label{sec:kep}
The kinetic energy is an important quantity in fluid flows and it is destroyed by the physical viscosity. In turbulent flows, the kinetic energy injected into the fluid at large scales cascades to smaller scales and is eventually destroyed by viscosity. Hence it is desirable that the numerical scheme faithfully represent the kinetic energy balance consistent with the Navier-Stokes equations. The kinetic energy per unit volume $K = \frac{1}{2}\rho u^2$ satisfies the following equation
\begin{equation}
\dd{}{t}\int_\Omega K \ud x =  \int_\Omega p \df{u}{x} \ud x - \frac{4}{3} \int_\Omega \mu \left( \df{u}{x} \right)^2 \ud x
\label{eq:kecont}
\end{equation}
where we have ignored boundary conditions. The first term on the right is the rate at which work is done by pressure forces and is present only for compressible flows. The second term represents the irreversible destruction of kinetic energy which is converted into internal energy due to viscous dissipation. We would like the numerical scheme to also satisfy this equation in a discrete sense which will then be refered to as kinetic energy preserving scheme. Consider the following approximation for the inviscid and viscous fluxes
\begin{equation}
\fx_\jph = \begin{bmatrix}
f^\rho \\
f^m \\
f^e
\end{bmatrix}_\jph = \begin{bmatrix}
f^\rho \\
\tp + \avg{u} f^\rho \\
f^e
\end{bmatrix}_\jph, \qquad g_\jph = \begin{bmatrix}
0 \\
\tau \\
\tilde{u} \tau - q
\end{bmatrix}_\jph
\end{equation}
where
\[
\avg{u}_\jph = \frac{1}{2}(u_j + u_{j+1}), \qquad \tau_\jph = \frac{4}{3}\mu \frac{u_{j+1} - u_j}{\Delta x}, \qquad q_\jph = -\kappa \frac{T_{j+1} - T_j}{\Delta x}
\]
Throughout the paper, we will use the overbar to denote the arithmetic average. The quantities $\tp, f^\rho, \tilde{u}, f^e$ are assumed to be consistent approximations but are yet to be specified. The kinetic energy varies in time according to
\begin{equation}
\df{K}{t} = -\frac{1}{2}u^2 \df{\rho}{t} + u \df{m}{t}
\end{equation}
We now derive the global kinetic energy balance equation for the finite volume scheme by adding the equation from each cell and summing over all cells in the grid.
\begin{eqnarray*}
\sum_j \Delta x \dd{K_j}{t} &=& \sum_j \left[ -\frac{1}{2} u_j^2 \dd{\rho_j}{t} + u_j \dd{m_j}{t} \right] \Delta x \\
&=& \sum_j \left[ \frac{1}{2}u_j^2 (f^\rho_\jph - f^\rho_\jmh) - u_j (f^m_\jph - f^m_\jmh) + u_j (g^m_\jph - g^m_\jmh) \right] \\
&=& \sum_j \left[ \frac{1}{2}(u_j^2 - u_{j+1}^2) f^\rho_\jph - (u_j - u_{j+1}) f^m_\jph + (u_j - u_{j+1}) g^m_\jph \right] \\
&=& \sum_j \left[ (u_j - u_{j+1})\left( \frac{1}{2}(u_j + u_{j+1}) f^\rho_\jph - f^m_\jph \right) + (u_j - u_{j+1}) g^m_\jph \right] \\
&=& \sum_j \left[ \frac{\Delta u_\jph}{\Delta x} \tp_\jph  - \frac{4}{3}\mu \left(\frac{\Delta u_\jph}{\Delta x} \right)^2 \right] \Delta x
\end{eqnarray*}
This is consistent with the continuous kinetic energy equation. The crucial property used in the above proof was that the momentum flux has the form $f^m = \tp + \avg{u} f^\rho$ which leads to the disappearance of the convective flux from the kinetic energy equation. However, we still have freedom in the choice of $f^\rho, \tp, \tilde{u}, f^e$. Jameson~\cite{Jameson:2008:FKE:1331553.1331559} makes the following choice:
\[
f^\rho = \bar{\rho} \bar{u}, \qquad \tp = \bar{p}, \qquad f^e = \bar{\rho} \bar{u} \bar{H}
\]
where $H$ is the enthalpy. We will refer to this as the KEP flux. In the present work we will determine the numerical flux in a unique manner from entropy considerations.

\section{Entropy condition}
\label{sec:entropy}
The concept of entropy condition is borrowed from the second law of thermodynamics and generalized to any arbitrary system of hyperbolic conservation laws~\cite{Amiram1983151}. Let $U(\con)$ be a strictly convex function; then $U(\con)$, $F(\con)$ is said to be an entropy-entropy flux pair if they satisfy
\[
U'(\con) \fx'(\con) = F'(\con)
\]
Then smooth solutions of the inviscid equations satisfy an additional conservation law
\[
\df{U}{t} + \df{F}{x} = 0
\]
But for discontinuous solutions we can only satisfy the entropy inequality
\begin{equation}
\df{U}{t} + \df{F}{x} \le 0
\label{eq:econd}
\end{equation}
where the inequality is satisfied in a weak sense. Define the entropy variables as
\[
\ent(\con) = U'(\con)
\]
and since $U(\con)$ is strictly convex, we can invert the above equation to obtain $\con = \con(\ent)$. Define the quantity $\psi(\ent)$ which is the dual of the entropy flux $F(\con)$ by the relation
\[
\psi(\ent) = \ent \cdot \fx(\con(\ent)) - F(\con(\ent))
\]
The finite volume scheme for the hyperbolic conservation law is given by
\[
\Delta x \dd{\con_j}{t} + \fx_\jph - \fx_\jmh = 0
\]
Taking the dot product of the above equation with $\ent_j$, we obtain the entropy equation
\[
\Delta x \dd{U_j}{t} + \ent_j \cdot (\fx_\jph - \fx_\jmh) = 0
\]
Tadmor~\cite{tadmor1987,tadmor2004} introduced the idea of an {\em entropy conservative numerical flux} which should satisfy the following condition
\begin{equation}
(\ent_{j+1} - \ent_j) \cdot \fx_\jph = \psi_{j+1} - \psi_j
\end{equation}
Then we can write
\[
\ent_j \cdot (\fx_\jph - \fx_\jmh) = F_\jph - F_\jmh
\]
where
\[
F_\jph = \avg{\ent}_\jph \cdot \fx_\jph - \avg{\psi}_\jph
\]
is a consistent numerical entropy flux. The semi-discrete scheme satisfies the entropy conservation equation
\[
\Delta x \dd{U_j}{t} + F_\jph - F_\jmh = 0
\]
An entropy conservative flux is given by Tadmor~\cite{tadmor1987} as
\[
\fx_\jph = \int_0^1 \fx(\ent_\jph(\theta)) \ud \theta, \qquad \ent_\jph(\theta) = \ent_j + \theta(\ent_{j+1} - \ent_j)
\]
which usually requires some numerical quadrature to approximate the flux. Once an entropy conservative flux has been constructed, we can add dissipative terms to the flux which leads to the satisfaction of the entropy condition as given by equation~(\ref{eq:econd}), i.e., the dissipative flux must lead to generation of entropy.
\subsection{Entropy condition for Euler equations}
For the Euler equations we can take the entropy-entropy flux pair to be
\begin{equation}
U = -\frac{\rho s}{\gamma-1}, \qquad F = -\frac{\rho u s}{\gamma-1}
\label{eq:eulentropy}
\end{equation}
where $s$ is the physical entropy given by
\begin{equation}
s=\ln(p) - \ln(\rho^\gamma) + \textrm{const} = -(\gamma-1) \ln(\rho) - \ln(\beta) + \textrm{const}, \qquad \beta = \frac{1}{2RT}
\end{equation}
and the constant term can be ignored. There are many other possible choices for the entropy function $U(\con)=\rho \eta(s)$ where $\eta$ is any convex function~\cite{Amiram1983151}, but the  above choice is the only one which is consistent with the entropy condition from thermodynamics~\cite{Hughes1986223} in the presence of heat transfer. Since we work with the correct choice of the entropy function, the schemes we develop will satisfy entropy condition for the Navier-Stokes equations also. The entropy variables $\ent$ and the Legendre transform $\psi$ are given by
\begin{equation}
\ent = \begin{bmatrix}
\frac{\gamma - s}{\gamma-1} - \beta u^2 \\
2\beta u \\
-2\beta
\end{bmatrix}, \qquad \psi = m = \rho u
\end{equation}
Hence an entropy conservative numerical flux for the Euler equations has to satisfy the following condition
\begin{equation}
(\ent_{j+1} - \ent_j) \cdot \fx_\jph = m_{j+1} - m_j
\label{eq:tadcond}
\end{equation}
This provides only one equation whereas the flux $\fx$ has more than one component; we can expect that there are many possible entropy conservative fluxes. 
\subsection{Roe's entropy conservative flux}
\label{sec:roe}
Roe~\cite{roehyp2006} has constructed explicit numerical fluxes for the Euler equations which satisfy condition~(\ref{eq:tadcond}) for the entropy given in equation~(\ref{eq:eulentropy}) but do not have the kinetic energy preservation property. Introduce the set of independent state variables or parameter vector
\[
z = \begin{bmatrix} z_1 \\ z_2 \\ z_3 \end{bmatrix} = \sqrt{\frac{\rho}{p}} \begin{bmatrix}
1 \\
u \\
p
\end{bmatrix}
\]
Roe also introduces the logarithmic average $\alog{\varphi}$ of any strictly positive quantity $\varphi_l$, $\varphi_r$ which is defined as
\[
\alog{\varphi} = \frac{\varphi_r - \varphi_l}{\ln \varphi_r - \ln \varphi_l} = \frac{\Delta\varphi}{\Delta\ln\varphi}
\]
A numerically stable procedure to compute the average when $\varphi_l \approx \varphi_r$ is given in~\cite{Ismail:2009:AEE:1552569.1552681}. Then the entropy conservative numerical flux at any cell face $\jph$  is given by
\[
\fx^* = \begin{bmatrix}
\tilde{\rho} \tilde{u} \\
\tilde{p}_1 + \tu f^\rho \\
\tilde{H} f^\rho
\end{bmatrix}
\]
where
\[
\tilde{\rho}= \avg{z}_1 \alog{z}_3, \quad \tu = \frac{\avg{z}_2}{\avg{z}_1}, \quad \tilde{p}_1 = \frac{\avg{z}_3}{\avg{z}_1}, \quad \tilde{p}_2 = \frac{\gamma + 1}{2\gamma} \frac{\alog{z}_3}{\alog{z}_1} + \frac{\gamma-1}{2\gamma} \frac{\avg{z}_3}{\avg{z}_1}
\]
\[
\tilde{a} = \left( \frac{\gamma \tilde{p}_2}{\tilde{\rho}} \right)^\frac{1}{2}, \quad \tilde{H} = \frac{\tilde{a}^2}{\gamma-1} + \frac{1}{2} \tu^2
\]
and all the averages are evaluated using the state in cell $j$ and $j+1$. This flux is entropy conservative and it can be made entropy stable~\cite{roehyp2006} by adding entropy variable based matrix dissipation terms which is described in a later section. The resulting numerical flux will be refered to as ROE-ES flux. Note that the momentum flux contains the weighted average velocity $\tu$ while kinetic energy preservation requires the presence of the arithmetic average $\avg{u}$. Hence the above entropy conservative flux is not kinetic energy preserving. In the following sections, we construct numerical fluxes by approximately or exactly satisfying the entropy condition for the Euler equations which also preserves the kinetic energy.
\subsection{Derivation of approximately entropy consistent flux}
We now derive kinetic energy preserving flux which is approximately entropy consistent. In order to simplify the notations, we will drop the subscripts on the cell indices, and use the convention that $\Delta(\cdot) = (\cdot)_{j+1} - (\cdot)_j$ denotes the jump across the cell face $\jph$. Equation (\ref{eq:tadcond}) can then be written as
\[
\left[ -\frac{1}{\gamma-1} \Delta s - \Delta(\beta u^2) \right] f^\rho + 2 \Delta(\beta u) (\tp + \avg{u} f^\rho) - 2 \Delta(\beta) f^e = \Delta (\rho u)
\]
We can easily check that the following exact linearizations are true.
\[
\Delta(\beta u) = \avg{\beta} \Delta u + \avg{u} \Delta \beta, \quad
\Delta(\rho u) = \avg{\rho} \Delta u + \avg{u} \Delta \rho, \quad
\Delta(\beta u^2) = \avg{\beta} \Delta u^2 + \avg{u^2} \Delta \beta = 2 \avg{\beta} \avg{u} \Delta u + \avg{u^2} \Delta\beta
\]
The entropy difference is
\[
\Delta s = -(\gamma-1) \Delta \ln\rho - \Delta \ln\beta
\]
Assuming smooth solutions we approximate the above differences by a Taylor formula, e.g.,
\begin{eqnarray*}
\Delta \ln \rho &=& \ln(\avg{\rho} + \frac{1}{2}\Delta\rho) - \ln(\avg{\rho} - \frac{1}{2}\Delta\rho) = \ln(1 + \frac{\Delta\rho}{2\avg{\rho}}) - \ln(1 - \frac{\Delta\rho}{2\avg{\rho}}) \\
&=& \left[ \frac{\Delta\rho}{2\avg{\rho}} - \frac{1}{2} \left(\frac{\Delta\rho}{2\avg{\rho}}\right)^2 + O(\Delta\rho)^3 \right] - \left[ -\frac{\Delta\rho}{2\avg{\rho}} - \frac{1}{2} \left(\frac{\Delta\rho}{2\avg{\rho}}\right)^2 + O(\Delta\rho)^3 \right] \\
&=& \frac{\Delta\rho}{\avg{\rho}} + O(\Delta\rho)^3
\end{eqnarray*}
and similarly for $\Delta\ln\beta$. We then obtain
\[
-\frac{1}{\gamma-1}\Delta s = \frac{\Delta\rho}{\avg{\rho}} + \frac{1}{\gamma-1} \frac{\Delta\beta}{\avg{\beta}} + O(\Delta\rho)^3 + O(\Delta\beta)^3
\]
For smooth solutions, the higher order terms are of $O(\Delta x^3)$. Neglecting the third order terms, we try to satisfy condition~(\ref{eq:tadcond}) which becomes
\[
 \left[ \frac{f^\rho}{\avg{\rho}} - \avg{u} \right] \Delta\rho + \left[ -2 \avg{\beta} \avg{u} f^\rho + 2  \avg{\beta} (\tp + \avg{u} f^\rho) - \avg{\rho} \right] \Delta u + \left[ \frac{1}{\gamma-1} \frac{f^\rho}{\avg{\beta}} - \avg{u^2} f^\rho + 2 \avg{u} f^m - 2 f^e\right] \Delta\beta = 0
\]
Since $\rho, u, \beta$ are independent variables, the above equation is satisfied if we choose
\[
f^\rho = \avg{\rho} \ \avg{u}, \qquad \tp = \frac{\avg{\rho}}{2 \avg{\beta}}
\]
and
\[
f^e = \frac{1}{\gamma-1} \frac{f^\rho}{2\avg{\beta}} - \frac{1}{2}\avg{u^2} f^\rho + \avg{u} f^m = \left[ \frac{1}{2(\gamma-1) \avg{\beta}} - \frac{1}{2} \avg{u^2} \right] f^\rho + \avg{u} f^m
\]
With the above unique choice, the numerical flux satisfies condition~(\ref{eq:tadcond}) to third order accuracy, i.e.,
\begin{equation*}
(\ent_{j+1} - \ent_j) \cdot \fx_\jph = m_{j+1} - m_j + O(\Delta\rho)^3 + O(\Delta\beta)^3
\end{equation*}
It is easy to check that the above numerical fluxes are consistent. We remark that these flux formulae have simpler expressions since they make use of arithmetic averages while the entropy conservative fluxes have slightly more complicated expressions with logarithmic averages as we see from Roe's fluxes and also the new flux derived in the next section.

\subsection{Kinetic energy preserving and entropy conservative flux}
In this section we derive numerical fluxes which preserves the kinetic energy and exactly conserves the entropy. The jump in the entropy variables can be written as
\begin{eqnarray*}
\Delta v_1 &=& \Delta\ln(\rho) + \frac{1}{\gamma-1}\Delta\ln(\beta) - \Delta(\beta u^2) = \frac{\Delta\rho}{\alog{\rho}} + \left[\frac{1}{(\gamma-1)\alog{\beta}} - \avg{u^2} \right] \Delta\beta - 2 \avg{u} \avg{\beta} \Delta u \\
\Delta v_2 &=& 2 \avg{\beta} \Delta u + 2 \avg{u} \Delta \beta \\
\Delta v_3 &=& -2\Delta\beta
\end{eqnarray*}
where we have made use of the logarithmic averages $\alog{\rho}$ and $\alog{\beta}$ as introduced in section~(\ref{sec:roe}). The flux is entropy conservative if equation~(\ref{eq:tadcond}) is satisfied, i.e., if
\begin{equation}
f^\rho \Delta v_1 + f^m \Delta v_2 + f^e \Delta v_3 = \Delta(\rho u) = \avg{\rho} \Delta u + \avg{u} \Delta\rho
\label{eq:tadcond2}
\end{equation}
As in the previous section, the above equation is a finite difference equation involving differences in the independent variables $\rho,u,\beta$. Equating the terms containing $\Delta\rho$ on both sides yields the mass flux as
\[
f^\rho = \alog{\rho} \avg{u}
\]
Next equating the terms containing $\Delta u$ we get the momentum flux as
\[
-2 \avg{u} \avg{\beta} f^\rho + 2 \avg{\beta} f^m = \avg{\rho} \quad\Longrightarrow\quad f^m = \frac{\avg{\rho}}{2\avg{\beta}} + \avg{u} f^\rho = \tp + \avg{u} f^\rho \quad\textrm{where}\quad \tp = \frac{\avg{\rho}}{2\avg{\beta}}
\]
The above form of the momentum flux satisfies the requirement for kinetic energy preservation. Finally, from the terms containing $\Delta\beta$, we get the energy flux
\[
f^e = \left[ \frac{1}{2(\gamma-1) \alog{\beta}} - \frac{1}{2} \avg{u^2} \right] f^\rho + \avg{u} f^m
\]
These fluxes are consistent and have exactly the same form as the approximately consistent fluxes derived in the previous section except that the logarithmic averages $\alog{\rho}$, $\alog{\beta}$ are used instead of the arithmetic averages $\avg{\rho}$, $\avg{\beta}$ in the mass and energy fluxes. We can also see that the new entropy conservative flux is computationally less expensive as compared to Roe's entropy conservative flux given in section~(\ref{sec:roe}) since it does not require a parameter vector in its definition.
\subsection{Uniqueness of the numerical flux}
Roe~\cite{roehyp2006} has proposed explicit entropy conservative flux for Euler equations and in the present work we have proposed another one. Tadmor~\cite{tadmor1987} gives a recipe for the construction of entropy conservative fluxes for any hyperbolic system endowed with an entropy-entropy flux pair.  However, the new flux function proposed here also satisfies kinetic energy preservation. Since there are many numerical flux functions that satisfy entropy conservation, it is interesting to ask if there is a unique numerical flux function that satisfies both conservation properties and here we give some justification in favour of this for the case of two point numerical fluxes. The crucial property for kinetic energy preservation is that the momentum flux must be of the form $f^m = \tp + \avg{u} f^\rho$ where $\tp$ is any consistent approximation to the pressure and $f^\rho$ is any consistent approximation to the mass flux. Looking at equation~(\ref{eq:tadcond2}), we see that the term $\avg{u} f^\rho$ can only come from the term $\Delta(\beta u^2)$ in $\Delta v_1$. In the above derivation, we have written $\Delta(\beta u^2) = \avg{u^2}\Delta\beta + \avg{\beta} \Delta u^2 = \avg{u^2}\Delta\beta + 2 \avg{\beta} \avg{u} \Delta u$ which gives the desired form in the momentum flux but there are other ways to do the linearization. If instead we write this term as 
\[
\Delta(\beta u^2) = \avg{\beta u}\Delta u + \avg{u}\Delta(\beta u) = (\avg{\beta u} + \avg{u} \avg{\beta}) \Delta u + \avg{u}^2 \Delta \beta
\]
then we do not get the correct form of the momentum flux. Roe writes it in the form $\Delta(\beta u^2) = \Delta(\sqrt{\beta}u)^2 = 2 \avg{\sqrt{\beta}u} \Delta( \sqrt{\beta}u)$ which also does not have the correct form for kinetic energy preservation.

The other type of non-uniqueness can arise with the choice of independent state variables. In the derivation of the new flux, we wrote condition~(\ref{eq:tadcond}) in terms of jumps in $\rho$, $u$, $\beta$. The other possibility is to use the variables $\rho$, $u$, $p$ and we show that this does not lead to a proper definition of the numerical fluxes. The jump in the physical entropy at a cell face $\jph$ can be written as
\[
\Delta s = \Delta \ln p - \gamma \Delta \ln \rho = \frac{\Delta p}{\alog{p}} - \gamma \frac{\Delta \rho}{\alog{\rho}}
\]
and the jump in $\beta$ can be written as
\[
2\Delta\beta = \Delta (\rho/p) = \avg{r}\Delta\rho - \avg{\rho} \frac{\Delta p}{p_j p_{j+1}}, \qquad r = \frac{1}{p}
\]
Then the jump in the entropy variables is given by
\begin{eqnarray*}
\Delta v_1 &=& \left[ -\frac{1}{(\gamma-1)\alog{p}} + \frac{\avg{\rho} \avg{u^2}}{2p_j p_{j+1}} \right] \Delta p + \left[ \frac{\gamma}{(\gamma-1)\alog{\rho}} - \frac{1}{2} \avg{r} \avg{u^2} \right] \Delta\rho - 2 \avg{u} \avg{\beta} \Delta u \\
\Delta v_2 &=& 2 \avg{\beta} \Delta u + \avg{u} \ \avg{r}\Delta\rho - \avg{\rho} \ \avg{u} \frac{\Delta p}{p_j p_{j+1}} \\
\Delta v_3 &=& - \avg{r}\Delta\rho + \avg{\rho} \frac{\Delta p}{p_j p_{j+1}}
\end{eqnarray*}
If we now use condition (\ref{eq:tadcond}), then the terms containing $\Delta u$ lead to the momentum flux
\[
f^m = \frac{ \avg{\rho}}{2\avg{\beta}} + \avg{u} f^\rho
\]
which satisfies the condition for kinetic energy preservation. The terms containing $\Delta p$ yield the energy flux
\[
f^e =  \left[ \frac{p_j p_{j+1}}{(\gamma-1) \avg{\rho} \alog{p}} - \frac{1}{2} \avg{u^2} \right] f^\rho + \avg{u} f^m
\]
while the terms containing $\Delta\rho$ yield
\[
\left[ \frac{\gamma}{(\gamma-1)\alog{\rho}} - \frac{1}{2} \avg{r} \avg{u^2} \right] f^\rho + \avg{u} \ \avg{r} f^m - \avg{r} f^e = \avg{u}
\]
From the last two equations and using the relation $\avg{r}=\avg{p}/(p_j p_{j+1})$, we obtain the mass flux as
\[
f^\rho = \frac{\alog{\rho}\avg{u}}{ \frac{\gamma}{(\gamma-1)} - \frac{\avg{p} \alog{\rho}}{(\gamma-1) \avg{\rho} \alog{p}} }
\]
The above mass flux depends on the ratio of specific heats $\gamma$ which is not physically meaningful.

We can also use the variables $p, u, \beta$ as independent variables. In this case we have the following exact linearizations
\begin{eqnarray*}
\Delta v_1 &=& \frac{ \Delta p}{\alog{p}} + \frac{\gamma}{\gamma-1} \frac{\Delta\beta}{\alog{\beta}} - 2 \avg{u} \avg{\beta} \Delta u \\
\Delta v_2 &=& 2(\avg{u}\Delta\beta + \avg{\beta}\Delta u) \\
\Delta v_3 &=& -2\Delta\beta \\
\Delta(\rho u) &=& \avg{\rho}\Delta u + 2\avg{u} \ \avg{p} \Delta\beta + 2 \avg{u} \avg{\beta} \Delta p
\end{eqnarray*}
Then condition~(\ref{eq:tadcond}) gives the fluxes as
\[
f^\rho = 2 \alog{p} \avg{\beta} \avg{u}, \qquad f^m = \frac{ \avg{\rho}}{2\avg{\beta}} + \avg{u} f^\rho, \qquad f^e = \left[ \frac{\gamma}{2(\gamma-1) \alog{\beta}} - \frac{1}{2} \avg{u^2} \right] f^\rho + \avg{u} f^m
\]
While this has the kinetic energy preservation property for the momentum flux, we find that the energy flux is not consistent. This indicates that the variables $\rho$, $u$, $\beta$ are the appropriate set for the satisfaction of the entropy conservative flux condition~(\ref{eq:tadcond}) for the Euler equations.
\subsection{Entropy equation}
Now we can derive the global entropy balance equation for the semi-discrete finite volume scheme for the Navier-Stokes equations. Taking dot product of equation~(\ref{eq:sdfvm}) with the entropy variables $\ent_j$ and summing up over all the finite volumes yields
\begin{eqnarray*}
\sum_j  \ent_j \cdot \dd{\con_j}{t} \Delta x + \sum_j \ent_j \cdot (\fx_\jph - \fx_\jmh) &=& \sum_j \ent_j \cdot (\gx_\jph - \gx_\jmh) \\
\sum_j  \dd{U_j}{t}\Delta x +  \sum_j (\ent_j - \ent_{j+1}) \cdot \fx_\jph &=& \sum_j (\ent_j - \ent_{j+1}) \cdot \gx_\jph \\
\sum_j  \dd{U_j}{t} \Delta x +  \sum_j [F_\jph - F_\jmh + O(\Delta x^3)] &=& -\sum_j \Delta \ent_\jph \cdot \gx_\jph 
\end{eqnarray*}
and the $O(\Delta x)^3$ terms are not present if the entropy conservative flux is used. The terms involving $F_\jph$ represent convection of entropy and cancel one another when we sum over all cells. The terms on the right which consist of viscous shear stress and heat flux can be shown to lead to entropy generation.
\begin{eqnarray*}
-\Delta\ent \cdot \gx &=&  -2 (\avg{\beta} \Delta u + \avg{u} \Delta\beta) \tau + 2 \Delta(\beta) (\tilde{u} \tau - q) \\
&=& - 2 \avg{\beta} \tau \Delta u + 2 (-\avg{u} + \tilde{u}) \tau \Delta\beta - 2 q \Delta\beta \\
&=& - 2 \avg{\beta} \frac{4}{3}\mu\frac{\Delta u}{\Delta x} \Delta u + 0 + 2 \kappa \frac{\Delta T}{\Delta x} \Delta\beta, \qquad \textrm{if} \quad \tilde{u} = \avg{u}\\
&=& - \frac{8\mu\avg{\beta}}{3} \left( \frac{\Delta u}{\Delta x} \right)^2 \Delta x + 2 \kappa \frac{\Delta T}{\Delta x} \frac{(-\Delta T)}{2 R T_j T_{j+1}} \\
&=& - \frac{8\mu\avg{\beta}}{3} \left( \frac{\Delta u}{\Delta x} \right)^2 \Delta x -  \frac{\kappa}{R T_j T_{j+1}} \left(\frac{\Delta T}{\Delta x}\right)^2 \Delta x \le 0
\end{eqnarray*}
Hence the entropy equation becomes
\[
\sum_j \Delta x \dd{U_j}{t} +  \sum_j O(\Delta x^3) = - \sum_j \left[ \frac{8\mu\avg{\beta}_\jph}{3} \left( \frac{\Delta u_\jph}{\Delta x} \right)^2  +  \frac{\kappa}{R T_j T_{j+1}} \left(\frac{\Delta T_\jph}{\Delta x}\right)^2 \right] \Delta x \le 0
\]
If the entropy conservative flux is used, the $O(\Delta x)^3$ terms are not present and the entropy condition is satisfied exactly for any mesh size $\Delta x$, which is consistent with the entropy condition from the second law of thermodynamics.
\subsection{Summary of flux formulae}
We now list the approximately and exactly entropy consistent centered fluxes for the 1-D Euler equations. \\

\noindent
(1) The centered, kinetic energy preserving and {\em entropy consistent} numerical flux which will be denoted by $\fx^*$ is given by
\begin{eqnarray}
f^{*,\rho}_\jph &=& \avg{\rho}_\jph \avg{u}_\jph \\
f^{*,m}_\jph &=& \tp_\jph + \avg{u}_\jph f^{*,\rho}_\jph \\
f^{*,e}_\jph &=& \left[ \frac{1}{2(\gamma-1) \avg{\beta}_\jph} - \frac{1}{2} \avg{u^2}_\jph \right] f^{*,\rho}_\jph + \avg{u}_\jph f^{*,m}_\jph
\end{eqnarray}
where
\begin{equation}
\tp_\jph = \frac{\avg{\rho}_\jph}{2 \avg{\beta}_\jph}
\label{eq:tp}
\end{equation}
This is equivalent to using the harmonic average for the temperature, i.e.,
\[
\tp_\jph = R \avg{\rho}_\jph \hat{T}_\jph, \qquad \hat{T}_\jph = \frac{2 T_j T_{j+1}}{T_j + T_{j+1}}
\]
(2) The centered, kinetic energy preserving and {\em entropy conservative} numerical flux is given by
\begin{eqnarray}
f^{*,\rho}_\jph &=& \alog{\rho}_\jph \avg{u}_\jph \\
f^{*,m}_\jph &=& \tp_\jph + \avg{u}_\jph f^{*,\rho}_\jph \\
f^{*,e}_\jph &=& \left[ \frac{1}{2(\gamma-1) \alog{\beta}_\jph} - \frac{1}{2} \avg{u^2}_\jph \right] f^{*,\rho}_\jph + \avg{u}_\jph f^{*,m}_\jph
\end{eqnarray}
where $\tp_\jph$ is given by equation~(\ref{eq:tp}).
\section{Scalar artificial dissipation flux}
\label{sec:scal}
The entropy consistent fluxes have a central character and will lead to a stable scheme for Navier-Stokes equations only on highly resolved meshes. For the scheme to be stable on coarse meshes and inviscid problems, let us introduce second order, scalar artificial dissipation term into the inviscid fluxes, so that the numerical flux becomes
\[
\fx_\jph = \fx_\jph^* - \frac{1}{2} \lambda_\jph \Dx_\jph, \qquad \Dx=[D^\rho, \ D^m, \ D^e]^\top, \qquad \lambda \ge 0
\]
We will derive the form of the dissipation from kinetic energy and entropy stability consoderations.
\subsection{Kinetic energy stability}
We make the following obvious choice for the dissipation in mass and momentum fluxes
\[
D^\rho = \Delta\rho, \qquad D^m = \Delta(\rho u)
\]
The kinetic energy equation becomes
\begin{eqnarray*}
\sum_j  \dd{K_j}{t} \Delta x
&=& - \frac{1}{2} \sum_j \lambda_\jph \left[ (u_j - u_{j+1})\left( \frac{1}{2}(u_j + u_{j+1}) D^\rho_\jph - D^m_\jph \right) \right] \\
&& + \sum_j \left[ \frac{\Delta u_\jph}{\Delta x} \tp_\jph  - \frac{4}{3}\mu \left(\frac{\Delta u_\jph}{\Delta x} \right)^2 \right] \Delta x \\
&=& - \frac{1}{2} \sum_j \lambda_\jph \avg{\rho}_\jph (\Delta u_\jph)^2 + \sum_j \left[ \frac{\Delta u_\jph}{\Delta x} \tp_\jph  - \frac{4}{3}\mu \left(\frac{\Delta u_\jph}{\Delta x} \right)^2 \right] \Delta x
\end{eqnarray*}
The artificial dissipation terms in the flux lead to dissipation of kinetic energy since the first term on the right is negative while the remaining terms are consistent with the continuous equation~(\ref{eq:kecont}). Thus stability in the sense that the kinetic energy does not grow spuriously is maintained by the choice of the scalar dissipation flux.
\subsection{Entropy condition}
The dissipation in the energy flux will be determined by satisfying entropy condition approximately or exactly as before. The natural choice is to take $D^e=\Delta E$ but this does not allow us to show entropy stability. Instead, the dissipation in the energy equation is taken to be of the form
\[
D^e = \left\{ \frac{1}{2(\gamma-1)\avg{\beta}} + \frac{1}{2}u_j u_{j+1} \right\} \Delta\rho + \avg{\rho} \ \avg{u} \Delta u + \frac{\avg{\rho}}{2(\gamma-1)} \Delta(1/\beta)
\]
The entropy equation including the dissipation term is
\[
\begin{aligned}
\sum_j  \dd{U_j}{t} \Delta x +  \sum_j O(\Delta x^3) = &- \sum_j \left[ \frac{8\mu\avg{\beta}_\jph}{3} \left( \frac{\Delta u_\jph}{\Delta x} \right)^2  +  \frac{\kappa}{R T_j T_{j+1}} \left(\frac{\Delta T_\jph}{\Delta x}\right)^2 \right] \Delta x \\
 &- \frac{1}{2} \sum_j \lambda_\jph \Delta \ent_\jph \cdot \Dx_\jph
\end{aligned}
\]
The first term on the right is consistent with the continuous entropy equation while the second term is due to the dissipative flux. The contribution from the dissipation terms can be calculated as follows.
\begin{eqnarray*}
\Delta \ent \cdot \Dx &=& + \left[ \frac{\Delta\rho}{\avg{\rho}} + \frac{1}{\gamma-1} \frac{\Delta\beta}{\avg{\beta}} + O(\Delta x^3) - 2 \avg{\beta} \avg{u} \Delta u - \avg{u^2} \Delta\beta \right] \Delta\rho  \\
&& + 2 ( \avg{\beta} \Delta u + \avg{u} \Delta \beta ) (\avg{\rho} \Delta u + \avg{u} \Delta\rho)  \\
&& - 2 \Delta(\beta) \left[ \left\{ \frac{1}{2(\gamma-1)\avg{\beta}} + \frac{1}{2}u_j u_{j+1} \right\} \Delta\rho + \avg{\rho} \ \avg{u} \Delta u +  \frac{\avg{\rho}}{2(\gamma-1)} \Delta(1/\beta) \right] \\
&=& \underbrace{\frac{(\Delta\rho)^2}{\avg{\rho}} + 2 \avg{\rho} \avg{\beta} (\Delta u)^2 + \frac{\avg{\rho}}{(\gamma-1)\beta_j \beta_{j+1}} (\Delta\beta)^2}_{\ge 0} + O(\Delta x^4)
\end{eqnarray*}
The leading order terms in the last equation are positive and hence lead to entropy generation; the scheme satisfies the entropy condition in the limit of $\Delta x \to 0$.

The scalar dissipation terms can also be constructed to be exactly entropy dissipative which can be used in combination with the entropy conservative flux. The dissipation in the mass and momentum fluxes are chosen as before which leads to the same kinetic energy stability as above, while the dissipation in the energy flux is taken to be
\[
D^e = \left\{ \frac{1}{2(\gamma-1)\alog{\beta}} + \frac{1}{2}u_j u_{j+1} \right\} \Delta\rho + \avg{\rho} \ \avg{u} \Delta u + \frac{\avg{\rho}}{2(\gamma-1)} \Delta(1/\beta)
\]
The only difference from the previous case is that the logarithmic average of $\beta$ has been used instead of the arithmetic average. Then the entropy production from the dissipative terms is 
\[
\Delta \ent \cdot \Dx = \frac{(\Delta\rho)^2}{\alog{\rho}} + 2 \avg{\rho} \avg{\beta} (\Delta u)^2 + \frac{\avg{\rho}}{(\gamma-1)\beta_j \beta_{j+1}} (\Delta\beta)^2 \ge 0
\]
and is non-negative which satisfies the entropy condition. We note that the dissipation is triggered if any one of the state variables $\rho, u, \beta$ is non-uniform. The physical dissipation on the other hand acts only in the presence of velocity and/or temperature gradients, and is unaffected by density gradients.
\subsection{Summary of scalar dissipation flux}
The second order kinetic energy preserving and approximately entropy consistent scalar dissipation terms are given by
\begin{eqnarray*}
D^\rho_\jph &=& \Delta\rho_\jph \\
D^m_\jph &=& \Delta(\rho u)_\jph = \avg{u}_\jph \Delta \rho_\jph + \avg{\rho}_\jph \Delta u_\jph\\
D^e_\jph &=& \left\{ \frac{1}{2(\gamma-1)\avg{\beta}_\jph} + \frac{1}{2}u_j u_{j+1} \right\} \Delta\rho_\jph + \avg{\rho}_\jph \avg{u}_\jph \Delta u_\jph + \frac{\avg{\rho}_\jph}{2(\gamma-1)} \Delta (1/\beta)_\jph
\end{eqnarray*}
and
\begin{equation}
\lambda_\jph = |\avg{u}_\jph| + \sqrt{\frac{\gamma}{2\avg{\beta}_\jph}}
\end{equation}
is the maximum wave speed at the interface. If we replace the arithmetic average $\avg{\beta}_\jph$ with the logarithmic average $\alog{\beta}_\jph$ in the above equations, then we obtain the exactly entropy consistent dissipation terms. Scalar artificial dissipation reduces the scheme to first order accuracy everywhere, so we can tune the dissipation to switch on only near shocks. Pressure gradients are a good indicator of the presence of shocks. However, lack of any dissipation in smooth regions can lead to instabilities. Hence, as in the JST scheme~\cite{jst}, we will add a blend of second order and fourth order dissipation by replacing the difference terms in the above equation, e.g.,
\begin{equation}
\Delta\rho_\jph \to \ept_\jph (\rho_{j+1} - \rho_j) - \epf_\jph (\rho_{j+2} - 3 \rho_{j+1} + 3 \rho_j - \rho_{j-1})
\end{equation}
where $\ept$ and $\epf$ are adapted to the flow, with similar expressions for $\Delta u$ and $\Delta T$. Define
\[
\nu_j = \frac{|p_{j-1} - 2p_j + p_{j+1}|}{|p_{j-1} + 2 p_j + p_{j+1}|}, \qquad \nu_\jph = \max(\nu_j, \nu_{j+1})
\]
\[
\ept_\jph = \min(1, \kapt \nu_\jph), \qquad \epf_\jph = \max(0, \kapf - \ept_\jph)
\]
with $\kapt \ge 0$ and $\kapf \ge 0$. In smooth regions of flow, $\ept = O(\Delta x^2)$, $\epf = O(1)$ and $\Dx=O(\Delta x^3)$, while near shocks $\ept=O(1)$, $\epf=0$ and $\Dx=O(\Delta x)$. Thus the scheme is second order accurate in smooth regions of the flow and if the mesh is well resolved, the scheme is second order accurate everywhere. For a very well resolved flow, the fourth order dissipation should be sufficient to stabilize the scheme as we already have consistency of kinetic energy and entropy condition which imparts some non-linear stability to the numerical scheme.
\section{Matrix dissipation flux}
\label{sec:mat}
Scalar dissipation schemes perform satisfactorily for weak shocks. For problems with strong shocks, one needs to add greater amount of scalar dissipation which unfortunately smears the contact discontinuities. An alternate way to add dissipation is via a matrix operator which leads to better shock capturing schemes. The classic example is the Roe scheme~\cite{Roe1981357} which is based on a linearization of the non-linear conservation law about some average state. The numerical flux of the Roe scheme has the form
\[
\fx_\jph = \frac{1}{2}(\fx_j + \fx_{j+1}) - \frac{1}{2} R_\jph |\Lambda_\jph| R^{-1}_\jph \Delta \con_\jph
\]
where $R$ is the matrix of right eigenvectors and $\Lambda$ is the diagonal matrix containing the eigenvalues, both of which are evaluated at the Roe average state. 
\[
R = \begin{bmatrix}
1 & 1  & 1 \\
u-a & u & u+a \\
H - ua & \frac{1}{2}u^2 & H+ua
\end{bmatrix}, \qquad |\Lambda| = |\Lambda|^{Roe} = \textrm{diag}\left[ |u-a|, \ |u|, \ |u+a| \right]
\]
The dissipative flux can also be written in terms of the jumps in the entropy variables. By a linearization we can write $\Delta \con = \con_\ent \Delta \ent$ and due to a theorem of Barth~\cite{barth1998}, there exists a scaling of the eigenvectors $R \to \tilde{R}$ such that $\con_\ent = \tilde{R} \tilde{R}^\top$. Then a Roe-type flux can be written by using jumps in the entropy variables instead of the conserved variables as
\[
\fx_\jph = \frac{1}{2}(\fx_j + \fx_{j+1}) - \frac{1}{2} \tilde{R}_\jph |\Lambda_\jph| \tilde{R}^\top_\jph \Delta \ent_\jph
\]
If we want to use the eigenvectors in the original form $R$, then the flux can be written as
\begin{equation}
\fx_\jph = \frac{1}{2}(\fx_j + \fx_{j+1}) - \frac{1}{2} R_\jph |\Lambda_\jph| S_\jph R^\top_\jph \Delta \ent_\jph
\label{eq:roeentflux}
\end{equation}
where the matrix $S$ provides the appropriate scaling and is given by
\[
S = \textrm{diag}\left[ \frac{\rho}{2\gamma}, \ \frac{(\gamma-1)\rho}{\gamma}, \ \frac{\rho}{2\gamma} \right]
\]
In fact the matrix $S$ is chosen so that $R^{-1}\ud\con = SR^\top\ud\ent$ which also motivates the form of the flux given by equation~(\ref{eq:roeentflux}). Following this approach we can add a matrix dissipation flux to our kinetic energy and entropy consistent/conservative centered flux $\fx^*$ to obtain the following dissipative numerical flux function
\begin{equation}
\fx_\jph = \fx^*_\jph - \frac{1}{2} R_\jph |\Lambda_\jph| S_\jph R^\top_\jph \Delta \ent_\jph
\end{equation}
The above numerical flux together with $|\Lambda|=|\Lambda|^{Roe}$ will be called the KEP-ES flux; when the approximately entropy consistent flux is used for the central part, it will be explicitly qualified as (AC). We note that the dissipation matrix $Q = R |\Lambda| S R^\top$ is positive definite which allows us to derive the entropy inequality as follows. By taking the dot product of the semi-discrete finite volume scheme for Euler equations with $\ent_j$ we obtain the entropy equation
\[
\Delta x \dd{U_j}{t} + \ent_j \cdot (\fx_\jph^* - \fx_\jmh^*) = \frac{1}{2}\left[ \ent_j^\top Q_\jph \Delta\ent_\jph - \ent_j^\top Q_\jmh \Delta\ent_\jmh \right]
\]
Define the numerical entropy flux to be
\[
F_\jph = \avg{\ent}_\jph \cdot \fx_\jph^* - \avg{\psi}_\jph + \frac{1}{2} \bar{\ent}_\jph^\top Q_\jph \Delta\ent_\jph
\]
which is a consistent flux, then the cell entropy equation can be written as
\[
\Delta x \dd{U_j}{t} + F_\jph - F_\jmh + O(\Delta x^3) = -\frac{1}{4}\left[ \Delta\ent_\jph^\top Q_\jph \Delta\ent_\jph + \Delta\ent_\jmh^\top Q_\jmh \Delta\ent_\jmh \right] \le 0
\]
which shows that the dissipation terms in the flux lead to entropy generation. The $O(\Delta x)^3$ terms are not present if we use the exactly entropy conservative flux for the centered flux $\fx^*$ and we then obtain strict entropy stability for any mesh size $\Delta x$.
\subsection{Resolution of stationary contact discontinuity}
Contact discontinuities occur because of the presence of linearly degenerate eigenvector fields. Because of their linear nature, the accurate resolution of contact waves is more difficult than shocks which have an inherent steepening mechanism due to their non-linear nature. The original Roe scheme exactly preserves stationary contact waves. This is a desirable property since it affects the accuracy with which shear layers and boundary layers are computed. The new flux functions with entropy variable based dissipation will be able to resolve stationary contacts exactly if we evaluate the enthalpy $H$ occuring in the eigenvector matrix $R$ in an appropriate manner. Consider an initial condition such that
\[
\rho_k = \begin{cases}
\rho_l & k \le j \\
\rho_r & k > j
\end{cases}, \qquad u_k = 0, \qquad p_k = p
\]
This initial discontinuity satisfies the RH jump conditions with zero speed and hence it represents a stationary solution. The numerical scheme will preserve this solution exactly if
\[
\ldots = \fx_{j-1} = \fx_j = \fx_\jph = \fx_{j+\frac{3}{2}} = \fx_{j+1} = \ldots \qquad\Longrightarrow\qquad \fx_\jph = \begin{bmatrix}
0 \\
p \\
0
\end{bmatrix}
\]
For the above initial condition, the centered flux $\fx^*$ already satisfies the above condition, i.e., $\fx^*_\jph = [0, \ p, \ 0]^\top$. Hence we need to have $Q_\jph \Delta \ent_\jph = 0$ in order to preserve the contact wave. The dissipation matrix $Q$ requires some average values of $u$, $a$, $H$ and $\rho$ which will be denoted by $u_\jph$, etc. For the above initial condition, any consistent averaging would yield $u_\jph=0$. We will compute the enthalpy as $H_\jph = a_\jph^2/(\gamma-1) + u_\jph^2/2$. The dissipation matrix then takes the form
\[
Q_\jph = \frac{1}{\gamma} \rho_\jph |a_\jph|
\begin{bmatrix}
1 & 0 & H_\jph \\
0 & a_\jph^2 & 0 \\
H_\jph & 0 & H_\jph^2
\end{bmatrix}, \quad \Delta \ent_\jph = \begin{bmatrix}
-\frac{1}{\gamma-1} \Delta s_\jph \\
0 \\
-2\Delta\beta_\jph
\end{bmatrix}, \quad H_\jph = \frac{a_\jph^2}{\gamma-1}
\]
so that
\[
Q_\jph \Delta \ent_\jph = -\frac{1}{\gamma} \rho_\jph |a_\jph| \left( \frac{1}{\gamma-1} \Delta s_\jph + 2 H_\jph \Delta\beta_\jph \right)
\begin{bmatrix}
1 \\
0 \\
H_\jph
\end{bmatrix}
\]
This is zero provided the term inside the brackets is zero. We will show that this is the case if we choose 
\begin{equation}
a_\jph = \sqrt{ \frac{\gamma}{2 \alog{\beta}_\jph} }
\label{eq:soundspeed}
\end{equation}
Then
\begin{eqnarray*}
\frac{1}{\gamma-1} \Delta s_\jph + 2 H_\jph \Delta\beta_\jph &=& -\frac{\Delta\rho_\jph}{\alog{\rho}_\jph} - \frac{1}{\gamma-1} \frac{\Delta\beta_\jph}{\alog{\beta}_\jph} + \frac{2}{\gamma-1} a_\jph^2 \Delta\beta_\jph \\
&=& -\frac{\Delta\rho_\jph}{\alog{\rho}_\jph} + \frac{\Delta\beta_\jph}{\alog{\beta}_\jph} = 0
\end{eqnarray*}
where the last equality follows due to the chosen initial condition. Thus exact resolution of stationary contacts is possible provided we evaluate the sound speed $a_\jph$ using the logarithmic average of $\beta$ as in equation~(\ref{eq:soundspeed}). The other quantities $\rho_\jph$, $u_\jph$ can be evaluated using arithmetic or logarithmic averaging. Since the velocity already appears in the central flux $\fx^*$ as an arithmetic average we can use the same choice in the dissipation matrix also. The density appears in the approximately consistent flux as an arithmetic average and in the exactly conservative flux as the logarithmic average; we will use the corresponding approximations in the dissipation matrix which avoids additional computations.
\subsection{Kinetic energy stability}
Consider the finite volume scheme with entropy consistent/conservative flux together with entropy variables based matrix dissipation. The semi-discrete kinetic energy equation for this scheme is
\begin{eqnarray*}
\sum_j  \dd{K_j}{t} \Delta x
&=& - \frac{1}{2} \sum_j (u_j - u_{j+1})\left( \frac{1}{2}(u_j + u_{j+1}) (Q\Delta\ent)^\rho_\jph - (Q\Delta\ent)^m_\jph \right) \\
&& + \sum_j \left[ \frac{\Delta u_\jph}{\Delta x} \tp_\jph  - \frac{4}{3}\mu \left(\frac{\Delta u_\jph}{\Delta x} \right)^2 \right] \Delta x \\
&=&  \frac{1}{2} \sum_j \Delta u_\jph \left[\avg{u}_\jph, \ -1, \ 0\right] Q_\jph \Delta\ent_\jph 
 + \sum_j \left[ \frac{\Delta u_\jph}{\Delta x} \tp_\jph  - \frac{4}{3}\mu \left(\frac{\Delta u_\jph}{\Delta x} \right)^2 \right] \Delta x
\end{eqnarray*}
The second term is consistent with the exact kinetic energy equation while the first term is due to the additional dissipation in the flux. This term will be dissipative for the kinetic energy provided
\begin{equation}
\left[\avg{u}_\jph, \ -1, \ 0\right] Q_\jph \Delta\ent_\jph = -\alpha_\jph \Delta u_\jph, \qquad \alpha_\jph \ge 0
\label{eq:kediss}
\end{equation}
Now dropping subscripts and overbars, we have
\[
\begin{aligned}
\left[u, \ -1, \ 0\right] Q \Delta\ent = \frac{a\rho}{2\gamma} [ (|\lambda_1| -|\lambda_3|)\Delta v_1 & - (a(|\lambda_1| + |\lambda_3|) + u(-|\lambda_1| + |\lambda_3|)) \Delta v_2 \\
&  - (au(|\lambda_1| + |\lambda_3|) + H(-|\lambda_1| + |\lambda_3|)) \Delta v_3 ]
\end{aligned}
\]
Note that since $\Delta v_1$ contains $\Delta\rho$ whereas $\Delta v_2$, $\Delta v_3$ do not, we can satisfy condition~(\ref{eq:kediss}) only if $|\lambda_1| = |\lambda_3| = \lambda$ for some $\lambda \ge 0$. Then the kinetic energy equation becomes
\[
\sum_j  \dd{K_j}{t} \Delta x =  -\frac{1}{2\gamma} \sum_j a^2_\jph \rho_\jph \avg{\beta}_\jph \lambda_\jph (\Delta u_\jph)^2 + \sum_j \left[ \frac{\Delta u_\jph}{\Delta x} \tp_\jph  - \frac{4}{3}\mu \left(\frac{\Delta u_\jph}{\Delta x} \right)^2 \right] \Delta x
\]
The first term on the right is dissipative and hence the scheme is stable for the kinetic energy. The Roe-type dissipation does not satisfy the required condition since $\lambda_1=u-a$ and $\lambda_3=u+a$. The condition that the first and third eigenvalues should be equal is satisfied if we choose the Rusanov form of the dissipation which corresponds to
\begin{equation}
|\Lambda| = |\Lambda|^{Rus} = \lambda I, \qquad \lambda = |u| + a
\label{eq:rus}
\end{equation}
However this scheme is very dissipative and leads to excessive smearing of shocks, contact discontinuities and shear layers. Since only the first and third eigenvalues have to be equal to have dissipation of kinetic energy, another obvious choice is
\begin{equation}
|\Lambda| = |\Lambda|^{KES} = \textrm{diag} \left[ \ \lambda, \ |u|, \ \lambda \ \right], \qquad \lambda = \max( |u-a|, |u+a| ) = |u| + a
\label{eq:kescontact}
\end{equation}
which resolves stationary contacts exactly. However steady shocks are smeared over many cells due to excessive dissipation. The robustness of the scheme due to its kinetic energy and entropy stability makes it attractive when used in combination with high resolution schemes like WENO or discontinuous Galerkin methods. The new entropy conservative flux toegther with the above forms of dissipation will be refered to as KEP-ES(Rus) and KEP-ES(KES) respectively.
\section{Monotone resolution of shocks}
\label{sec:monotone}
\begin{figure}
\begin{center}
Mach = 1.5 \\
\includegraphics[width=0.45\textwidth]{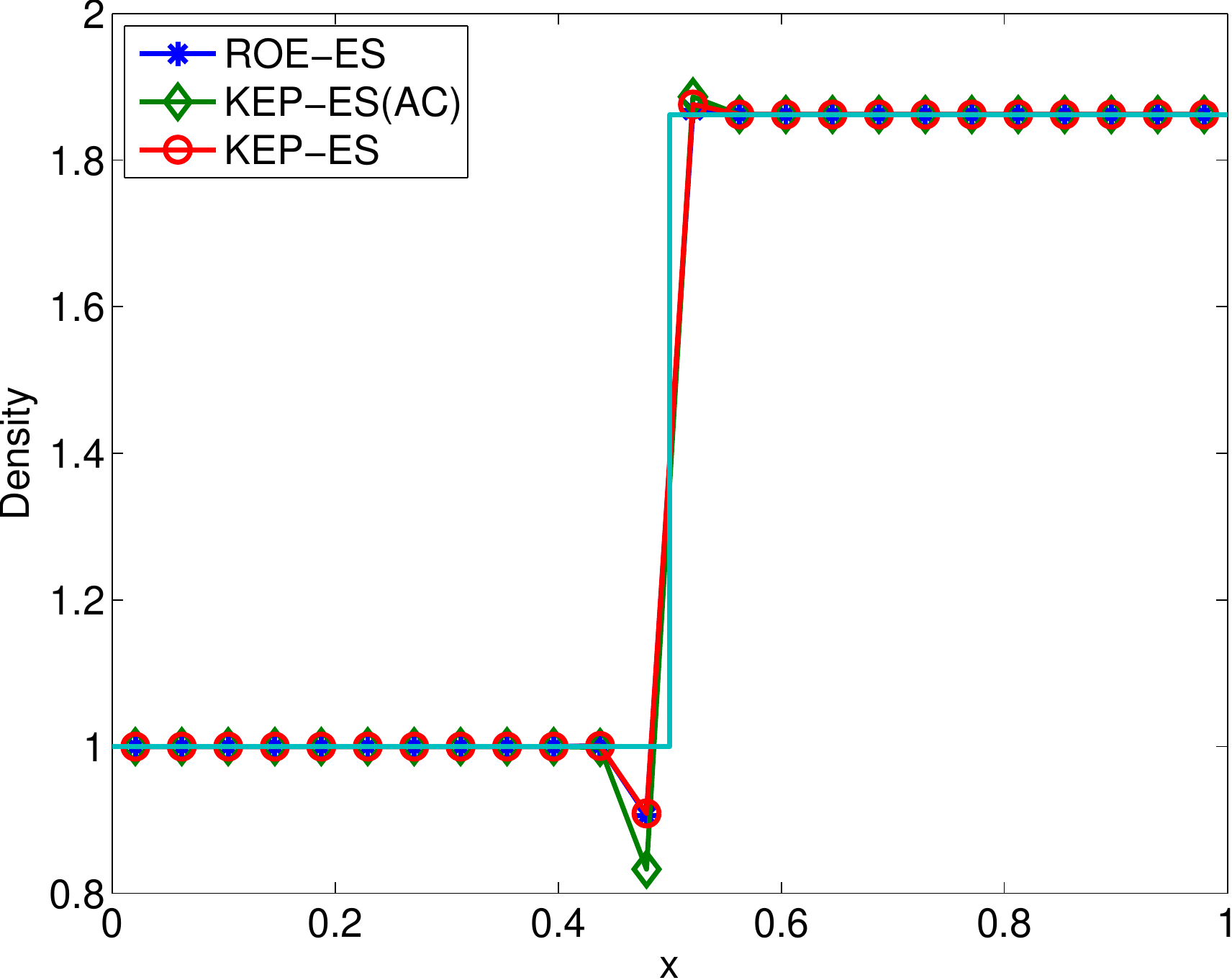}
\includegraphics[width=0.45\textwidth]{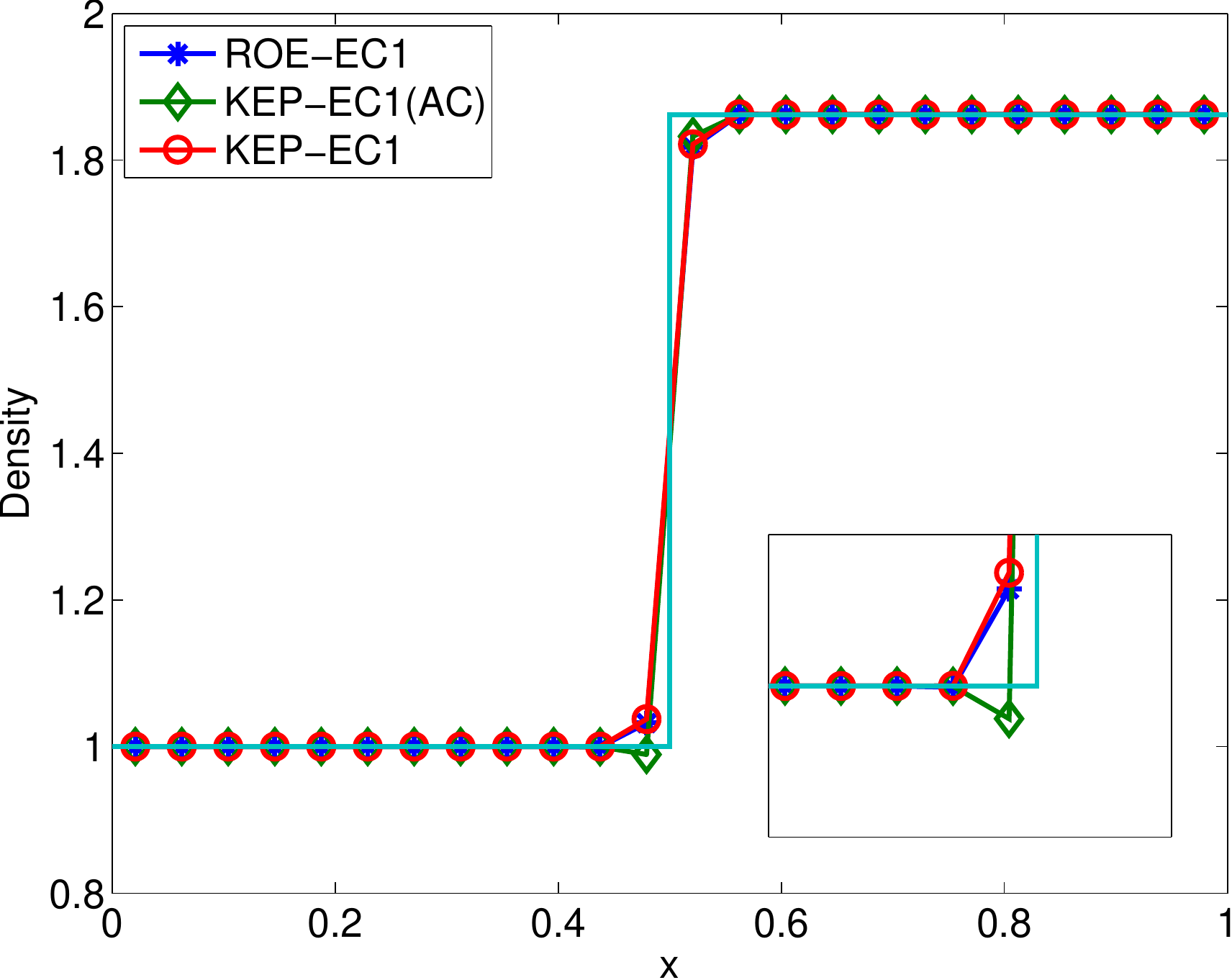} \\
Mach = 4 \\
\includegraphics[width=0.45\textwidth]{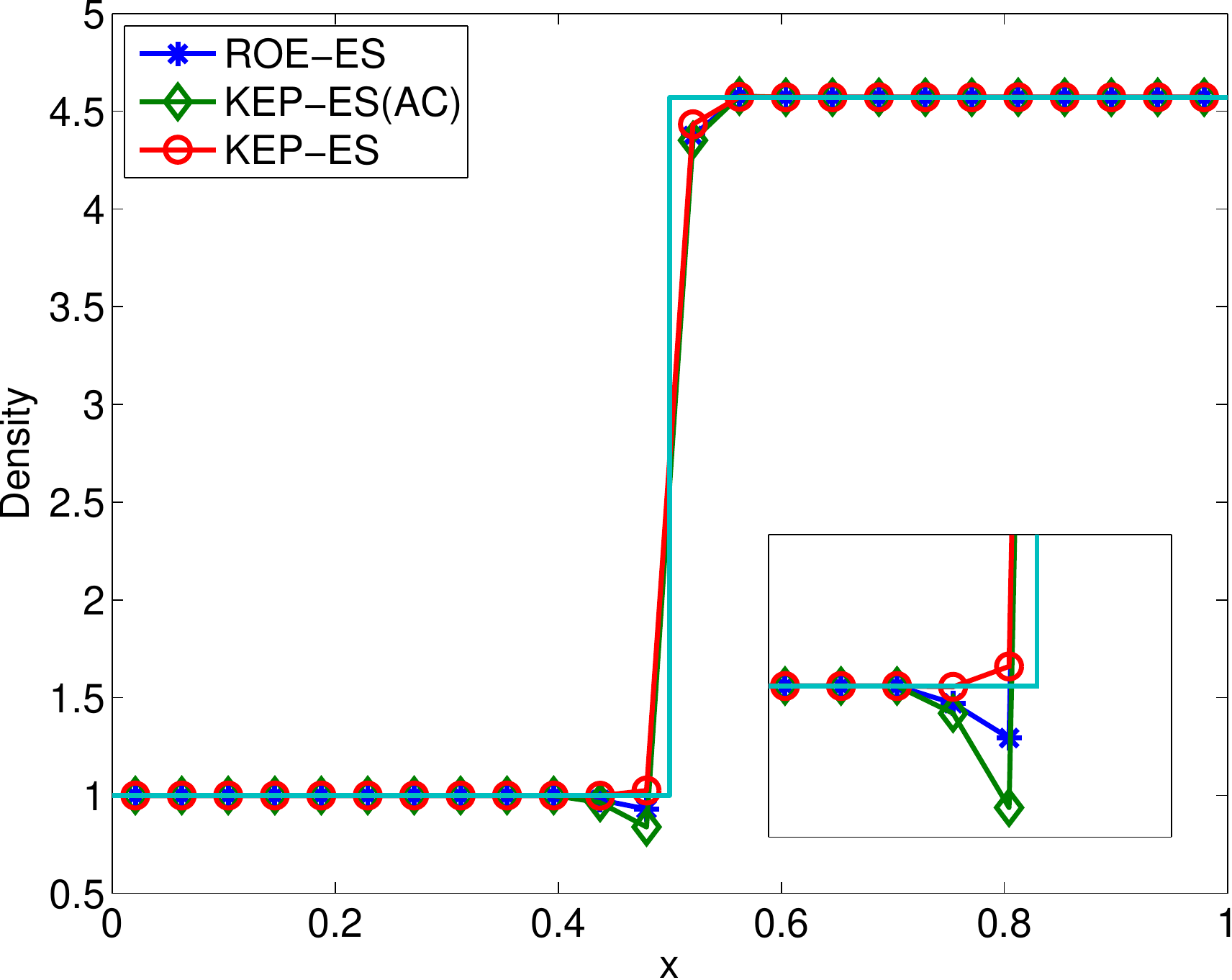}
\includegraphics[width=0.45\textwidth]{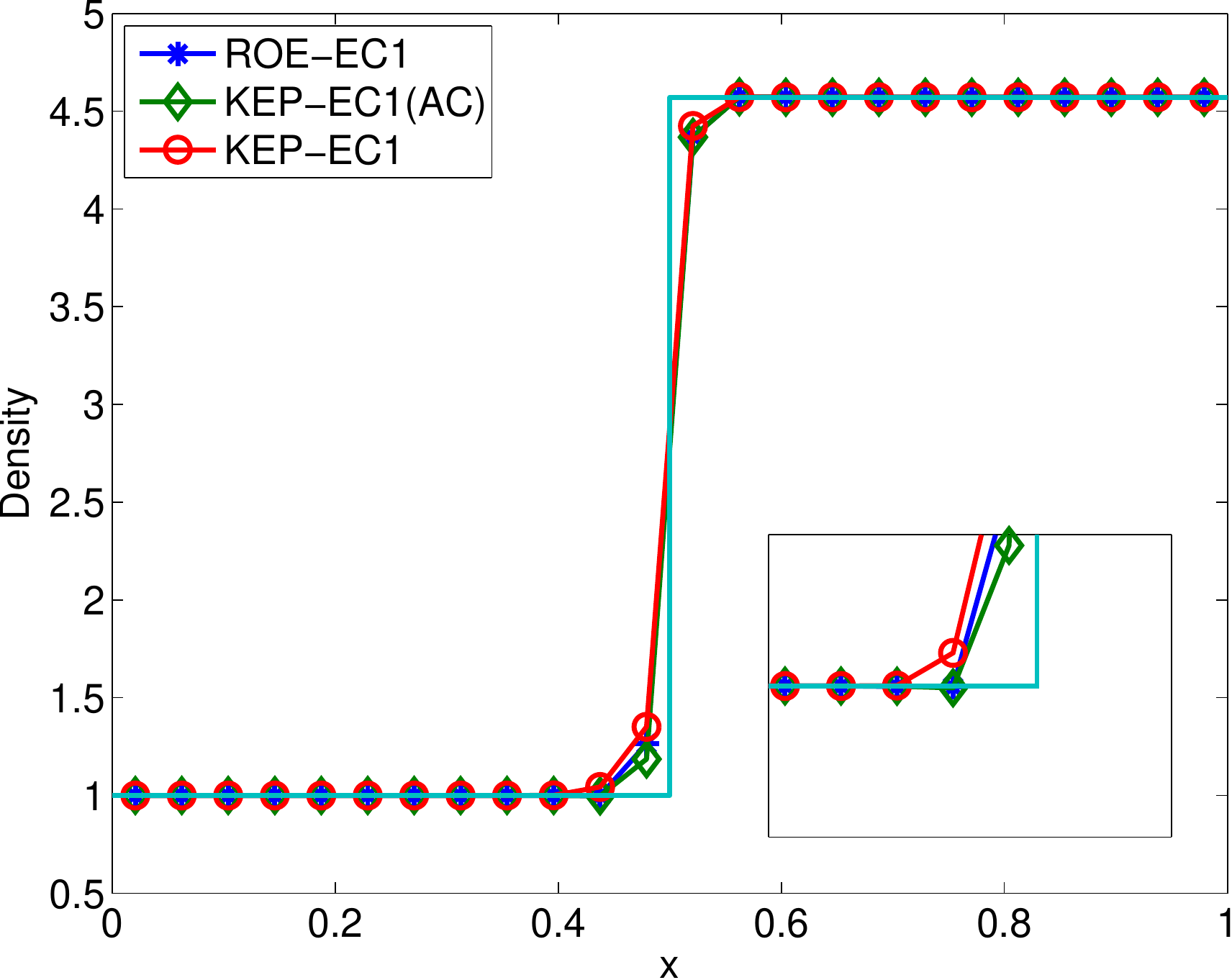} \\
Mach = 20 \\
\includegraphics[width=0.45\textwidth]{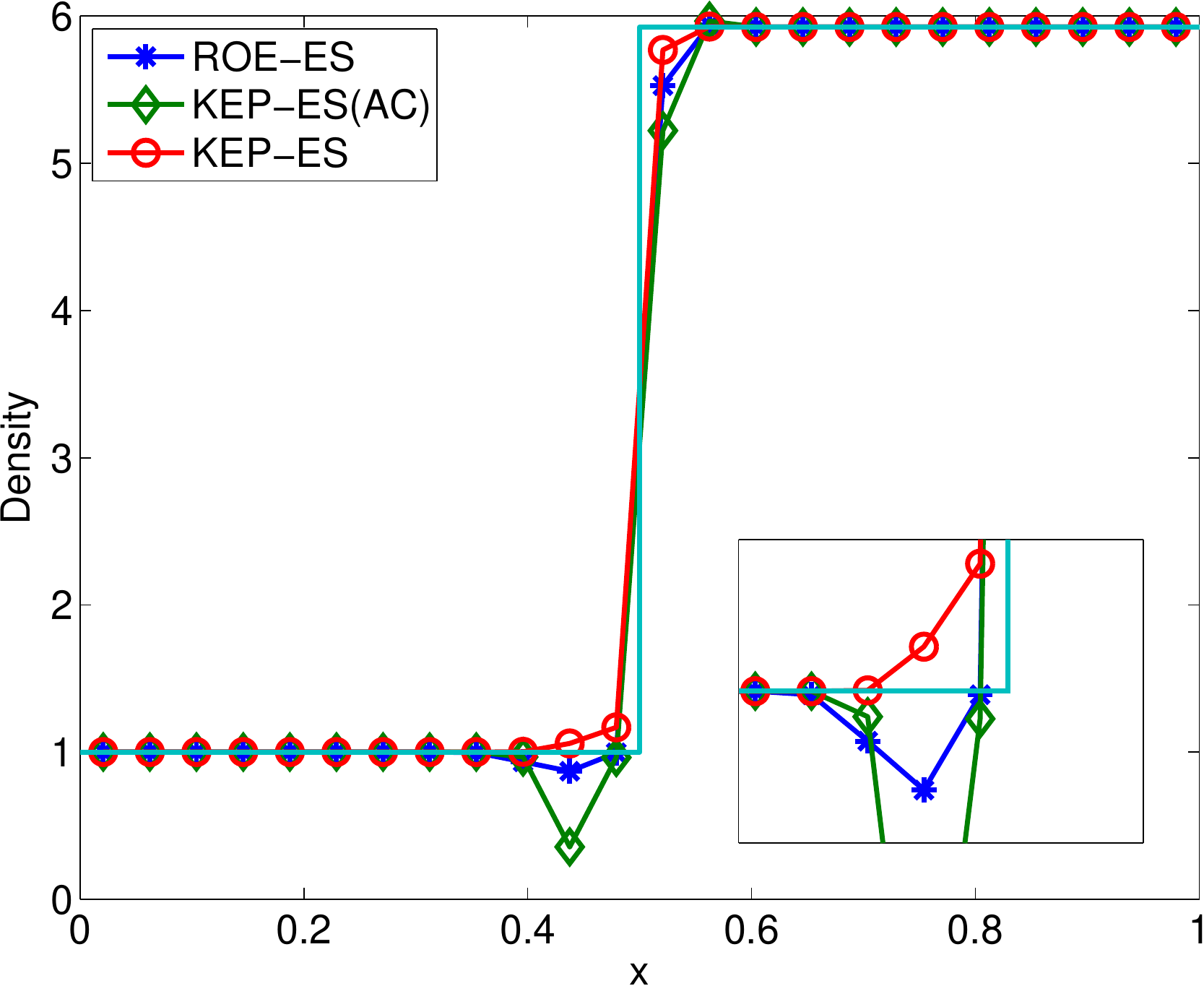}
\includegraphics[width=0.45\textwidth]{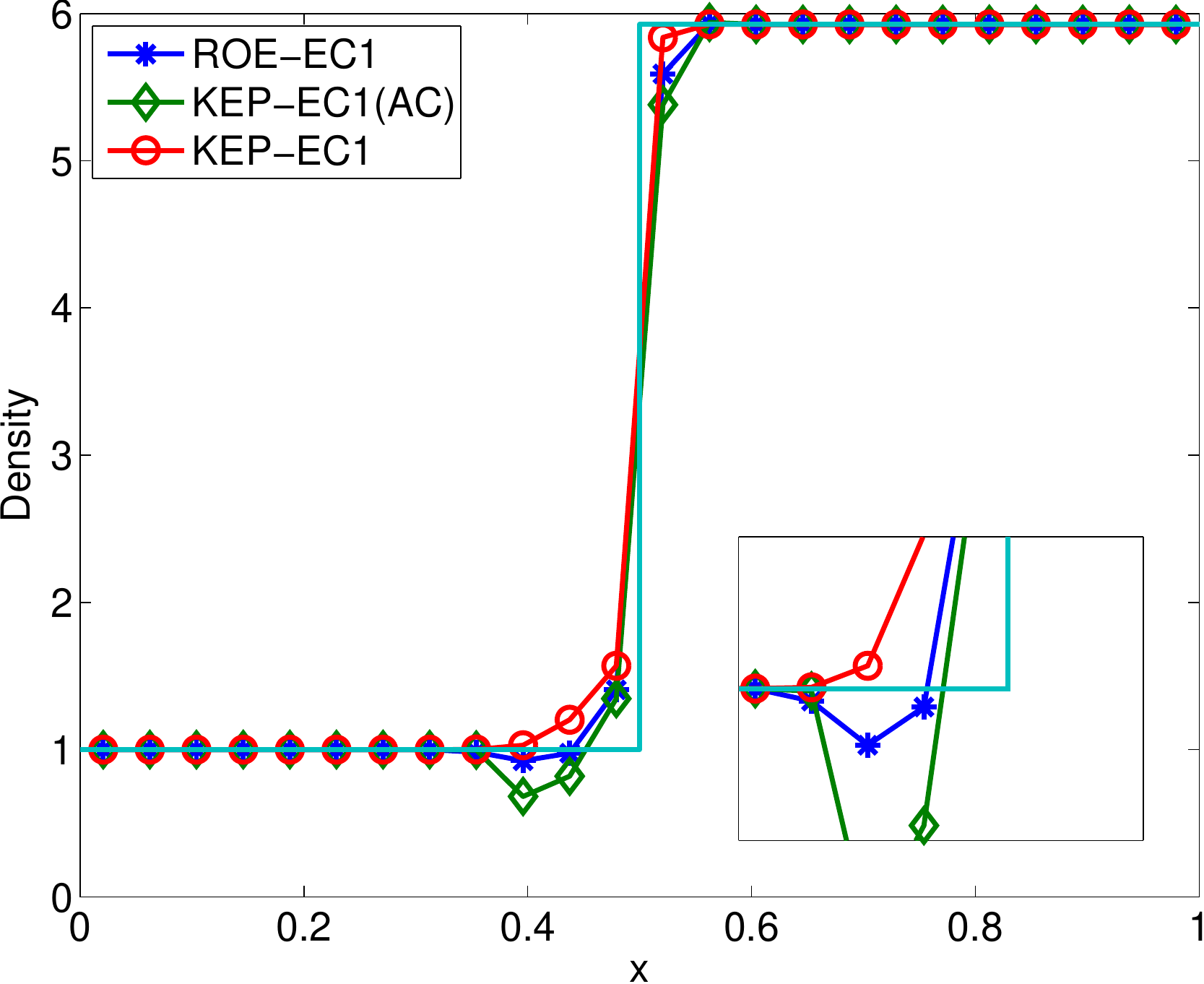} 
\caption{Stationary shock problem: $N=24$ cells}
\label{fig:stshock}
\end{center}
\end{figure}

An important desirable property of a numerical scheme for hyperbolic problems is that it should yield monotone discontinuous solutions. Central schemes yield highly oscillatory solutions near shocks and other discontinuities which are physically incorrect and can even lead to negative density and pressure. Upwind schemes which are based on wave propagation ideas and central schemes with explicitly added dissipation are able to yield non-oscillatory solutions. For scalar conservation laws, the TVD condition provides a complete solution for the design of accurate, non-oscillatory schemes. However the TVD property does not carry over to non-linear systems like the Euler equations and there is a lack of theoretical basis for the design of non-oscillatory schemes for systems of conservation laws. It is hoped that the entropy condition can be a  useful criterion for the design of non-oscillatory schemes. However this condition only specifies that there must be non-zero dissipation which leads to some entropy production but does not say how much entropy generation is necessary to yield monotone solutions. The matrix dissipation flux discussed in the previous sections may not yield monotone shock solutions even though they are entropy stable. In order to study this numerically, we follow Ismail and Roe~\cite{Ismail:2009:AEE:1552569.1552681} and consider a problem with a stationary shock solution. This corresponds to a Riemann problem with the following left and right states
\[
\rho_l = 1, \qquad u_l = 1, \qquad p_l = \frac{1}{\gamma M_l^2}
\]
and
\[
\rho_r = \left[ \frac{2}{(\gamma+1) M_l^2} + \frac{\gamma-1}{\gamma+1} \right]^{-1}, \qquad u_r = \frac{1}{\rho_r}, \qquad p_r = p_l \left[ \frac{2\gamma M_l^2}{\gamma + 1} - \frac{\gamma-1}{\gamma+1} \right]
\]
The Mach number on the left $M_l$ can be changed to produce shocks of different strengths and following~\cite{Ismail:2009:AEE:1552569.1552681} we consider $M_l=1.5, \ 4$ and 20. All computations are performed using a CFL=0.1 and a three stage Runge-Kutta scheme. The density obtained from the different schemes is shown on the left of figure~(\ref{fig:stshock}). The ROE-ES flux and the new approximately entropy consistent flux gives pre-shock oscillations in the density at all Mach numbers while the new entropy conservative flux yields oscillations at the lower Mach number but is monotone at higher Mach numbers. Ismail and Roe attribute the non-monotonicity to insufficient entropy production and hence insufficient dissipation at the shock. For a weak shock, the entropy production is $O(\Delta\rho)^3$ whereas the entropy production due to the dissipative flux is $O(\Delta\rho)^2$. Hence they propose modifying the acoustic eigenvalues so that the eigenvalue matrix $|\Lambda|$ becomes
\begin{equation}
|\Lambda| = |\Lambda|^{EC1} = \textrm{diag}\left[ |u-a| + \beta|\Delta\lambda_1|, \ |u|, \ |u+a| + \beta|\Delta\lambda_3| \right]
\label{eq:lec1}
\end{equation}
where $\Delta\lambda_1$ and $\Delta\lambda_3$ are the jumps in the corresponding eigenvalues across the cell face. If $\beta = \frac{1}{6}$ then the entropy production across a weak shock due to matrix dissipation corresponds to the correct entropy production of $O(\Delta\rho)^3$ as shown in~\cite{Ismail:2009:AEE:1552569.1552681}. With this modification of the eigenvalues, the entropy conservative Roe flux is called ROE-EC1 flux where EC stands for ``entropy consistent". We will adopt the same eigenvalue modifications as in equation~(\ref{eq:lec1}) in combination with the new entropy conservative fluxes refered to as KEP-EC1 and repeat the computation of the stationary shock problem. The solutions shown on the right of figure~(\ref{fig:stshock}) indicate that the ROE-EC1 flux is able to give monotone solution at the lower Mach number but still produces some oscillations at the higher Mach numbers. The new approximately entropy consistent flux produces oscillations at all Mach numbers\footnote{These oscillations disappear if we use $\beta=1$.} while the new entropy conservative flux yields monotone solutions at all Mach numbers. This shows that exact entropy conservation does yield better solutions and is to be prefered. In order to improve the monotonicity of the ROE-EC1 flux, the acoustic eigenvalues can be further augmented leading to the ROE-EC2 flux as discussed in~\cite{Ismail:2009:AEE:1552569.1552681}. However with the new entropy conservative flux, we find that there is no need to increase the eigenvalues for higher Mach numbers and the basic Roe-type dissipation yields monotone solutions. The behaviour of the new flux is consistent with the entropy production analysis of \cite{Ismail:2009:AEE:1552569.1552681} which is valid only for weak shocks.

\begin{figure}
\begin{center}
\includegraphics[width=0.45\textwidth]{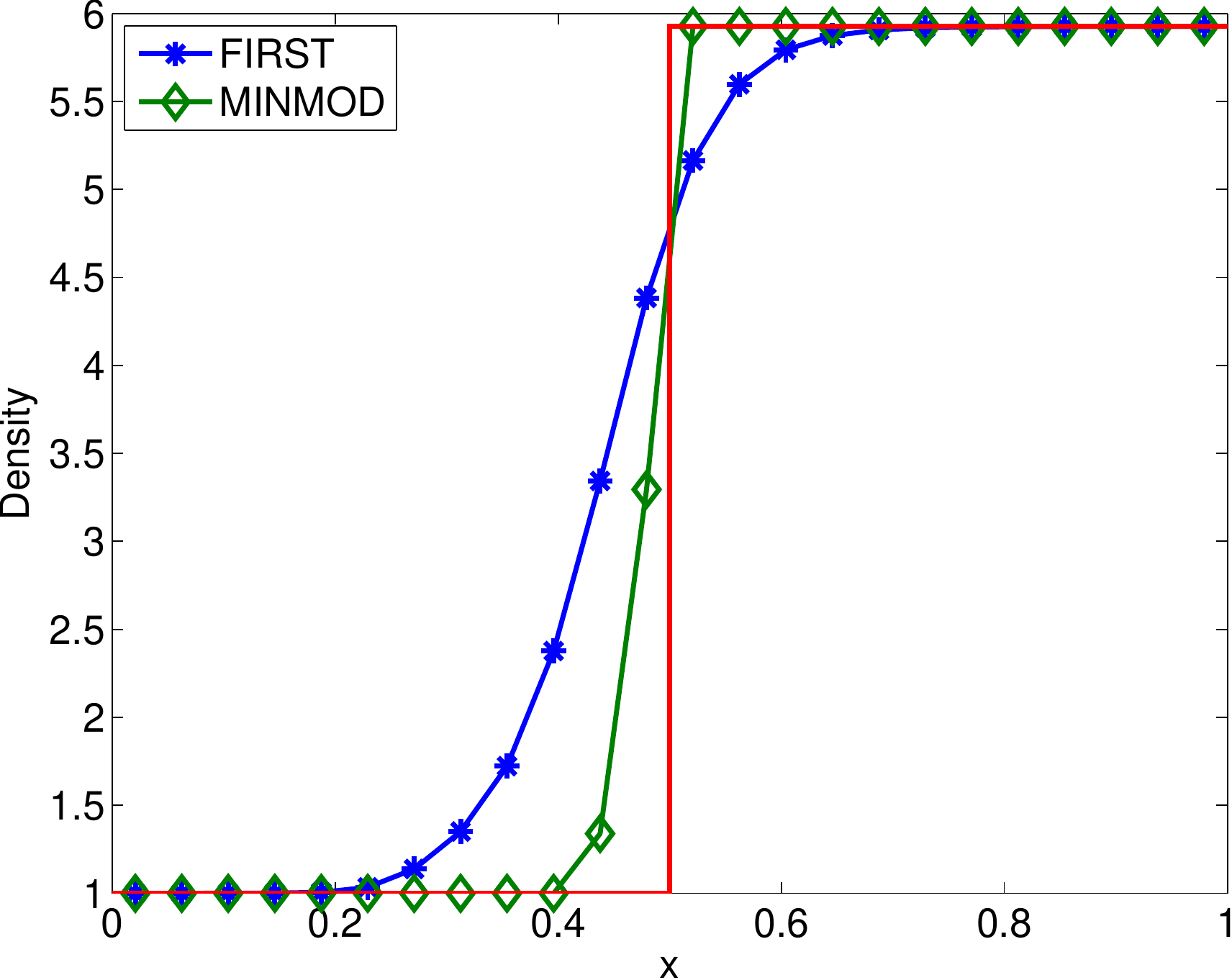}
\caption{Stationary shock problem: $N=24$ cells}
\label{fig:stshock_kes}
\end{center}
\end{figure}

We also compute the problem using the entropy dissipation that leads to kinetic energy stable scheme. Note that we do not use augmented eigenvalues as in the ROE-EC1 scheme but the first and third eigenvalues are taken to be equal to the maximum eigenvalue while the second eigenvalue is unchanged as in equation~(\ref{eq:kescontact}). 
The solution for Mach = 20 is shown in figure~(\ref{fig:stshock_kes}); the first order scheme is highly diffused due to the increased dissipation and the scheme is non-oscillatory. If we use the linear reconstruction scheme with minmod limiter, then a sharper shock resolution is obtained which is again non-oscillatory. Such kinetic energy and entropy stable schemes would be attractive for use with high order accurate schemes like WENO and discontinuous Galerkin methods.

\section{Shock instability, carbuncle and a hybrid scheme}
\label{sec:hyb}
Numerical schemes based on exact or approximate Riemann solvers suffer from a variety of pathological behaviours~\cite{FLD:FLD1650180603}. The entropy violating shock can be overcome by ensuring that the numerical scheme is consistent with the entropy condition. Other problems like the shock instability and carbuncle  do not have a proper rational fix available.

Consider the stationary shock problem as discussed in the previous section. If we set the initial conditions to exactly correspond to the discontinuity, then all the schemes yield stable, stationary solutions for which the residual converges to machine zero. However the Roe scheme and many other schemes are known to be unstable if the initial condition contains an intermediate point at the shock location. In order to induce the instability, the mass flux at the right boundary, which is an outflow boundary, is fixed to be one, while the momentum and energy fluxes are computed by zero gradient condition; this amounts to not updating the momentum and energy in the last cell. Small disturbances produced at the shock due to the intermediate point are propagated downstream, get reflected from the outflow boundary back into the shock and further disturb the shock location leading to a limit cycle oscillation. The new KEP-ES schemes do not suffer from the shock instability problem though the solutions are not monotone. The addition of extra dissipation in the KEP-EC1 scheme leads to the appearance of instability at higher Mach numbers.

It is thought that numerical schemes which suffer from this instability will also be prone to the carbuncle phenomenon which is found in flows with strong shocks over blunt bodies. However there is no guarantee that a scheme which is stable for the 1-d shock problem will not produce the carbuncle phenomenon. The ROE-EC2 flux which has increased dissipation avoids the 1-D shock instability but can still suffer from the carbuncle problem in multi-dimensions~\cite{ismailcfd2006}. We show in the later sections that the KEP-EC1 scheme also produces carbuncle effect on certain types of meshes where the shock is aligned with the mesh. Increasing the value of $\beta$ in equation~(\ref{eq:lec1}) even upto $\beta=1$ does not seem to eliminate this problem. The kinetic energy stable scheme given by equation~(\ref{eq:kescontact}) which has more dissipation in the acoustic waves also produces the carbuncle but the Rusanov scheme of equation~(\ref{eq:rus}) does not produce the carbuncle. In fact what we observe is that all the schemes which resolve grid aligned stationary contacts exactly seem to suffer from the carbuncle effect which is consistent with what is noticed in other schemes in the literature. The usual fix in such cases is to increase the amount of dissipation in the numerical scheme which however causes poor resolution of boundary layers. There is also the idea of switching the numerical scheme to a more dissipative one only near shocks and using a high resolution Riemann solver type scheme in smooth parts of the flow~\cite{FLD:FLD1650180603}. In the framework of the entropy conservative/stable scheme as discussed in this paper, we have the freedom in designing the eigenvalues which essentially control the amount of dissipation in the scheme. Since the Rusanov scheme is free of carbuncles, we propose a blending of the usual Roe scheme with the Rusanov scheme, i.e., the matrix $|\Lambda|$ appearing in the entropy dissipation flux is of the form
\[
|\Lambda| = |\Lambda|^{Hyb} = (1-\phi) |\Lambda|^{Roe} + \phi |\Lambda|^{Rus}
\]
where the switching function $\phi$ is based on the pressure jump 
\[
\phi = \left| \frac{\Delta p}{2\avg{p}} \right|^\frac{1}{2}
\]
Note that $0 \le \phi \le 1$; for a strong shock $\phi \approx 1$ and the scheme is close to the Rusanov scheme, while near a contact discontinuity $\phi \approx 0$ and the scheme is close to the more accurate Roe scheme. This hybrid scheme does not produce the 1-D shock instability problem. In the results section we show that this scheme avoids the carbuncle effect and still gives good resolution of boundary layers and shear layers. The new entropy conservative flux together with the above hybrid dissipation will be refered to as KEP-ES(Hyb) scheme. We also  remark that these modifications of the dissipation term still retain the entropy stability property of the scheme.
\section{2-D NS equations and finite volume method}
\label{sec:2d}
The two dimensional Navier-Stokes equations in conservation form can be written as
\[
\df{\con}{t} + \df{\fx_1}{x_1} + \df{\fx_2}{x_2} = \df{\gx_1}{x_1} + \df{\gx_2}{x_2}
\]
where
\[
\con = \begin{bmatrix}
\rho \\
\rho u_1 \\
\rho u_2 \\
E
\end{bmatrix}, \qquad \fx_1 = \begin{bmatrix}
\rho u_1 \\
p + \rho u_1^2 \\
\rho u_1 u_2 \\
(E+p) u_1
\end{bmatrix}, \qquad \fx_2 = \begin{bmatrix}
\rho u_2 \\
\rho u_1 u_2 \\
p + \rho u_2^2 \\
(E+p) u_2
\end{bmatrix}
\]
with $u=(u_1,u_2)$ being the Cartesian components of the fluid velocity and $E=p/(\gamma-1) + \rho |u|^2/2$ is the total energy per unit volume, while the viscous fluxes are given by
\[
\gx_1 = \begin{bmatrix}
0 \\
\tau_{11} \\
\tau_{21} \\
\tau_{11}u_1 + \tau_{12}u_2 - q_1
\end{bmatrix}, \qquad \gx_2 = \begin{bmatrix}
0 \\
\tau_{12} \\
\tau_{22} \\
\tau_{21}u_1 + \tau_{22}u_2 - q_2
\end{bmatrix}
\]
where $\tau_{ij}$ is the shear stress tensor and $q_i$ is the heat flux vector, for which we assume the Newtonian and Fourier constitutive laws, respectively. We will approximate these equations using a finite volume method on unstructured grids. Specifically we will use a vertex-based finite volume scheme in which a primal grid of triangles is used; the finite volumes are constructed around each vertex using either a {\em median dual cell} or a {\em voronoi cell}. In the median dual cell, the finite volumes are constructed by joining the cell centroids to the edge mid-points while in the voronoi dual, the circumcenter of the triangles are joined to form the finite volumes. If the triangle is obtuse angled so that the circumcenter lies outside the triangle, then the mid-point of the largest side is used instead of the circumcenter.

Let $A_i$ be the area of the $i$'th vertex/finite volume and let $N(i)$ denote the set of neigbouring finite volumes which share a common boundary with $A_i$. We also use the set $C(i)$ which is the set of all triangles $T$ containing the $i$'th vertex. The semi-discrete finite volume scheme is given by
\[
A_i \dd{\con_i}{t} + \sum_{j \in N(i)} \fx_{ij} \Delta s_{ij} = \sum_{T \in C(i)} \gx^T_i \Delta s_{i}^T
\]
where $\fx_{ij}$ is the flux from the $i$'th cell into the $j$'th cell across their common boundary whose length is $\Delta s_{ij}$. The discretization of viscous fluxes uses a $P_1$ finite element approach on triangular grids which leads to only a nearest neighour stencil (edge connected neighbours) and hence is very compact. The flux $\gx^T_i$ is the contribution of the dissipative flux to vertex $i$ coming from the portion of the boundary of $A_i$ lying inside triangle $T$ and $\Delta s^T_i$ is the length of this boundary. For a boundary finite volume, there will be flux contributions from the boundary edges also. We will approximate this flux using the entropy conservative flux together with some dissipation which takes the form
\[
\fx_{ij} = \fx_{ij}^* - \frac{1}{2} R_{ij} |\Lambda_{ij}| S_{ij} R_{ij}^\top \Delta\ent_{ij}
\]
The entropy conservative flux $\fx^*$ in the direction of the unit normal vector $n=(n_1, n_2)$ to the cell face is given by
\[
\fx^* = \begin{bmatrix}
\alog{\rho} (\avg{u} \cdot n) \\
\tp n_1 + \avg{u}_1 f^{*,\rho} \\
\tp n_2 + \avg{u}_2 f^{*,\rho} \\
\left[ \frac{1}{2(\gamma-1)\alog{\beta}} - \frac{1}{2}\avg{|u|^2} \right] f^{*,\rho} + \avg{u} \cdot f^{*,m}
\end{bmatrix}
\]
where $\alog{\rho}$, $\alog{\beta}$ are the logarithmic averages and $\tp$ is given by equation~(\ref{eq:tp}). The matrix of eigenvectors $R$ is given by
\[
R = \begin{bmatrix}
1 &             1 &            0 &                     1 \\
u_1 - a n_1 &  u_1 &           n_2 &  u_1 + a n_1 \\
u_2 - a n_2 &  u_2 &           -n_1 & u_2 + a n_2 \\
H - a (u \cdot n) &    \frac{1}{2}|u|^2 & u_1 n_2 - u_2 n_1 & H     + a (u \cdot n)
\end{bmatrix}
\]
and the scaling matrix $S$ is
\[
S = \textrm{diag}\begin{bmatrix} \frac{1}{2\gamma}\rho & \frac{(\gamma-1)}{\gamma}\rho & p & \frac{1}{2\gamma}\rho \end{bmatrix}
\]
The eigenvalue matrix for the basic Roe scheme (ES) is given by
\[
|\Lambda| = |\Lambda|^{Roe} = \textrm{diag}\begin{bmatrix} |u \cdot n - a| & |u \cdot n| & |u \cdot n| & |u \cdot n + a| \end{bmatrix}
\]
and the modifications for the EC1, KES, Rus and Hyb schemes are obvious from their 1-D counterparts. Finally, the entropy variables are given by
\[
\ent = \begin{bmatrix}
\frac{\gamma-s}{\gamma-1} - \beta|u|^2 & 2\beta u_1 & 2\beta u_2 & -2\beta
\end{bmatrix}^\top
\]
Second order accuracy is achieved by a MUSCL type reconstruction using nodal gradients of primitive variables $p, u_1, u_2, T$ which are calculated using Green-Gauss theorem while time-stepping is done using a 3-stage strong stability preserving Runge-Kutta scheme~\cite{Shu1988439}. Triangular grids used in the computations are generated with GMSH~\cite{NME:NME2579} and Delaundo~\cite{muller}.
\section{Numerical results}
\label{sec:results}
In all the test problems, we consider an ideal gas with $\gamma=1.4$ except when indicated otherwise.
\subsection{Modified Sod test case}

\begin{figure}
\begin{center}
\includegraphics[width=0.35\textwidth]{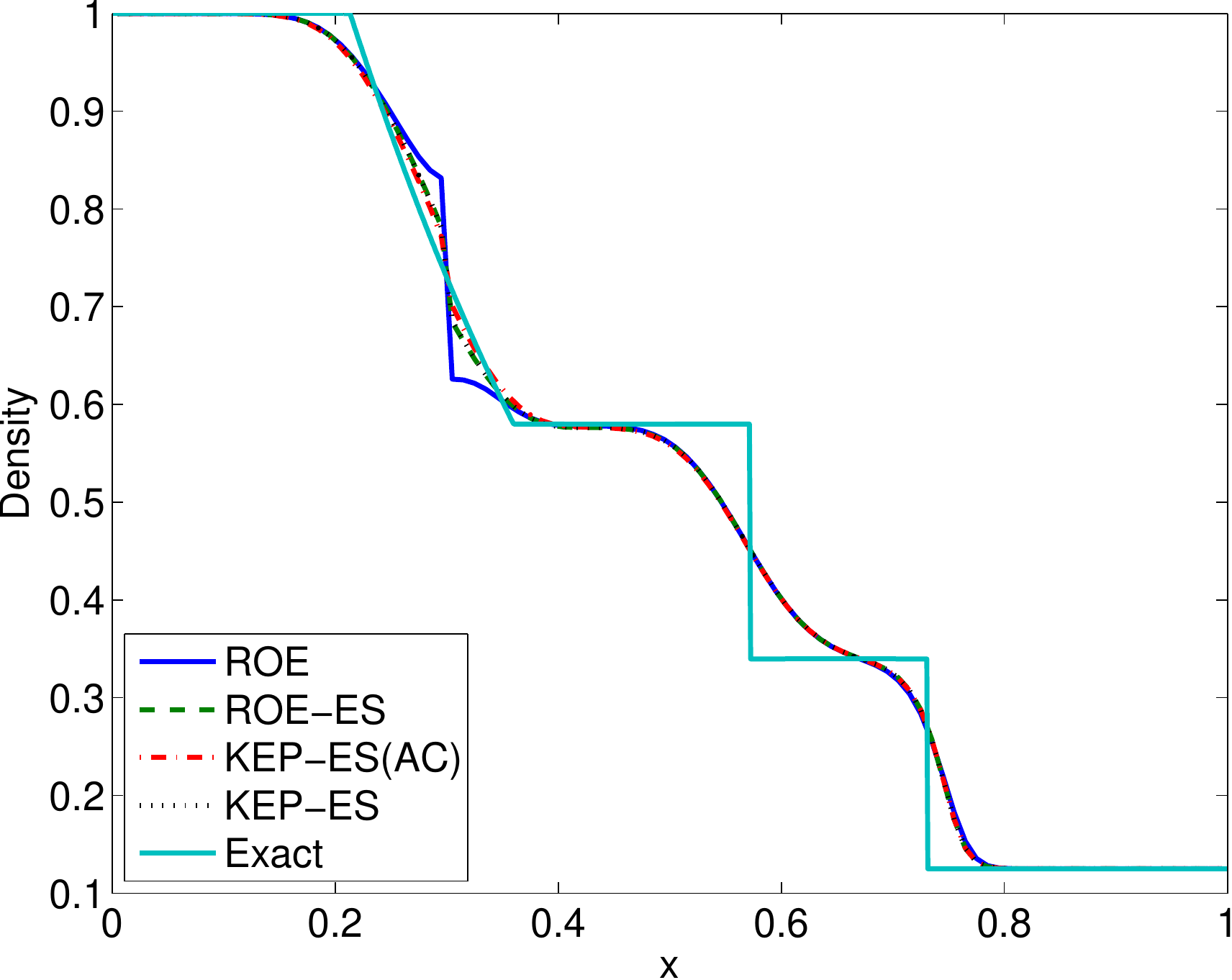}
\includegraphics[width=0.35\textwidth]{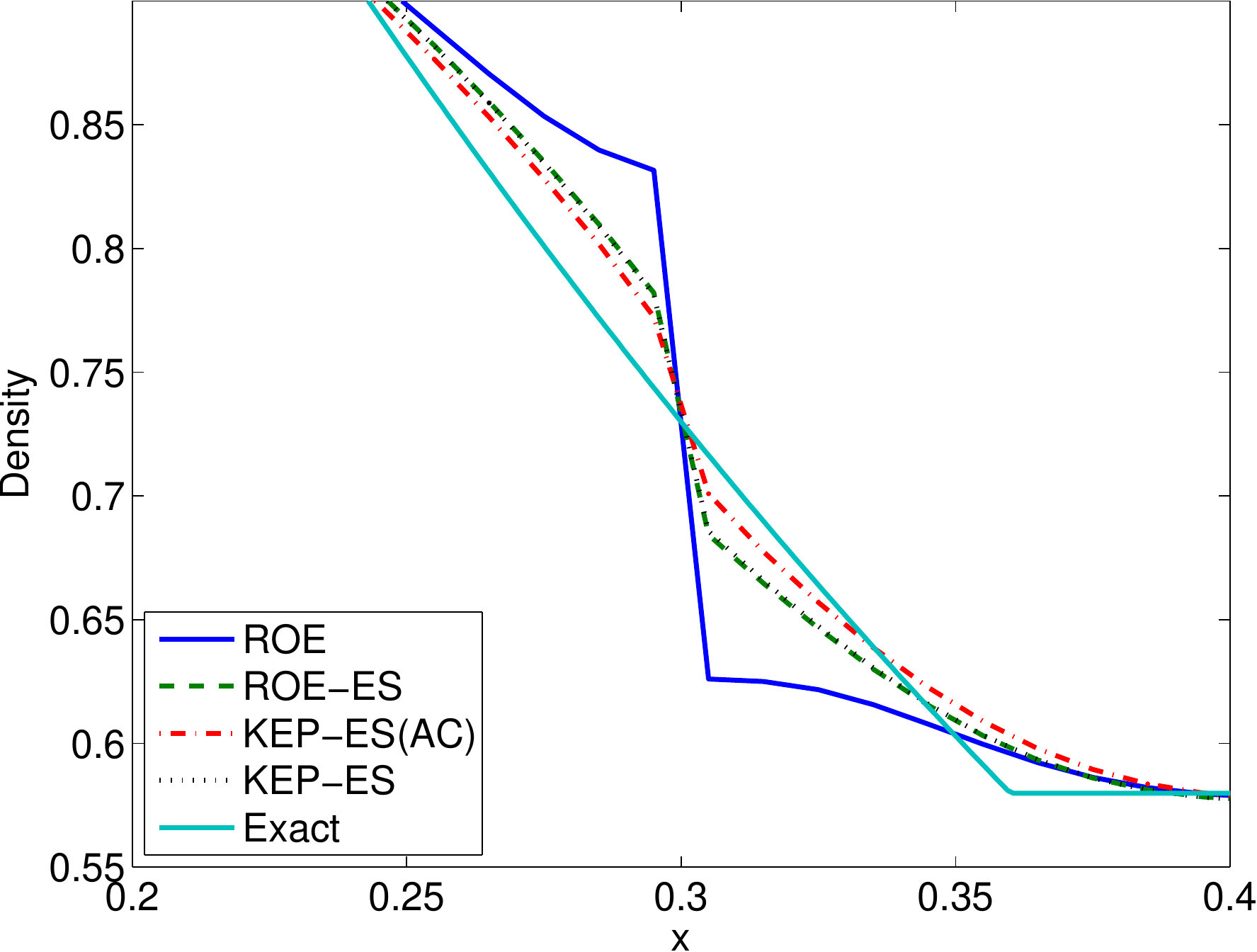} \\
\includegraphics[width=0.35\textwidth]{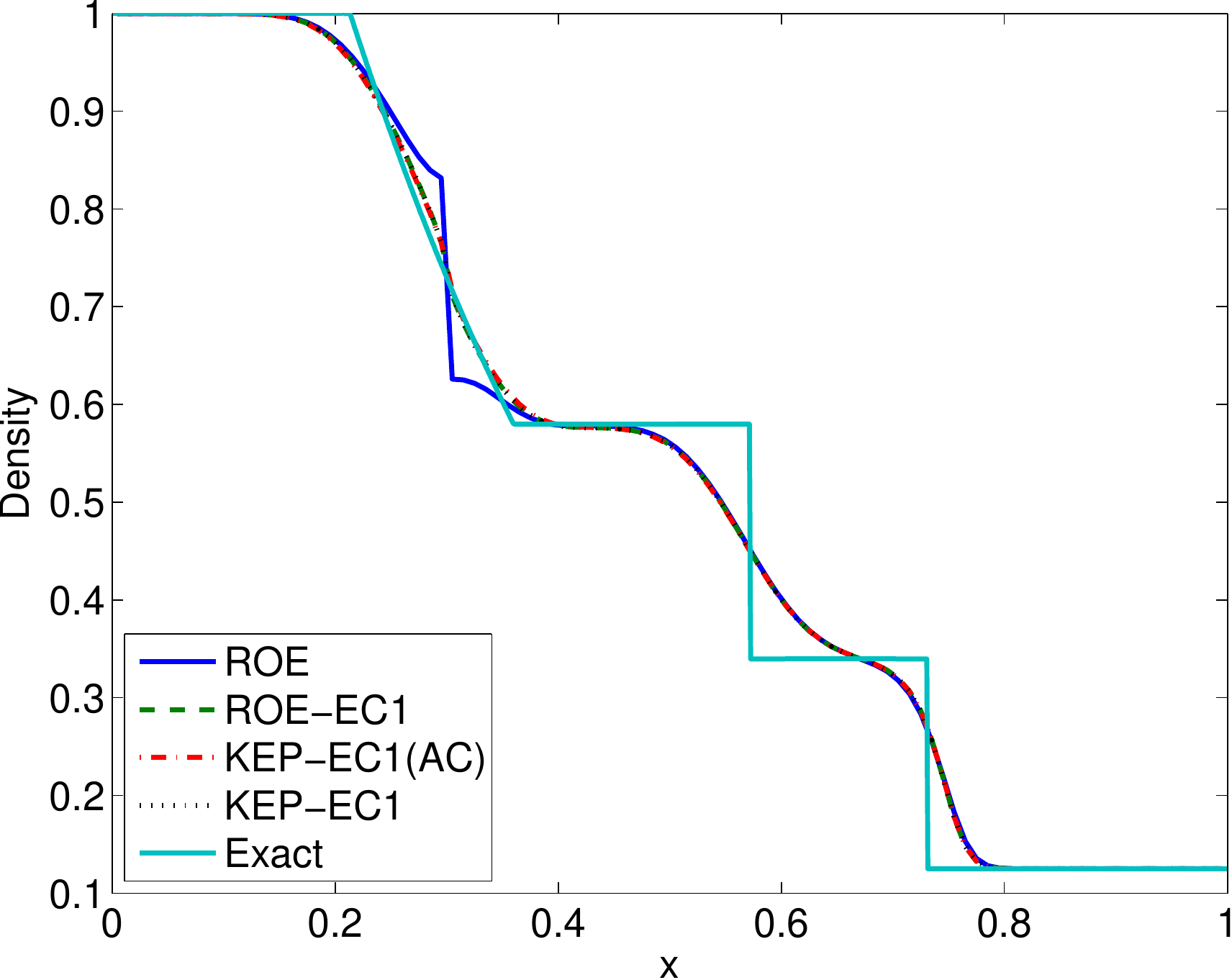}
\includegraphics[width=0.35\textwidth]{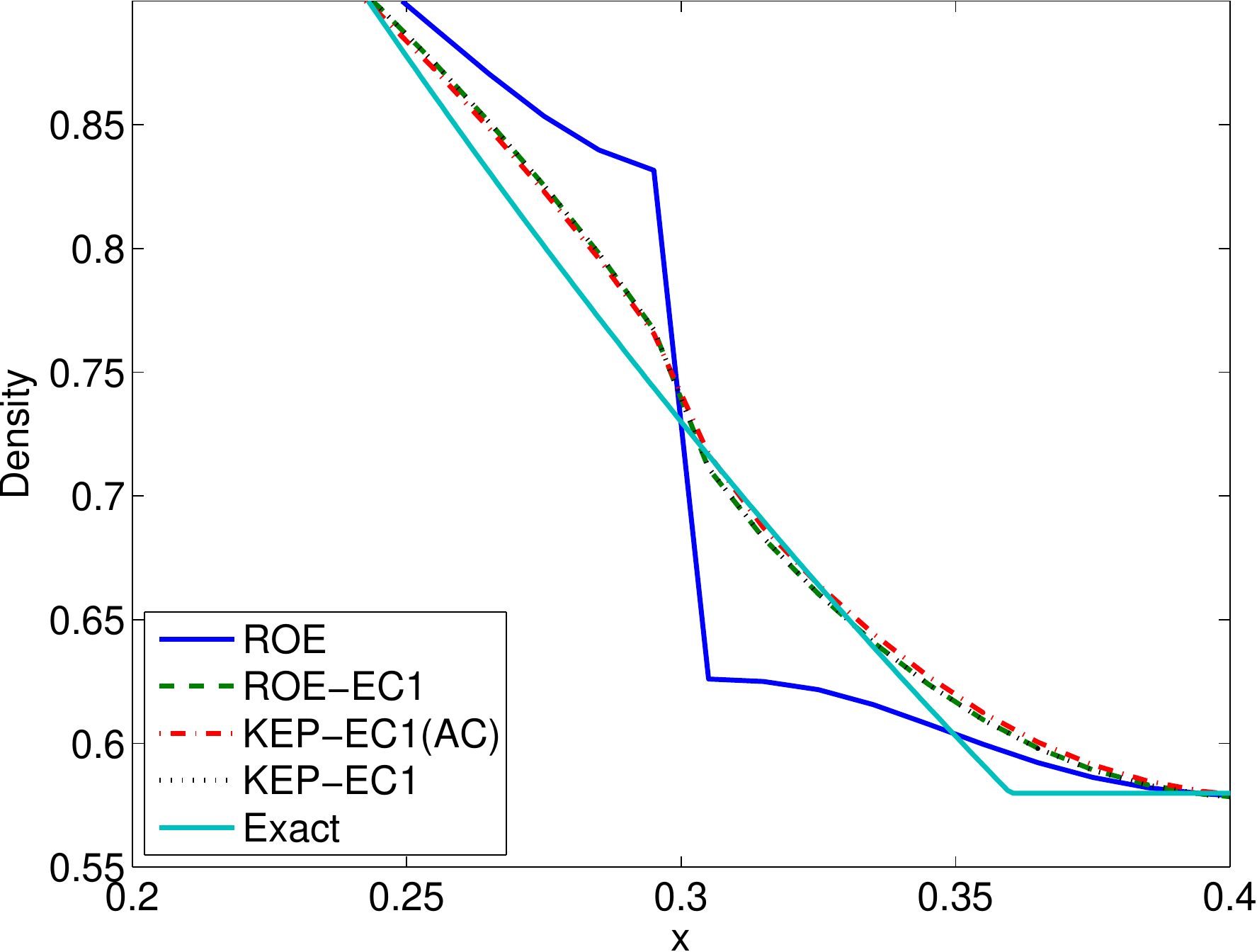} \\
\caption{Modified Sod problem: Density at time $t=0.2$. The figures on the right show a zoomed view of the rarefaction region.}
\label{fig:modsod}
\end{center}
\end{figure}

This is a shock tube problem with the left state being $(\rho,u,p)=(1.0, 0.75,  1.0)$ and the right state being $(\rho,u,p)=(0.125, 0.0, 0.1)$. All computations are made with $N=100$ cels and a CFL=0.4 upto a final time of $t=0.2$ units. The original    Roe scheme gives entropy violating jump in the expansion region where the flow becomes sonic as shown in figure~(\ref{fig:modsod}). The reason for this is the vanishing of one eigenvalue at the sonic point which leads to insufficient dissipation. All the entropy consistent schemes including the new approximately entropy consistent scheme together    with entropy variable based dissipation give solutions which do not suffer     from this problem as shown in figure~(\ref{fig:modsod}). The zoomed view of the figure on the right shows that in the sonic region of the expansion fan all of them behave in a very similar manner. Note that even the entropy stable schemes have a vanishing eigenvalue in the expansion fan. However due to the entropy conservative nature of the central part of the flux $\fx^*$, they do not give rise to the entropy violating shock unlike the original Roe scheme which is not entropy conservative.
\subsection{Stationary contact}
\begin{figure}
\begin{center}
\includegraphics[width=0.48\textwidth]{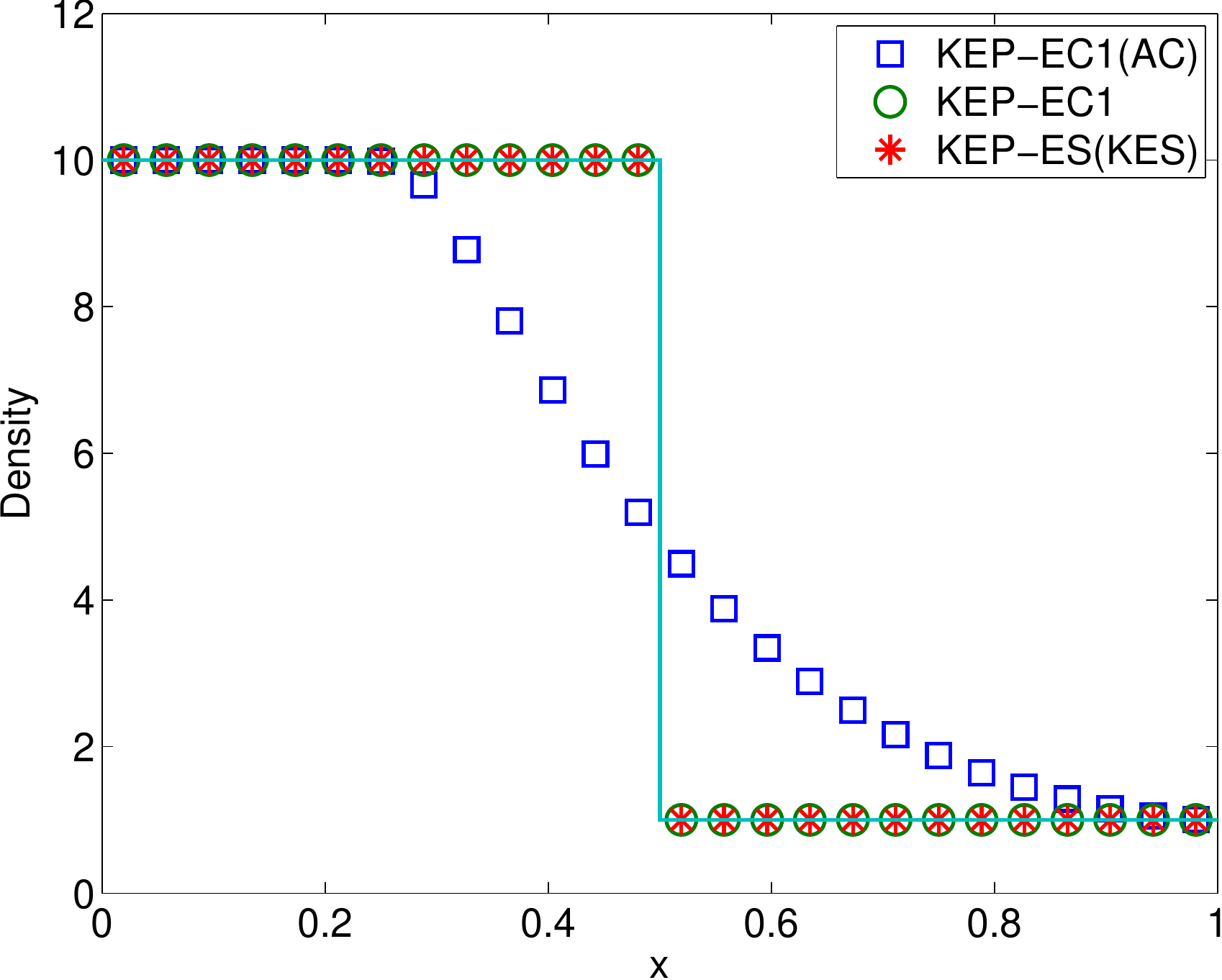}
\caption{Stationary contact}
\label{fig:stcontact}
\end{center}
\end{figure}

We consider inviscid Riemann problem with left state $(\rho_l,u_l,p_l)=(10,0,1)$ and right state $(\rho_r,u_r,p_r)=(1,0,1)$ whose solution is a stationary contact wave. The problem is solved with $N=26$ cells and the solution is shown in figure~(\ref{fig:stcontact}). The condition for the exact resolution of stationary contacts is the correct computation of the enthalpy in the dissipation matrix. For the approximately entropy consistent flux (AC), we use the arithmetic average of $\beta$ to compute the enthalpy which does not satisfy the exact contact resolution property. For the exactly conservative flux, we use the logarithmic average of $\beta$ which should resolve the contact wave exactly. The exactly conservative flux is also used in combination with the kinetic energy stable scheme KEP-ES(KES) which should resolve the contact discontinuity exactly. The results of these three schemes are consistent with the theoretical prediction as seen in the figure~(\ref{fig:stcontact}). These results show that contact waves can be highly dissipated if the numerical flux is not able to resolve contacts exactly.
\subsection{Sod test case}
\label{sec:nssod}
\begin{figure}
\begin{center}
\includegraphics[width=0.35\textwidth]{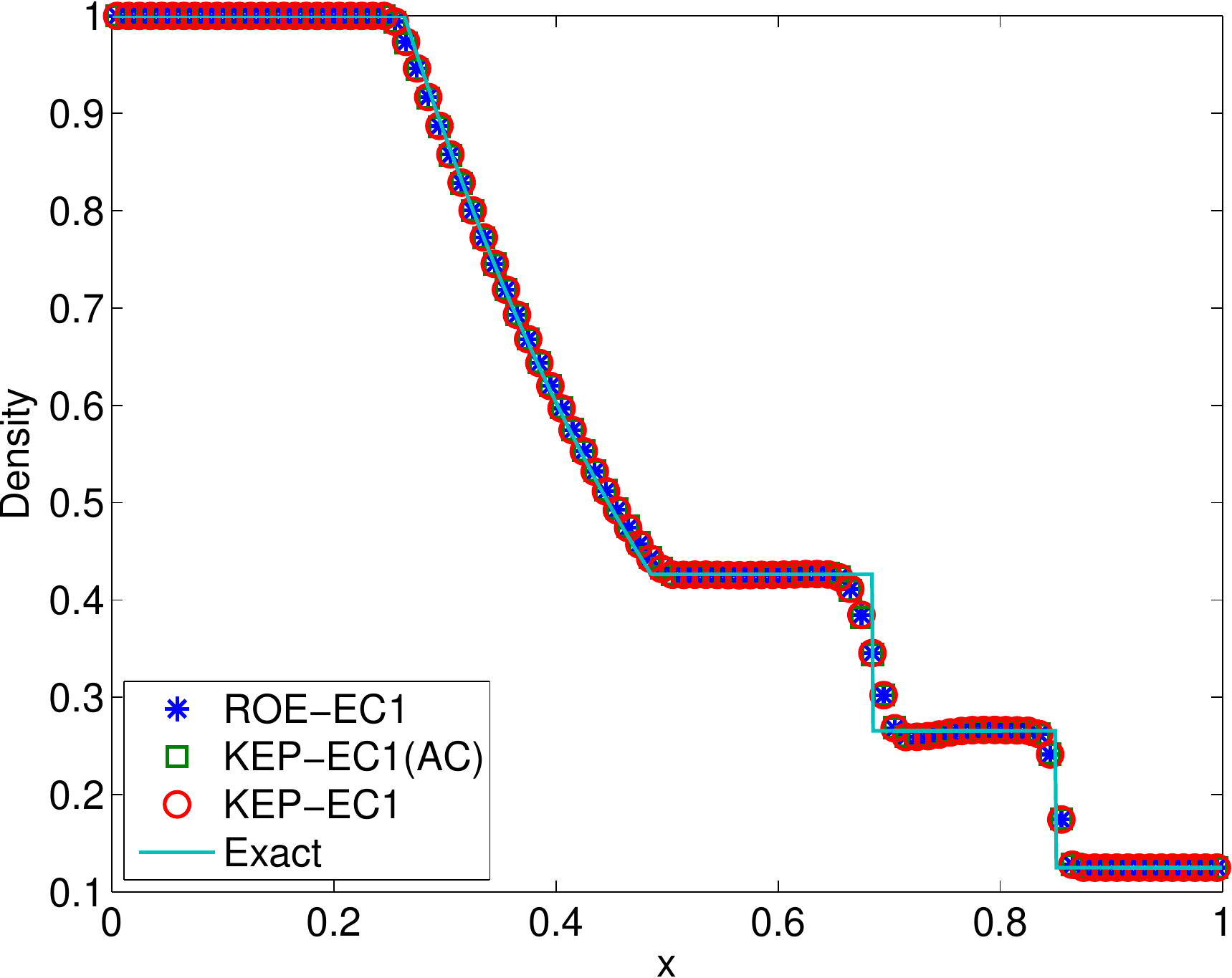}
\includegraphics[width=0.35\textwidth]{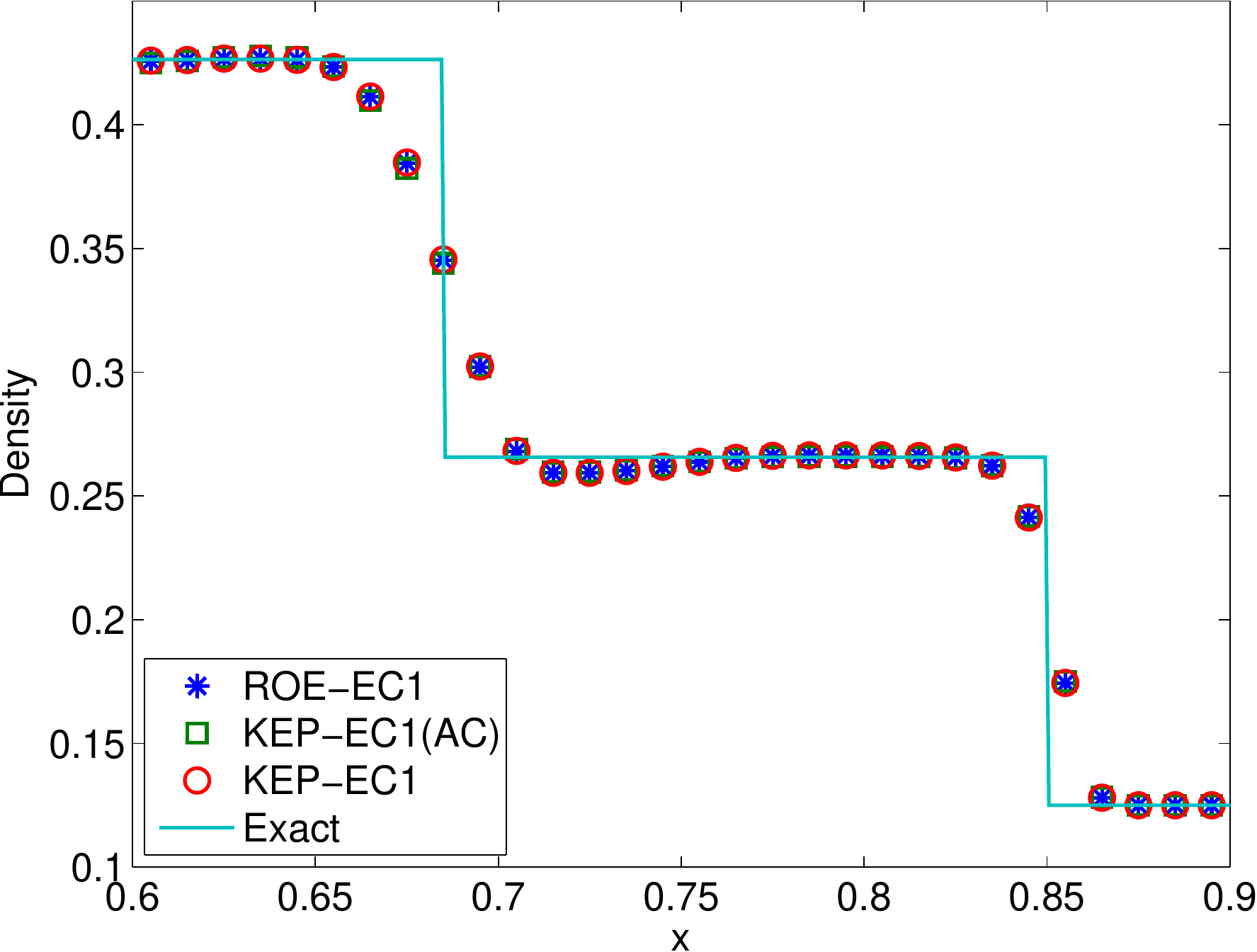} \\
\includegraphics[width=0.35\textwidth]{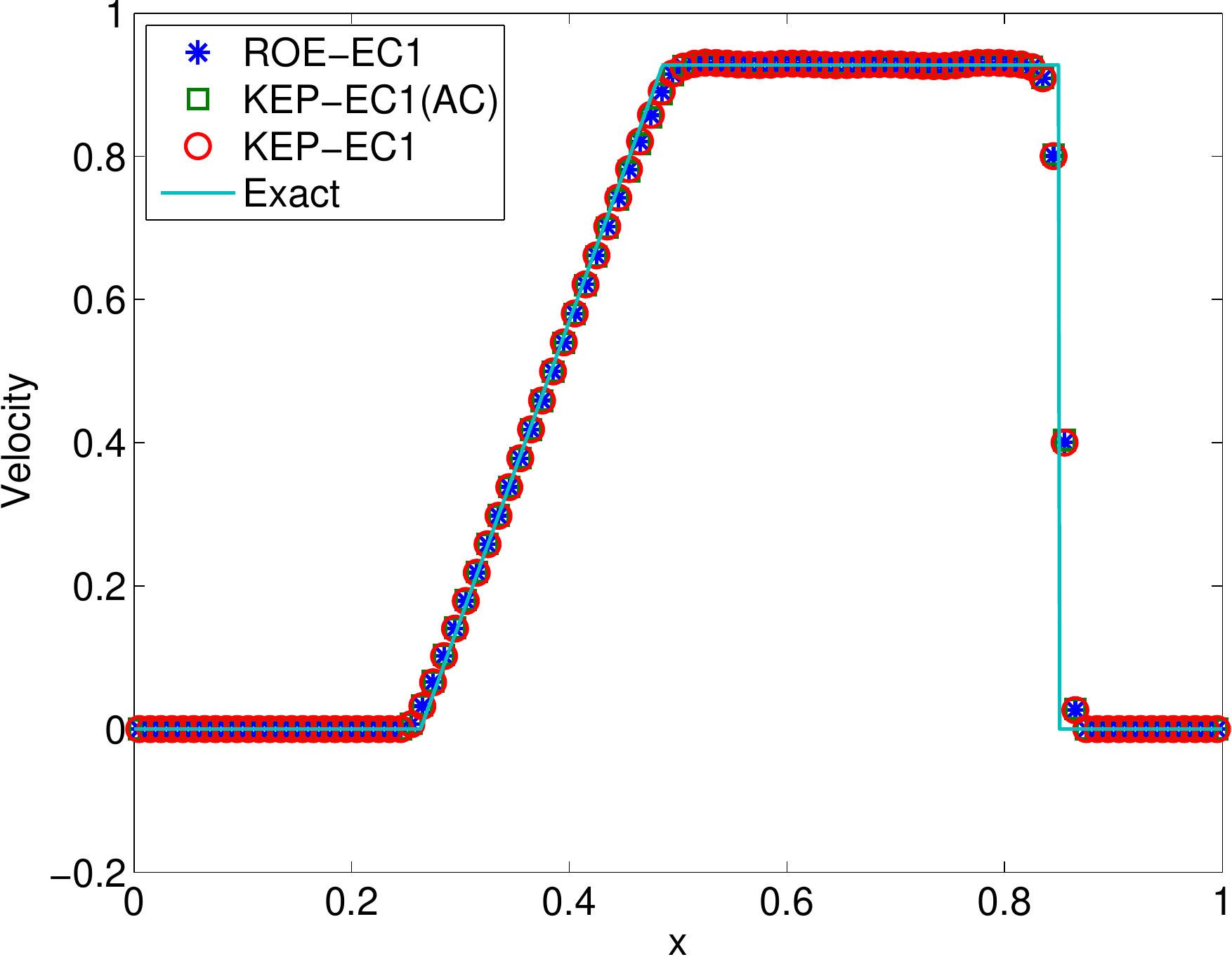}
\includegraphics[width=0.35\textwidth]{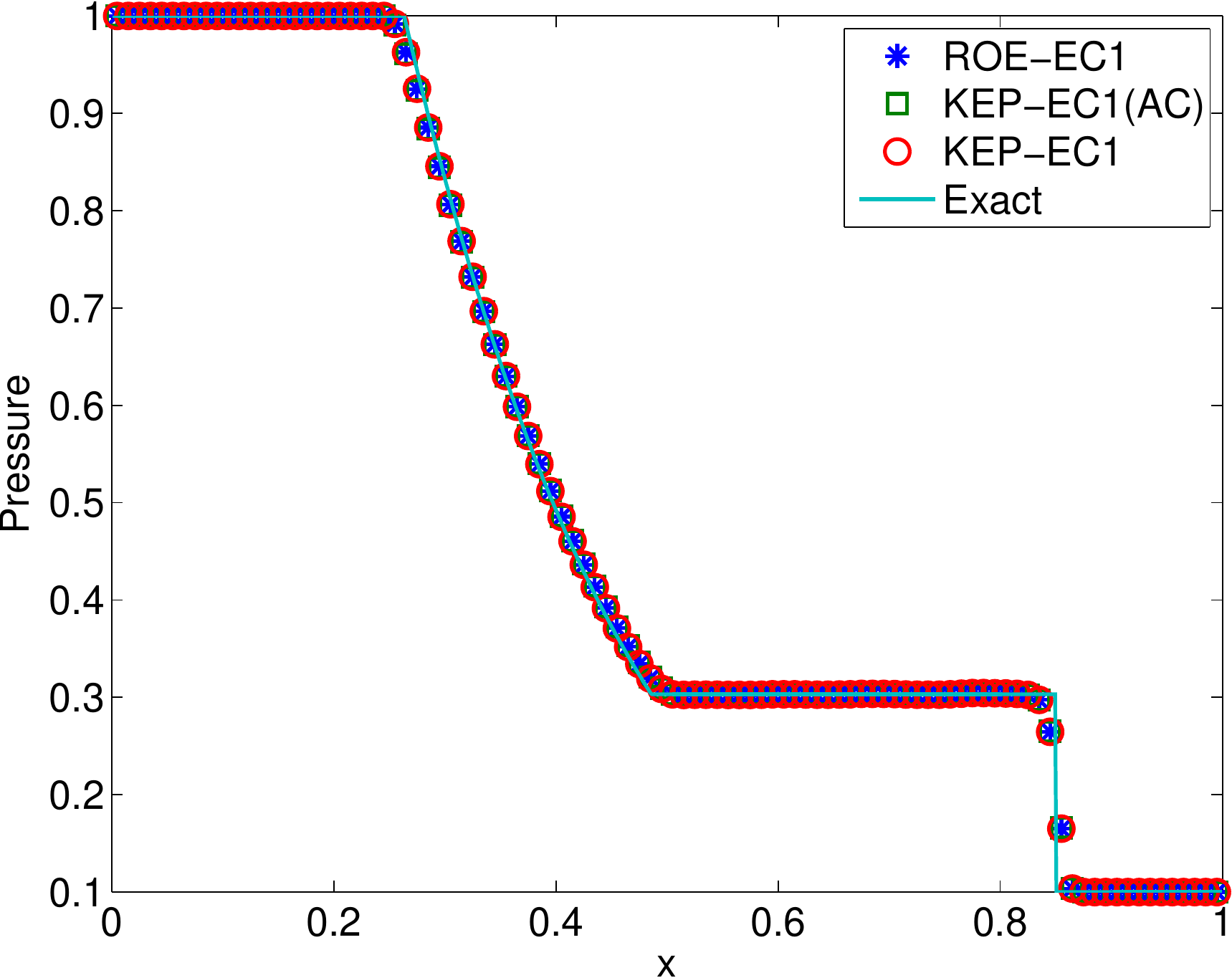} \\
\caption{Sod problem}
\label{fig:sod}
\end{center}
\end{figure}

{\em Inviscid case}: This is also a shock tube problem with the left state being $(\rho,u,p)=(1.0, 0.0, 1.0)$ and the right state being $(\rho,u,p)=(0.125, 0.0, 0.1)$. All computations are made with $N=100$ cells and a CFL=0.4 and a 3-stage Runge-Kutta scheme upto a final time of $t=0.2$ units. We compute   the solution using entropy consistent scheme with entropy dissipation and using MUSCL scheme and minmod limiter. The solution shown in figure~(\ref{fig:sod})   indicates that the sharp resolution of contact and shocks can be achieved with  the new schemes. The shock is resolved within two cells and the contact is resolved with about four cells. It can be seen that all the schemes give essentially the same solution on this problem.

\begin{figure}
\begin{center}
\begin{tabular}{cc}
\includegraphics[width=0.48\textwidth]{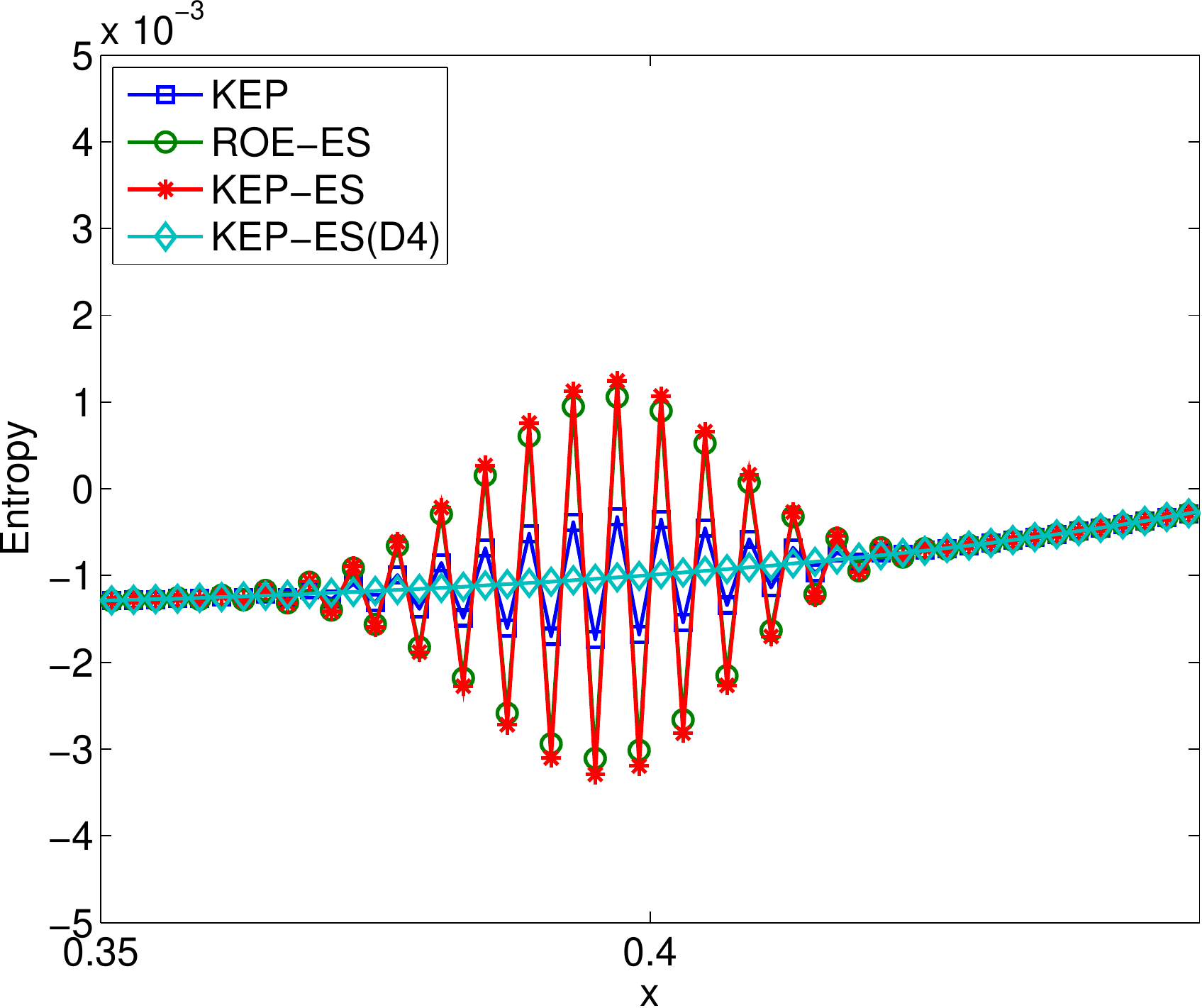} &
\includegraphics[width=0.48\textwidth]{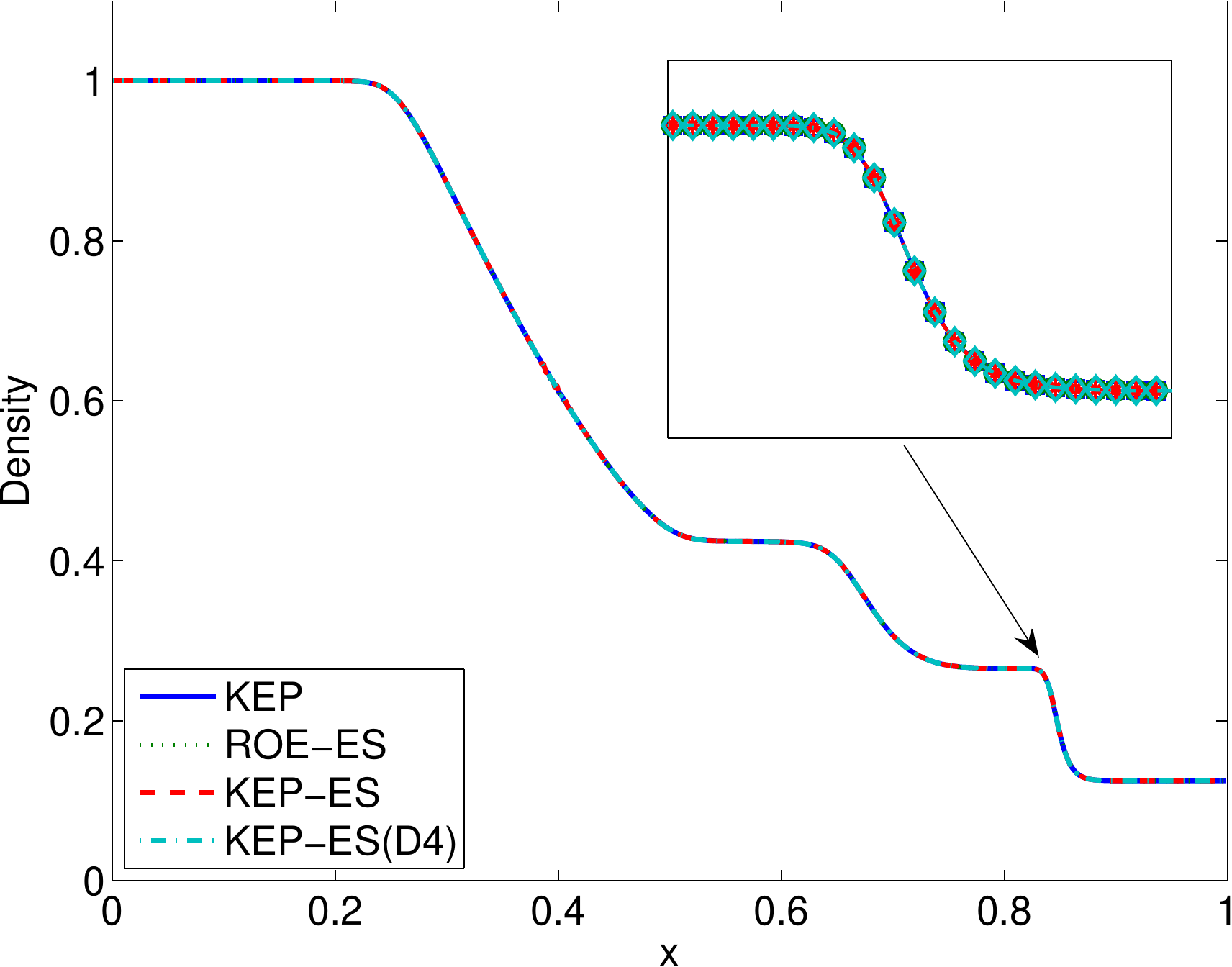} \\
(a) & (b)
\end{tabular}
\caption{Sod problem using Navier-Stokes equations: Solution at time $t=0.2$ using 500 cells (a) Zoomed view of entropy (b) Density}
\label{fig:nssod}
\end{center}
\end{figure}

\noindent
{\em Viscous case}: The same problem is solved using Navier-Stokes equations at a Reynolds number of 2000 based on the sound speed of the left state and the length of the domain using 500 cells and a CFL number of 0.1. We first use the central fluxes without any dissipation; the solutions are obtained with the kinetic energy preserving scheme of Jameson (KEP), Roe's entropy conservative scheme ROE-ES and the new KEP-ES scheme. The entropy $s = \ln(p \rho^{-\gamma})$ is shown in figure~(\ref{fig:nssod}a) in a zoomed section in the expansion region. The density shown in figure~(\ref{fig:nssod}b) also oscillates in the expansion region but the oscillations are too small to be seen on the scale of the figure. It is found that all the central schemes produce some oscillations in density and pressure which is also reflected in the entropy, while there are no oscillations in velocity and temperature. Reducing the CFL number does not eliminate these oscillations. These oscillations originate in the discontinuous initial condition, particularly density and pressure, and propagate upstream. The KEP scheme is seen to produce smaller oscillations compared to the entropy conservative schemes. This indicates that the KEP scheme has some inherent dissipation which is able to damp the oscillations in density but it is not able to eliminate them completely. The entropy conservative schemes do not have this additional dissipation probably because of their entropy conservative nature. Note that the physical viscosity which acts due to the velocity and temperature gradients is ineffective in damping oscillations in density. This problem is present even with implicit kinetic energy preserving schemes and in~\cite{Subbareddy20091347} the authors modify the pressure flux to be biased instead of centered; this adds some dissipation which becomes active if there are oscillations in pressure. In the present work, we also perform computations using the KEP-ES scheme but adding only fourth order scalar dissipation (D4) with $\kapf=\frac{1}{100}$. As shown in figure~(\ref{fig:nssod}), the fourth order dissipation is able to damp the density oscillations without affecting the accuracy of the solution in other regions. In particular the zoomed view in the shock region shows that the fourth order dissipation does not spoil the accuracy of the solution even inside the shock region.
\subsection{NS shock structure}

\begin{figure}
\begin{center}
\includegraphics[width=0.34\textwidth]{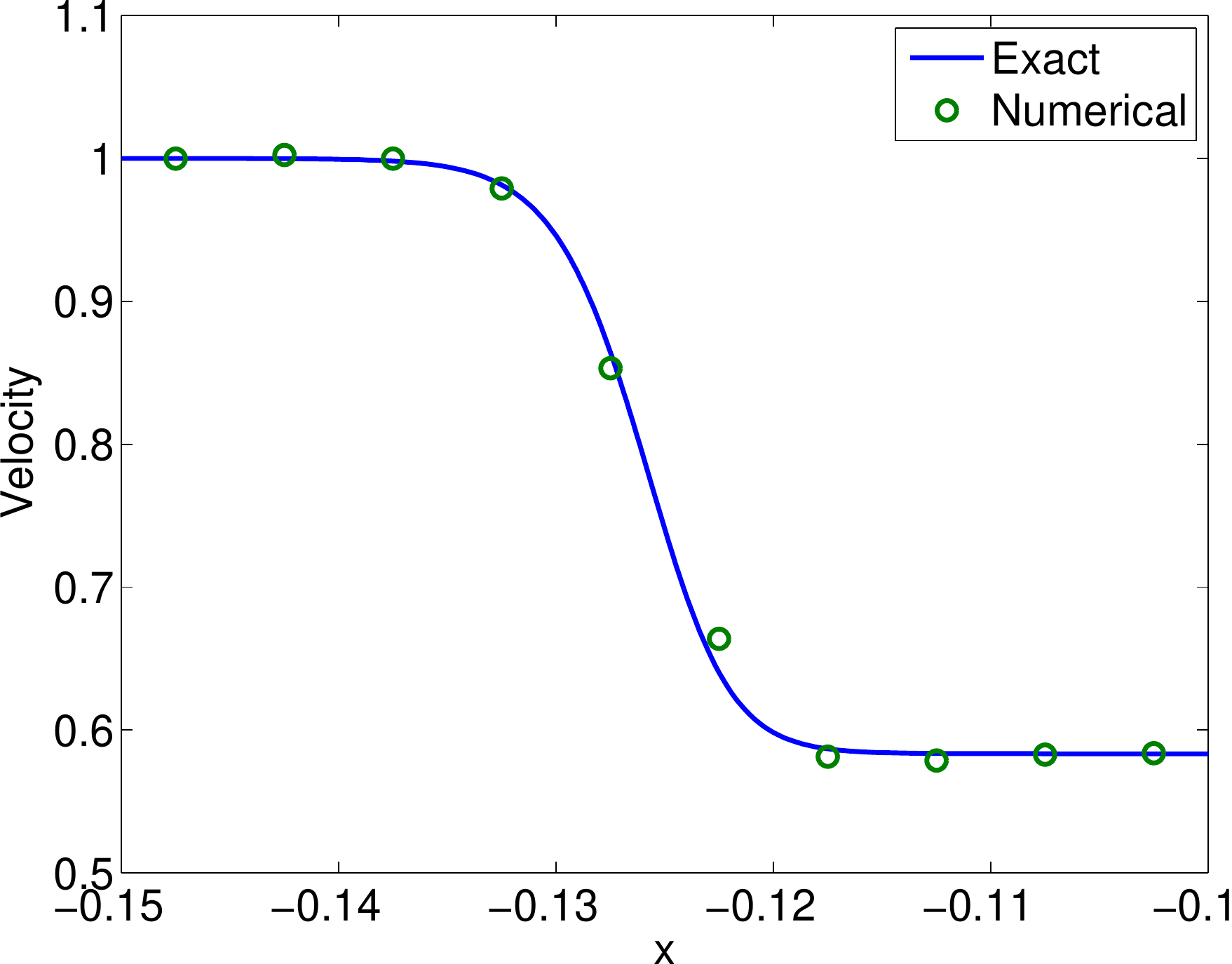}
\includegraphics[width=0.34\textwidth]{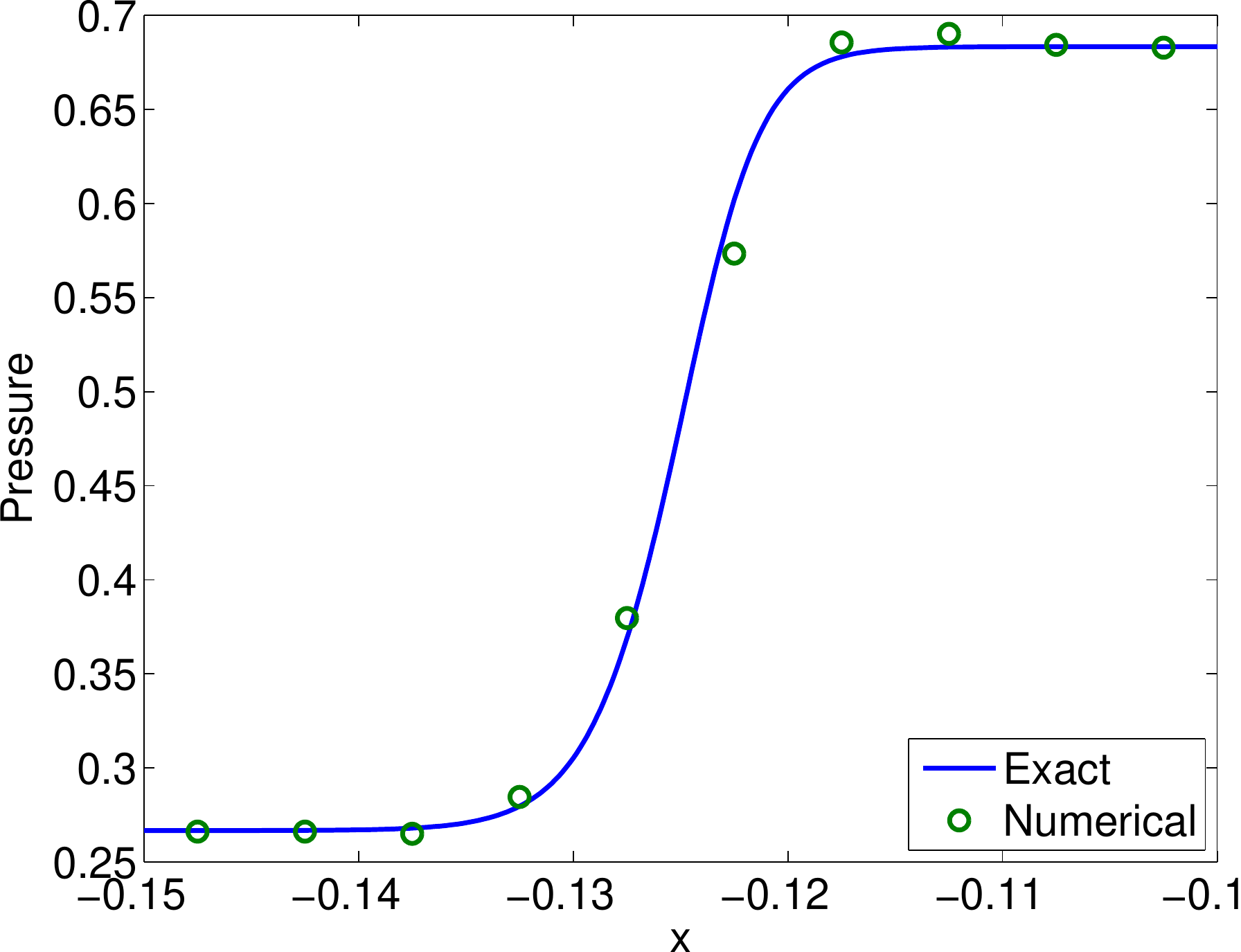} \\
\includegraphics[width=0.34\textwidth]{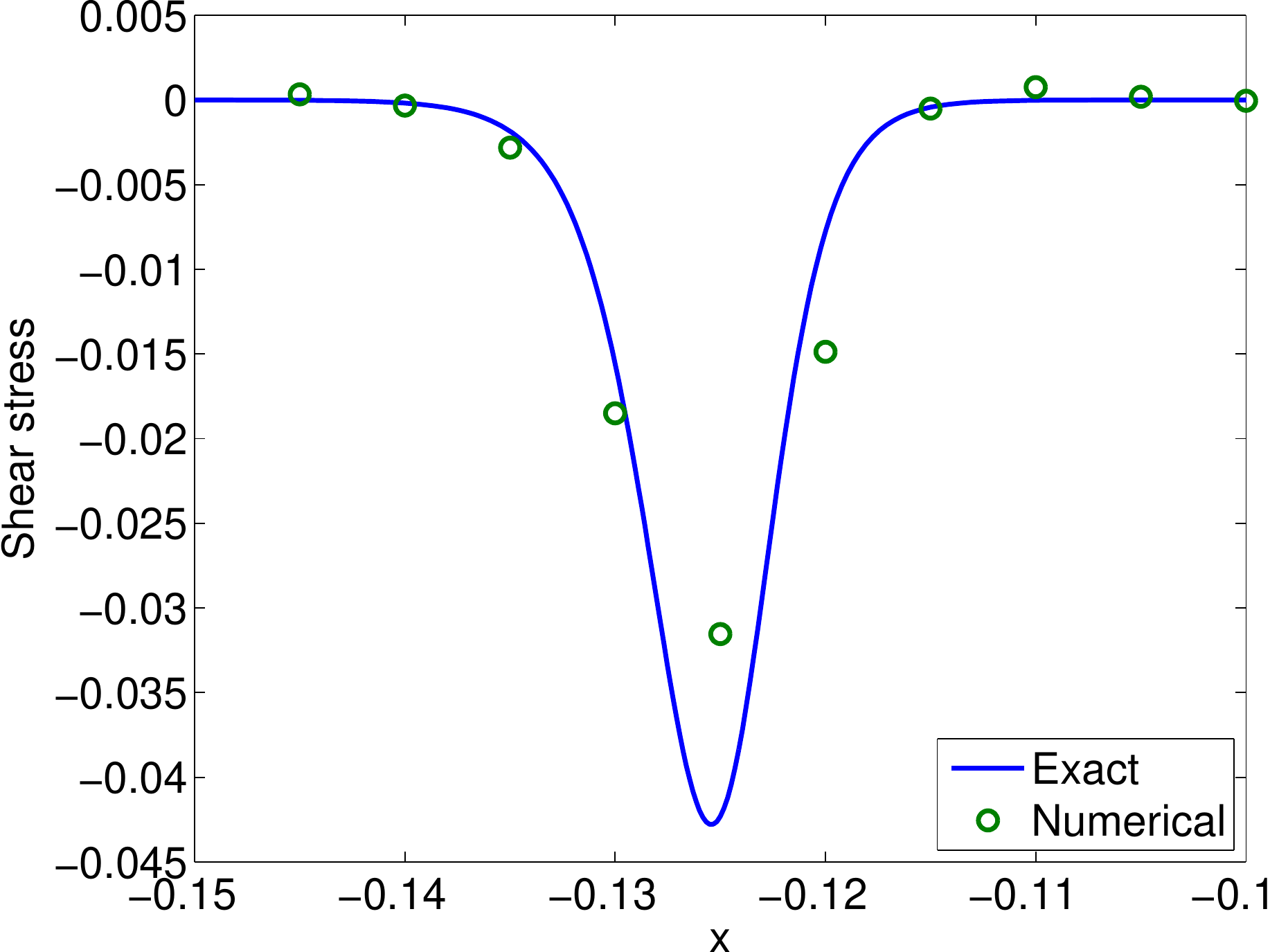}
\includegraphics[width=0.34\textwidth]{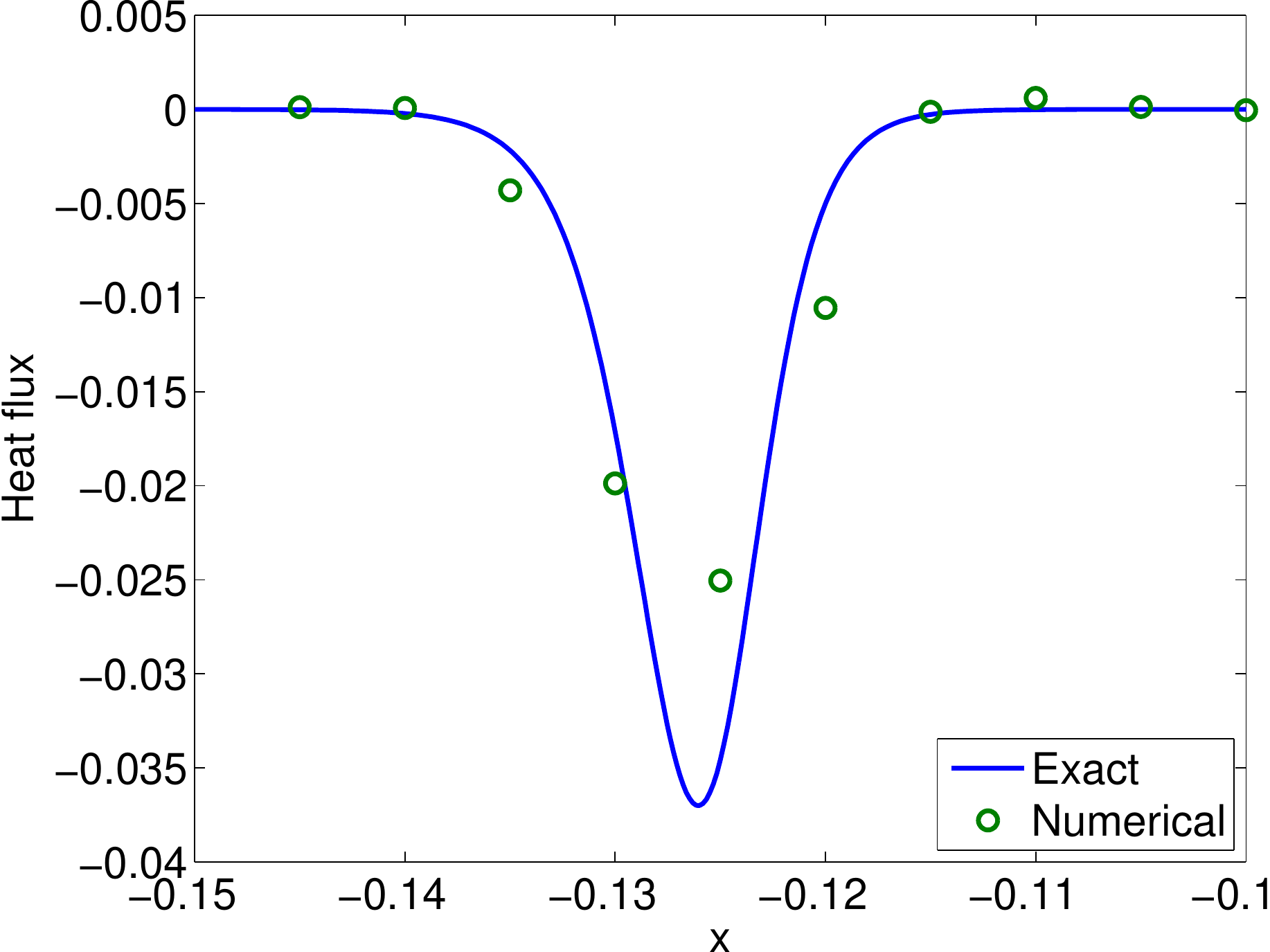}
\caption{NS shock structure: $N=50$ cells, KEP-ES(SD) flux with second and fourth order dissipation}
\label{fig:nsshock1}
\end{center}
\end{figure}

\begin{figure}
\begin{center}
\includegraphics[width=0.34\textwidth]{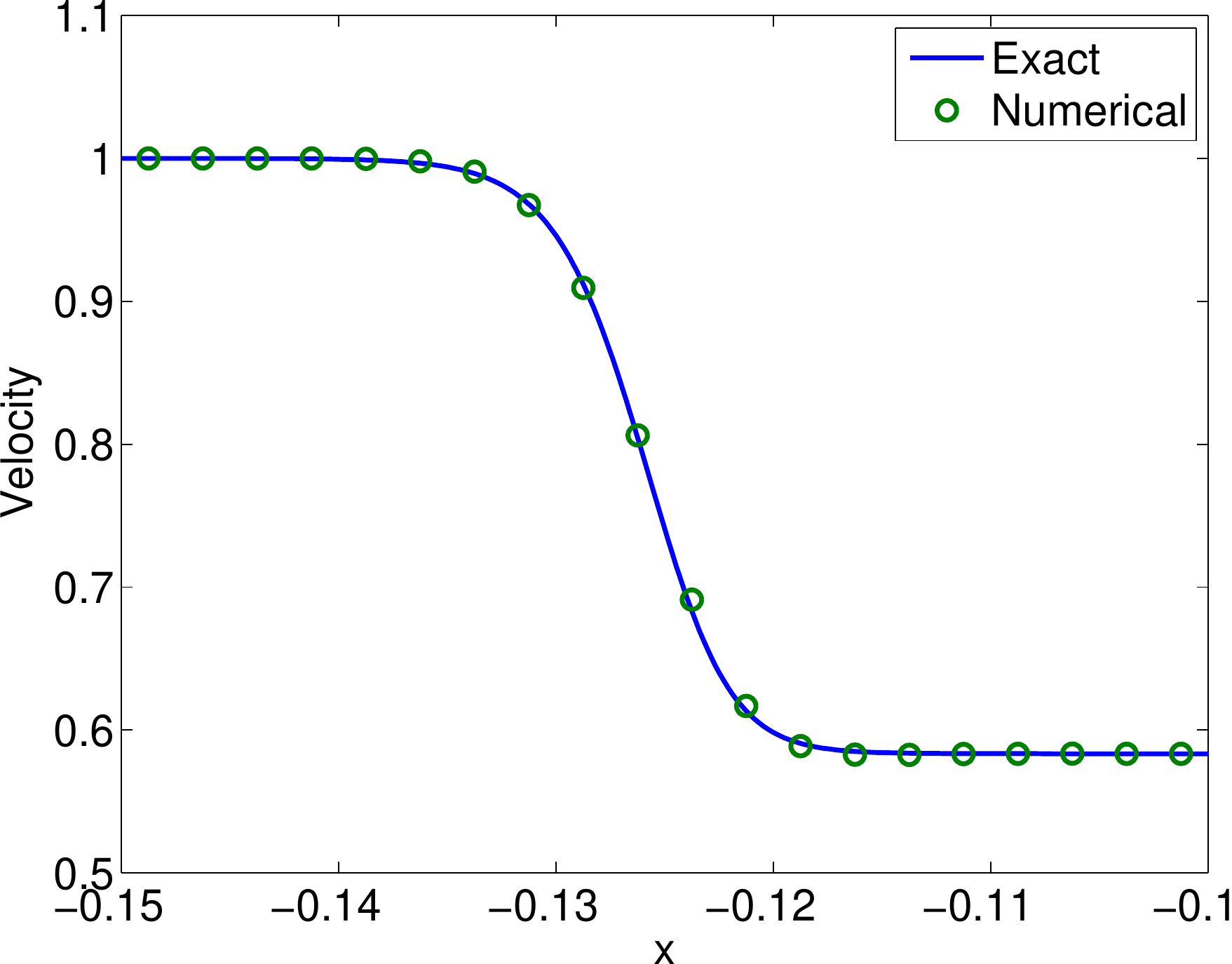}
\includegraphics[width=0.34\textwidth]{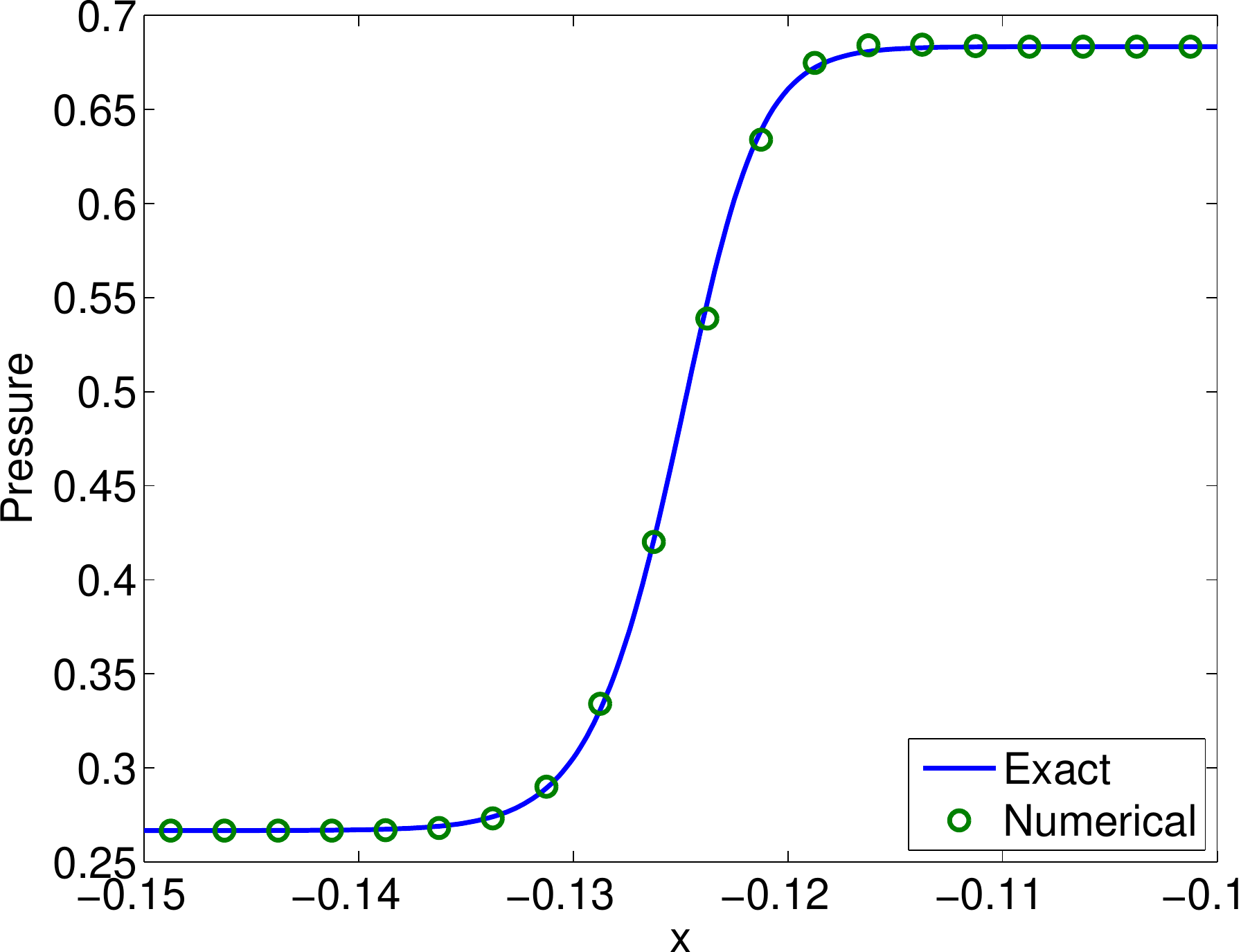} \\
\includegraphics[width=0.34\textwidth]{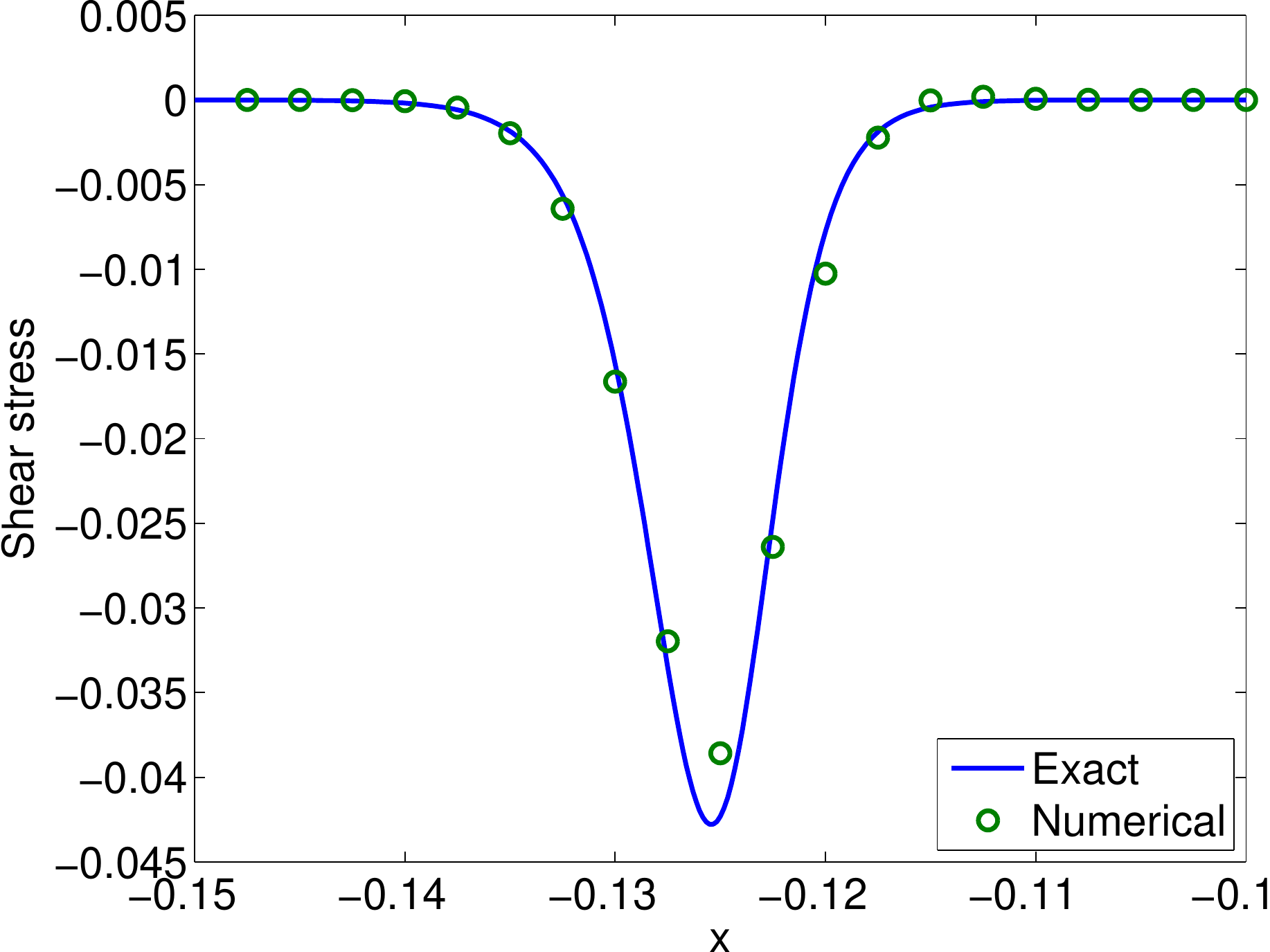}
\includegraphics[width=0.34\textwidth]{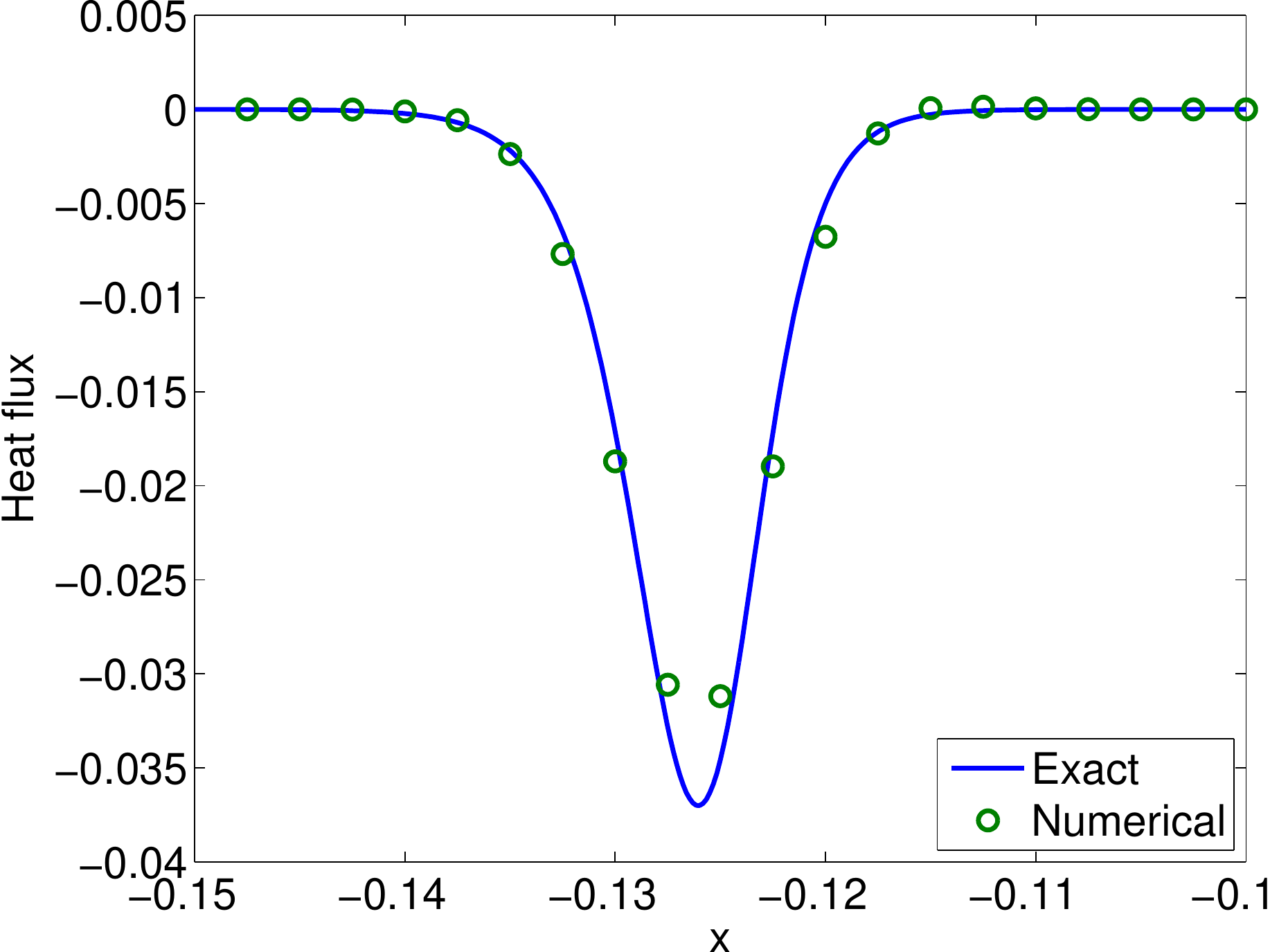}
\caption{NS shock structure: $N=100$ cells, KEP-ES(SD) flux with second and fourth order dissipation}
\label{fig:nsshock2}
\end{center}
\end{figure}

\begin{figure}
\begin{center}
\includegraphics[width=0.34\textwidth]{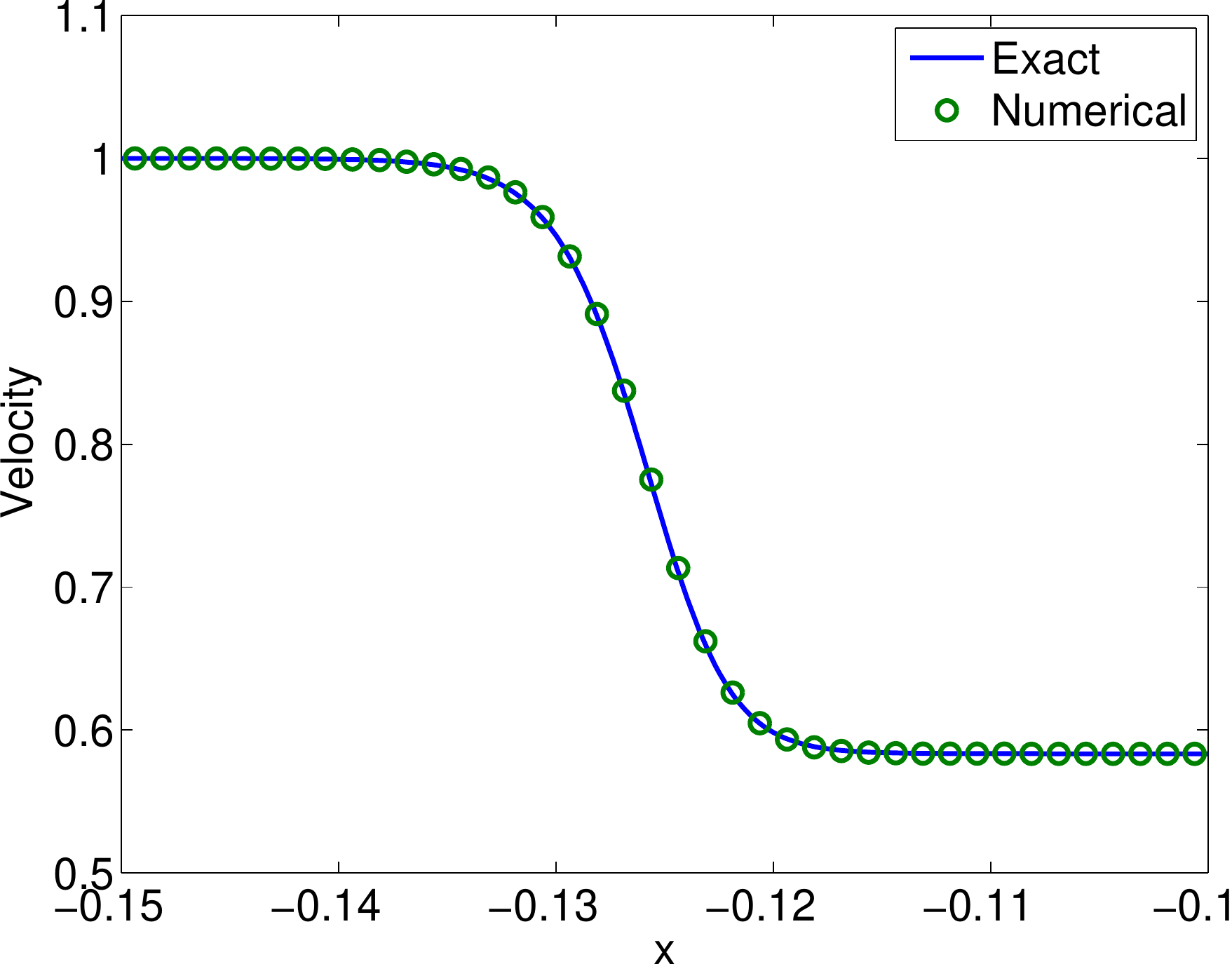}
\includegraphics[width=0.34\textwidth]{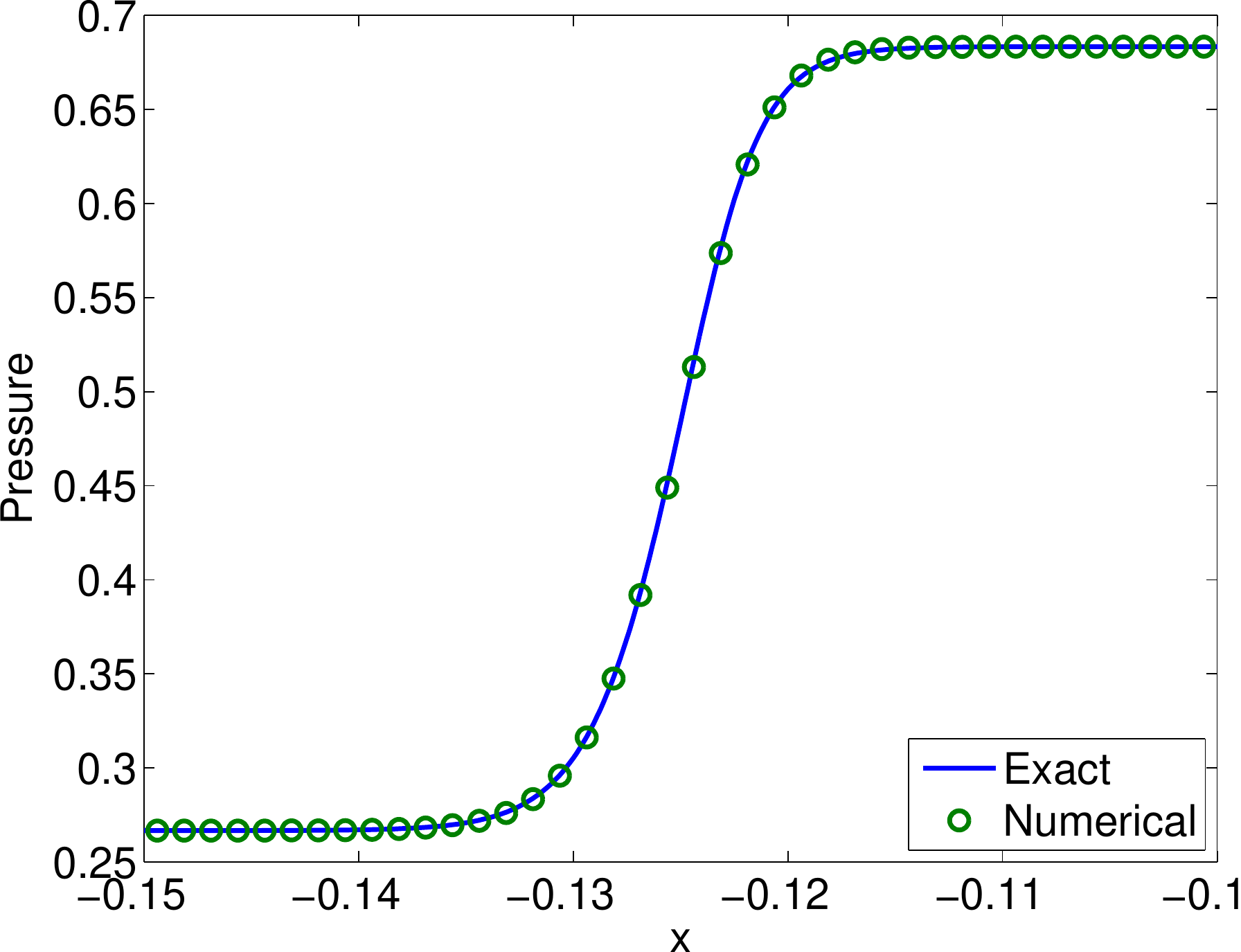}
\includegraphics[width=0.34\textwidth]{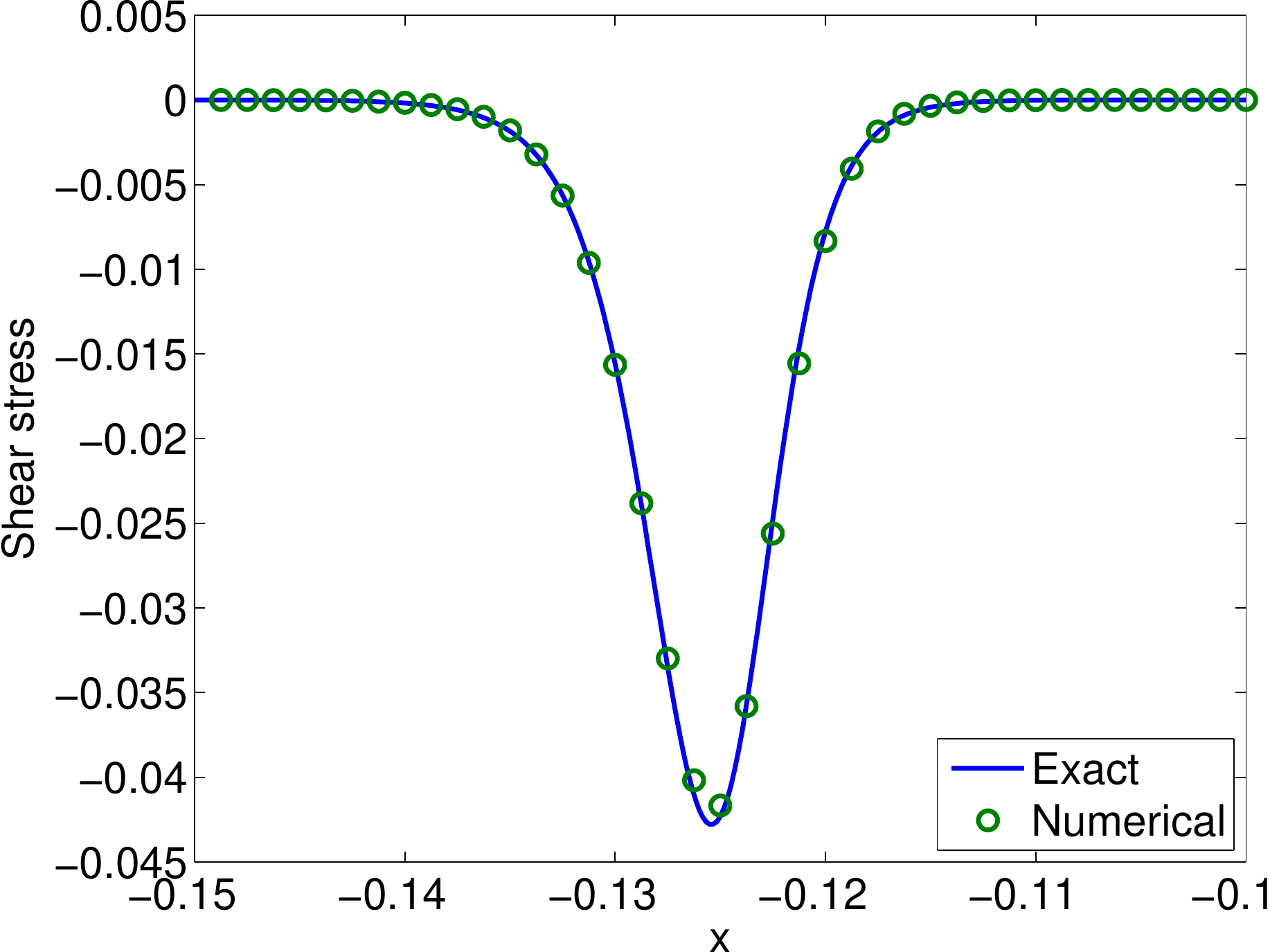}
\includegraphics[width=0.34\textwidth]{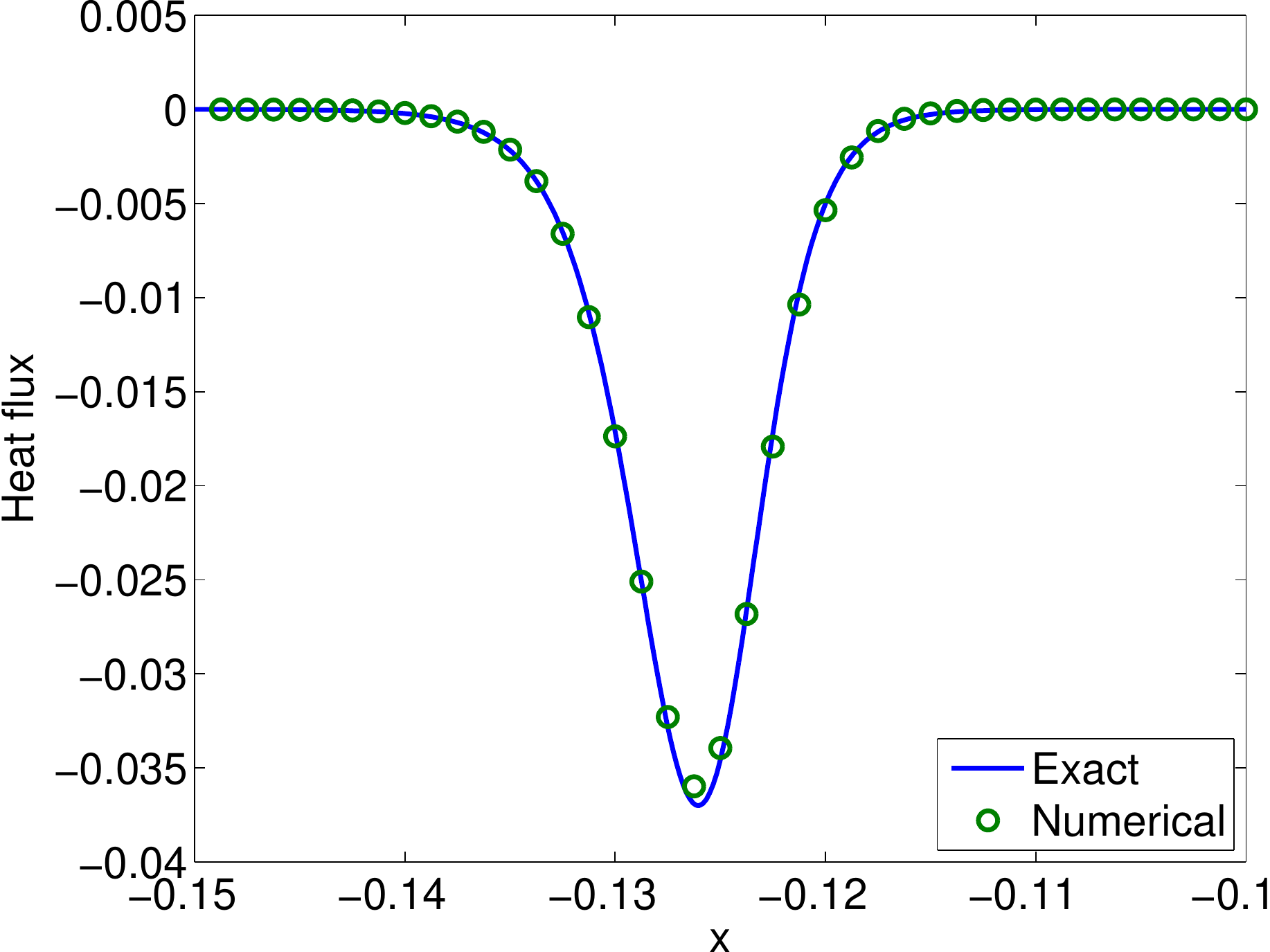}
\caption{NS shock structure: $N=200$ cells, KEP-ES(SD) flux with second and fourth order dissipation}
\label{fig:nsshock3}
\end{center}
\end{figure}

\begin{figure}
\begin{center}
\includegraphics[width=0.34\textwidth]{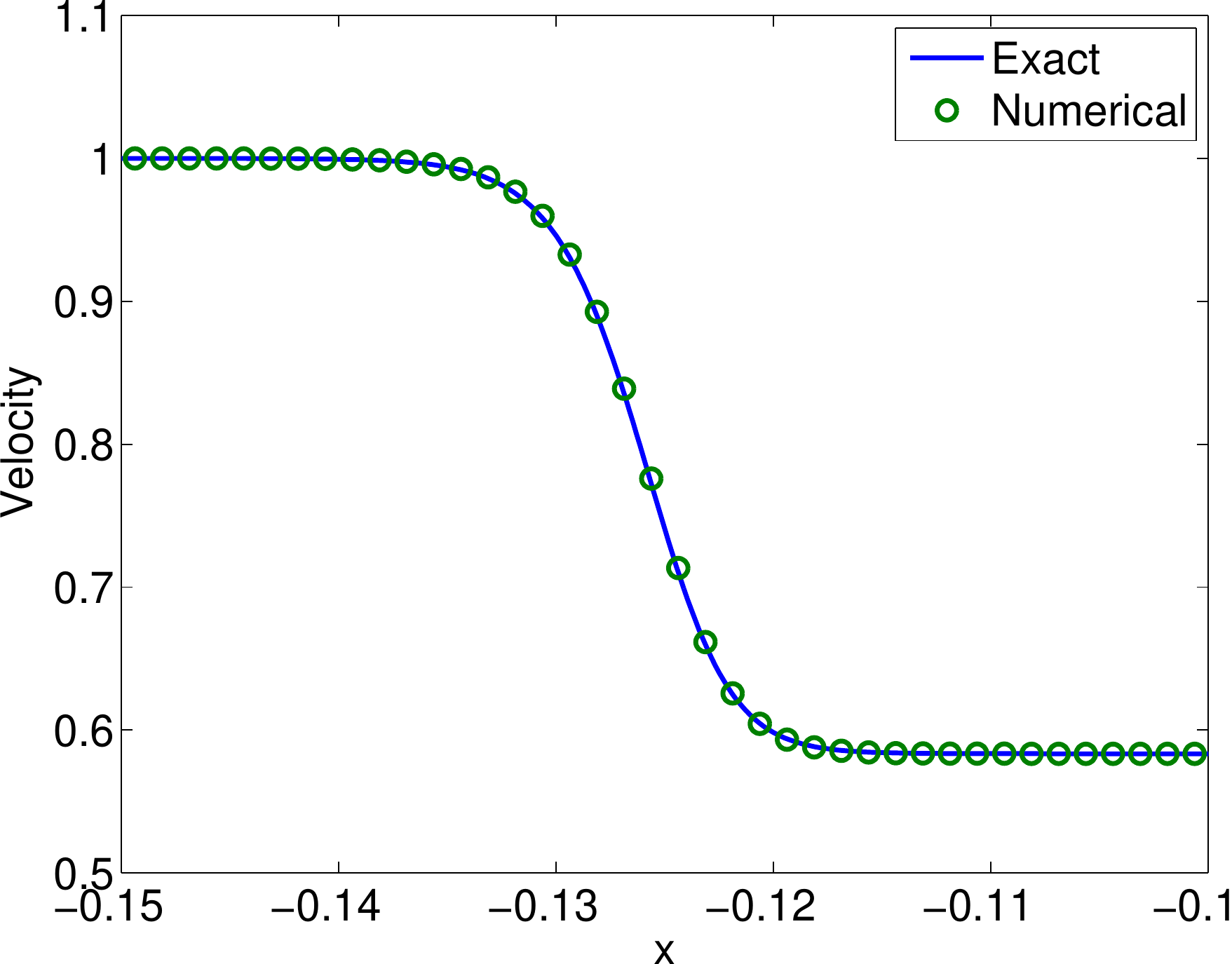}
\includegraphics[width=0.34\textwidth]{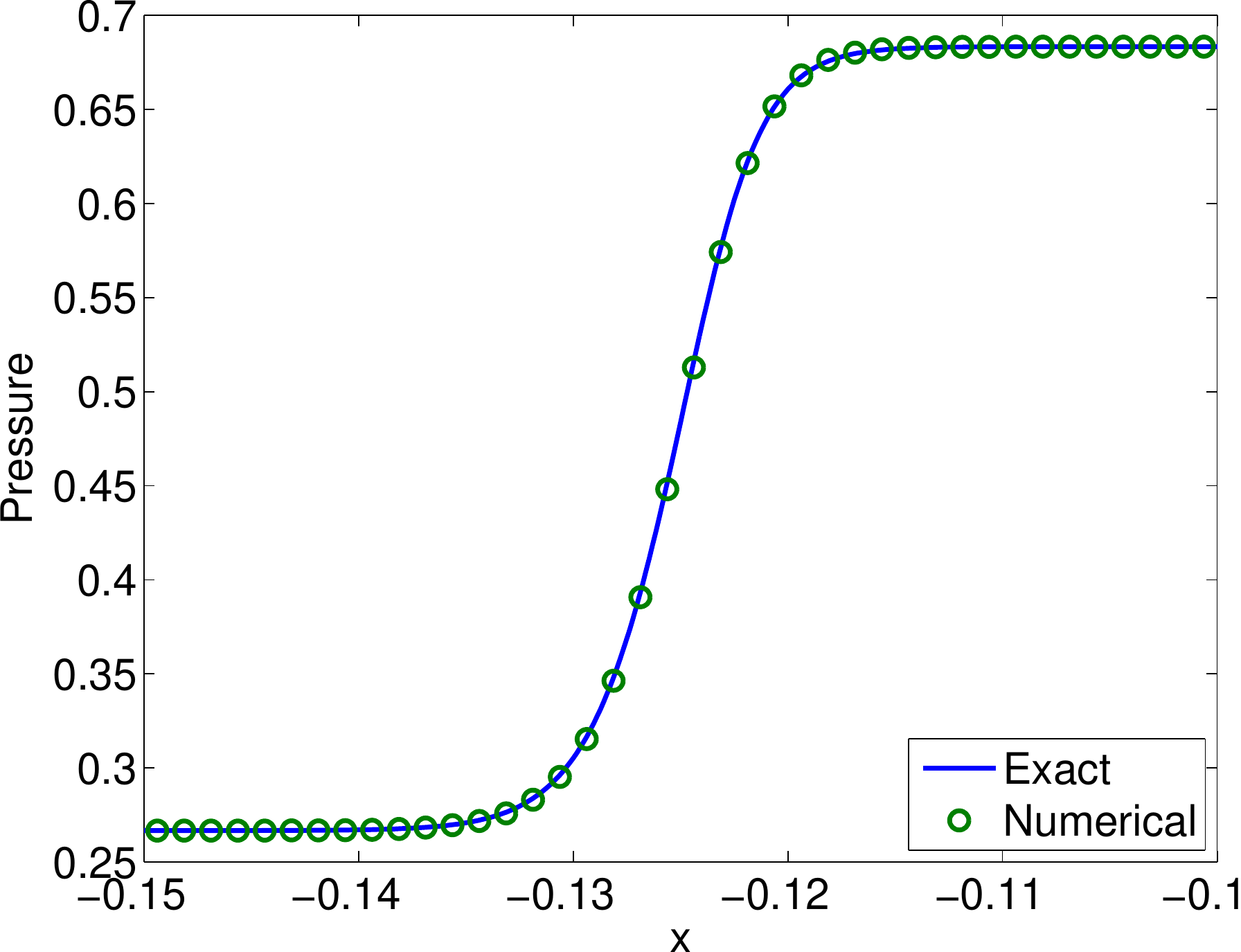}
\includegraphics[width=0.34\textwidth]{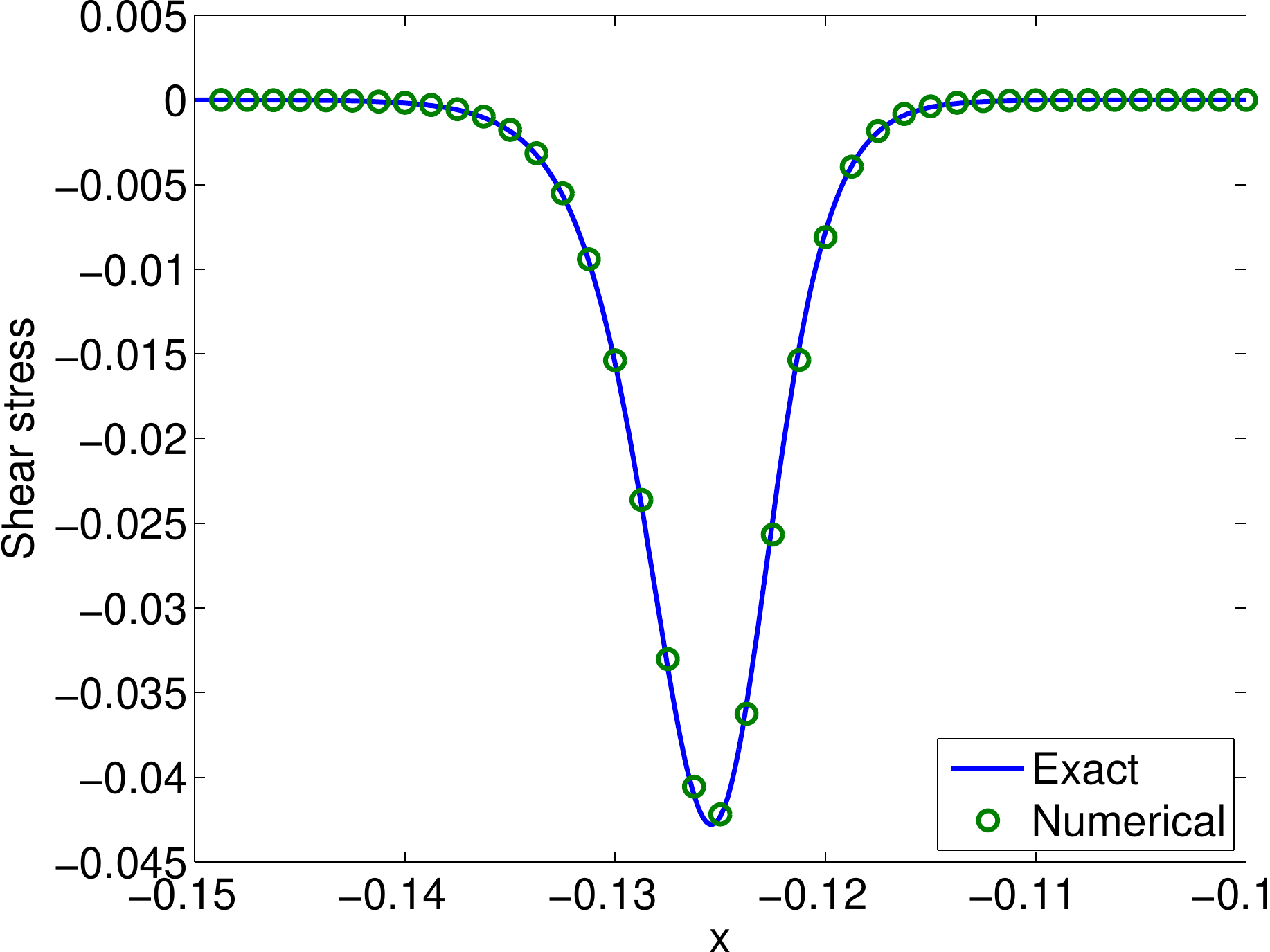}
\includegraphics[width=0.34\textwidth]{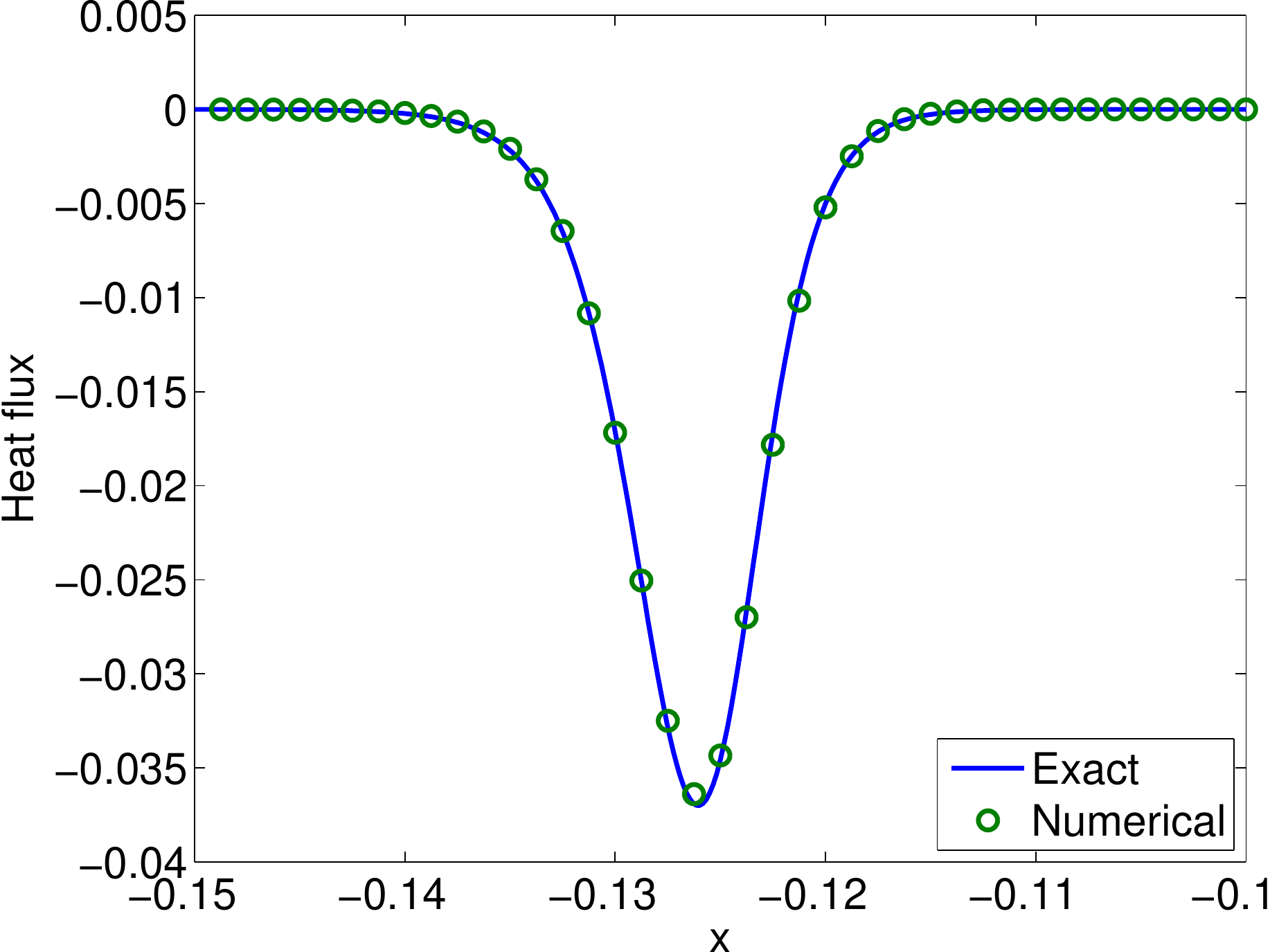}
\caption{NS shock structure: $N=200$ cells, KEP-ES(SD) flux with fourth order dissipation}
\label{fig:nsshock4}
\end{center}
\end{figure}

We compute the shock structure using the Navier-Stokes equations with the help  of the kinetic energy preserving and entropy conservative flux together with scalar artificial dissipation. In the artifificial dissipation we use $\kapt=\frac{1}{2}$ and $\kapf=\frac{1}{25}$. The           parameters defining the problem are: Mach number ahead of the     shock is      $M_1=1.5$, $\gamma=5/3$, $\Pr = 2/3$ while the viscosity law is given  by $\mu  = \mu_1 (T/T_1)^{0.8}$ where the subscript ``1" denotes pre-shock conditions and $\mu_1=0.0005$. Figure~(\ref{fig:nsshock1})-(\ref{fig:nsshock3}) shows the solutions   on $N=50,100,200$ cells. On the coarse mesh, the stress and heat flux cannot be computed accurately since there are too few points inside the shock region, but the solution is still non-oscillatory. On the finer meshes, the scheme is able  to compute the shock structure with good accuracy. Using $N=200$ cells, we compute the solution with only the fourth order dissipation, $\kapf=\frac{1}{200}$ and the solution is shown in figure~(\ref{fig:nsshock4}). We are thus able to obtain accurate solutions on fine meshes for the Navier-Stokes equations using the central kinetic energy and entropy conservative scheme. The fourth order dissipation helps to damp the oscillations in density/pressure which are created due to the initial discontinuity as discussed in the Sod test case.

\subsection{Supersonic flow past cylinder}
\begin{figure}
\begin{center}
\begin{tabular}{ccc}
\includegraphics[width=0.32\textwidth]{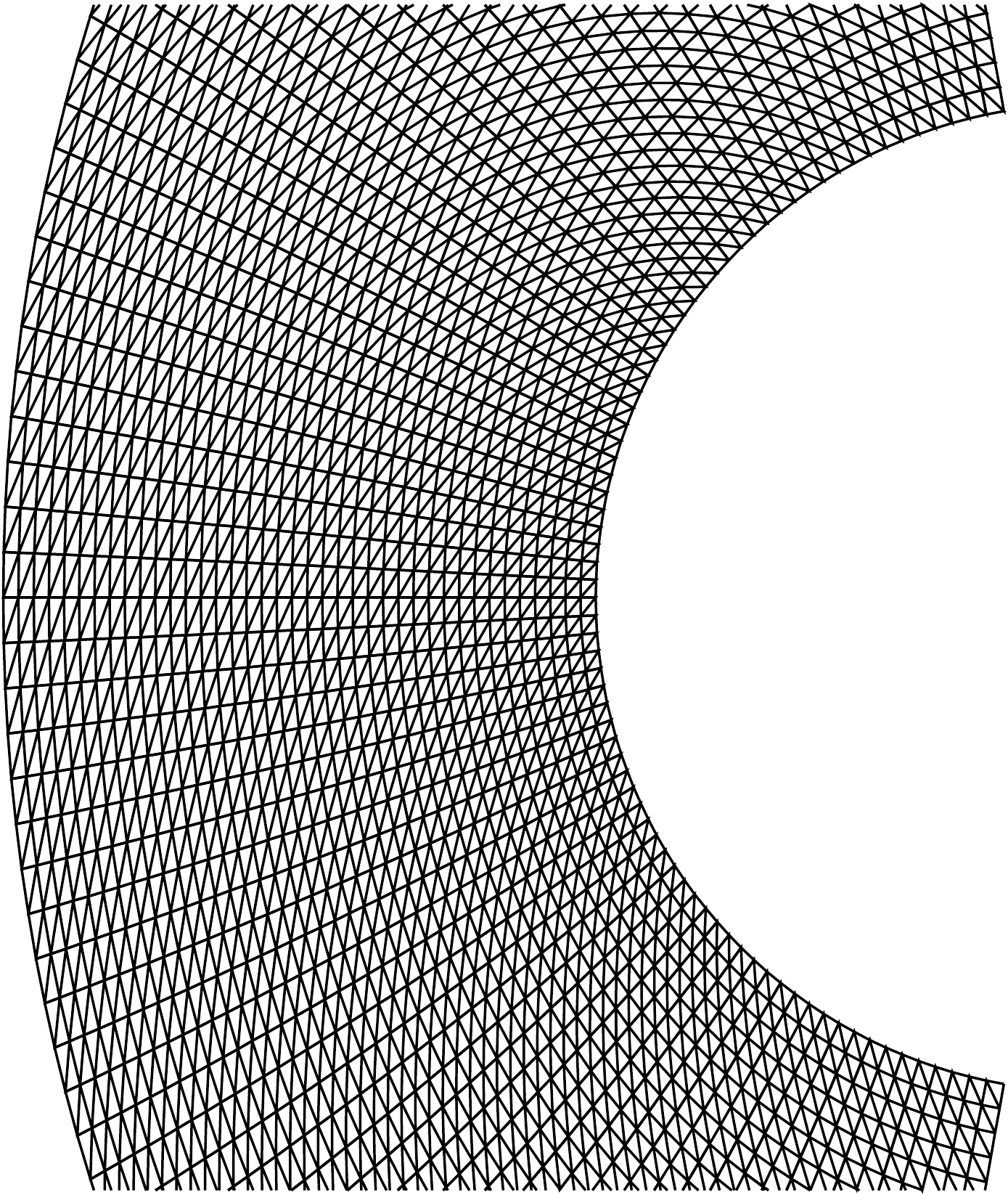} &
\includegraphics[width=0.32\textwidth]{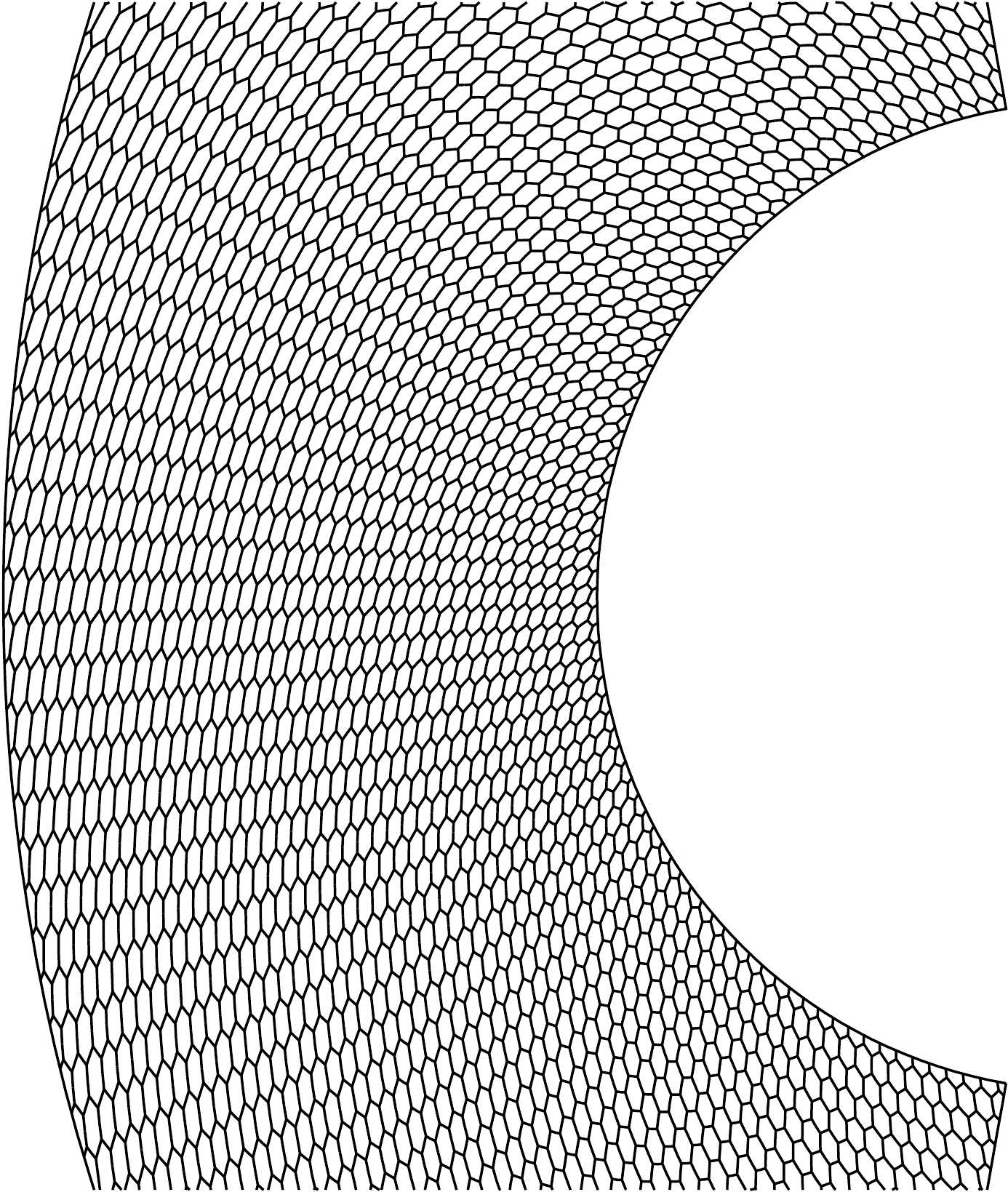} &
\includegraphics[width=0.32\textwidth]{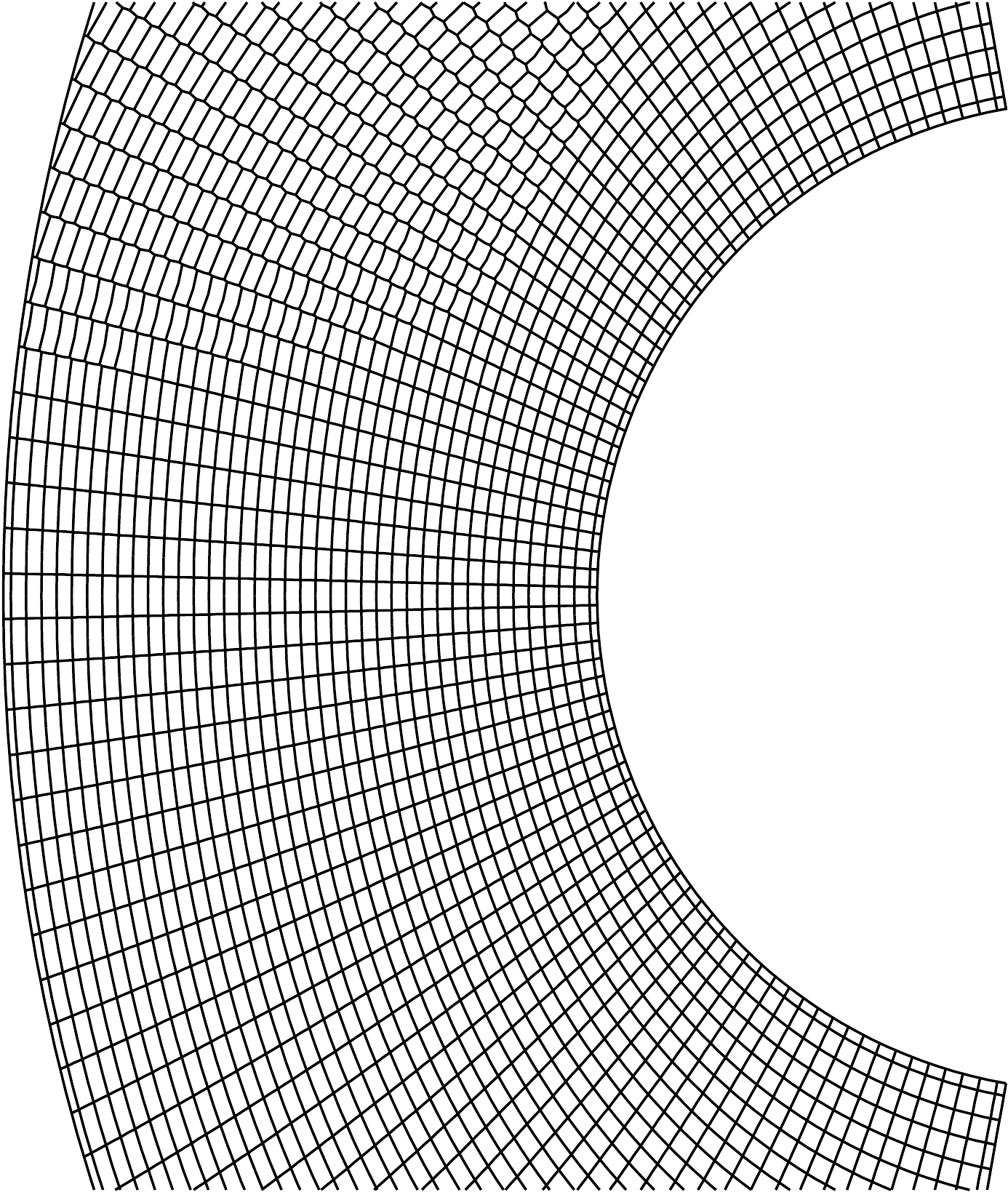} \\
(a) & (b) & (c)
\end{tabular}
\caption{Grid used for supersonic cylinder problem: (a) Primal grid, (b) median dual grid and (c) voronoi dual grid}
\label{fig:cylmesh}
\end{center}
\end{figure}

\begin{figure}
\begin{minipage}[b]{0.49\textwidth}
\begin{center}
\begin{tabular}{cc}
\includegraphics[width=0.28\textwidth]{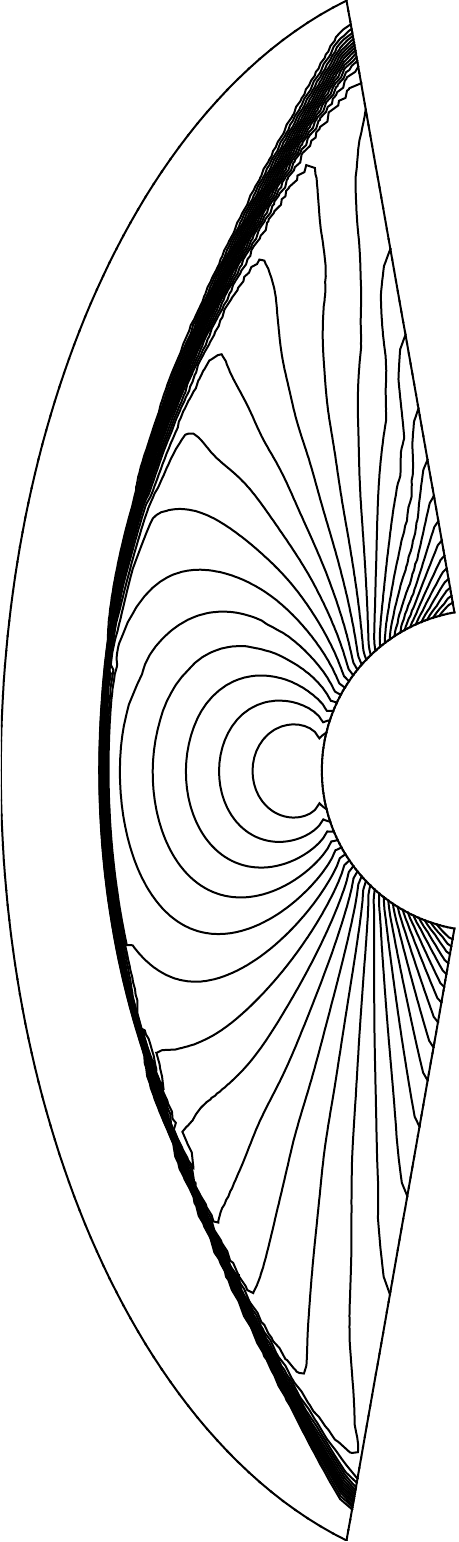} &
\includegraphics[width=0.28\textwidth]{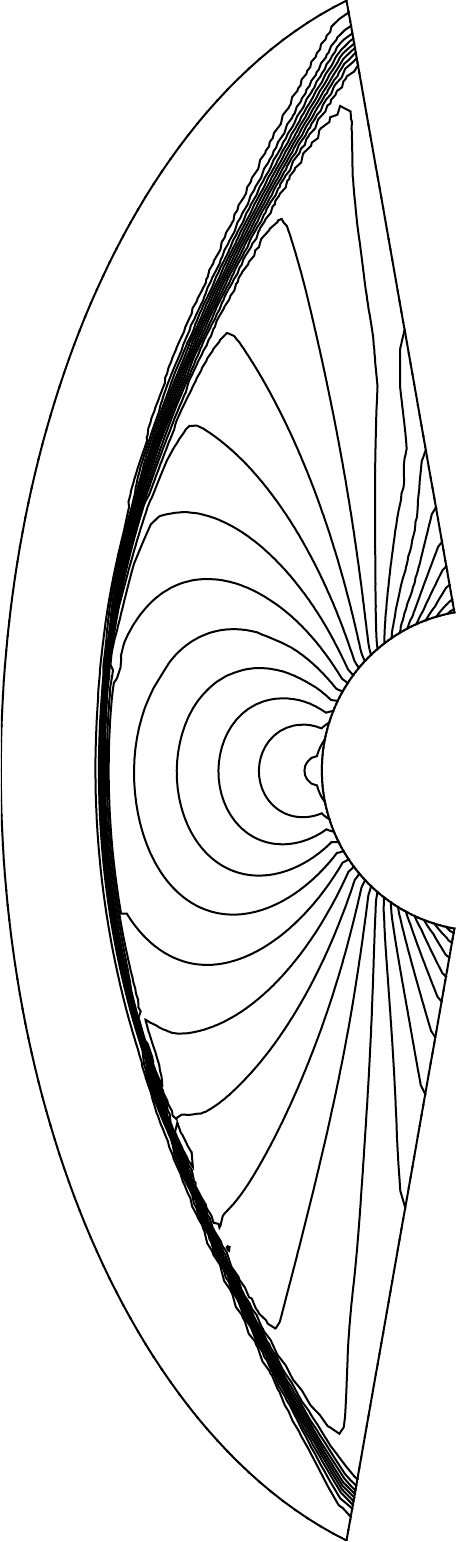} \\
(a) & (b)
\end{tabular}
\caption{Supersonic cylinder, Mach=2: KEP-EC1 flux, Voronoi dual grid, (a) density and (b) pressure}
\label{fig:ec1m2voronoi}
\end{center}
\end{minipage}
\hspace{0.5cm}
\begin{minipage}[b]{0.49\textwidth}
\begin{center}
\begin{tabular}{cc}
\includegraphics[width=0.32\textwidth]{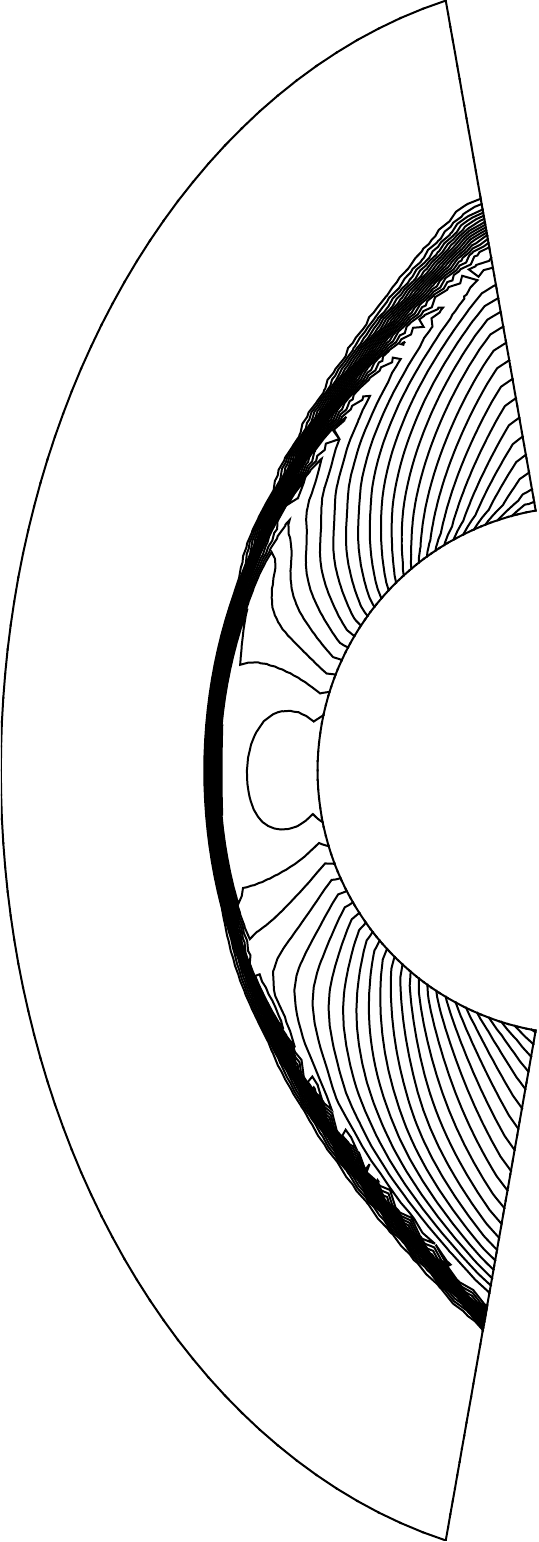}  &
\includegraphics[width=0.32\textwidth]{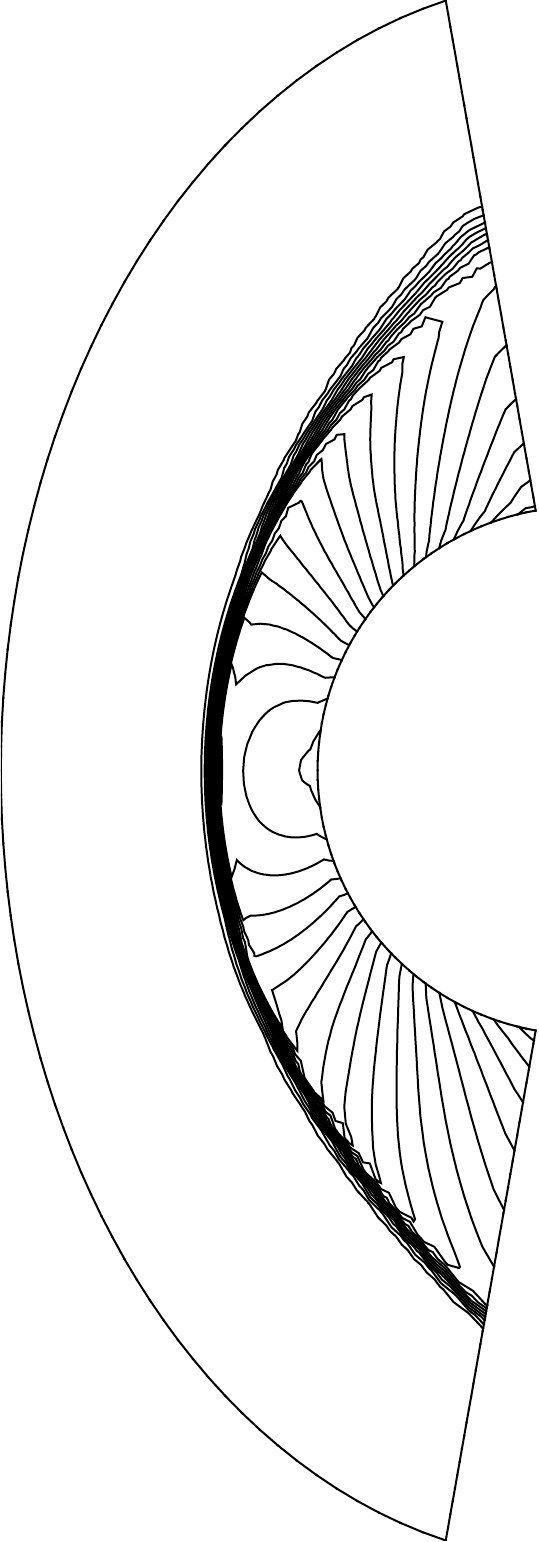} \\
(a) & (b) 
\end{tabular}
\caption{Supersonic cylinder, Mach=20: KEP-EC1 flux, Median dual grid, (a) density and (b) pressure}
\label{fig:ec1m20median}
\end{center}
\end{minipage}
\end{figure}

\begin{figure}
\begin{center}
\begin{tabular}{cccc}
\includegraphics[width=0.18\textwidth]{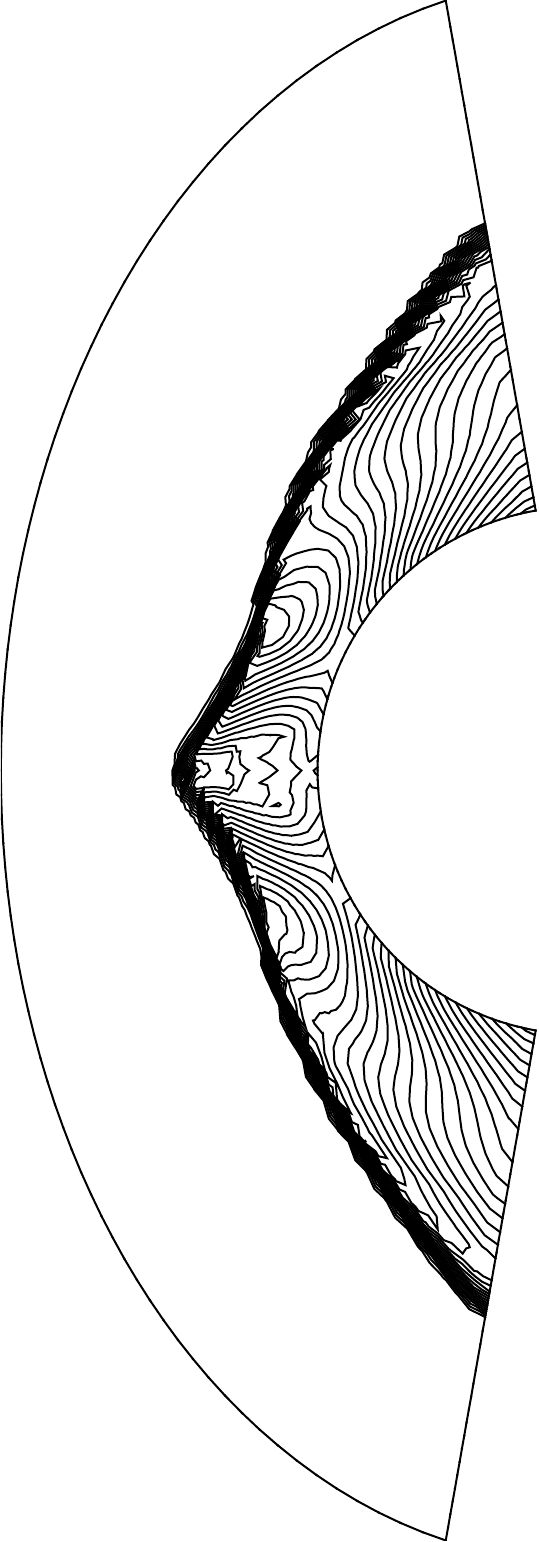} &
\includegraphics[width=0.18\textwidth]{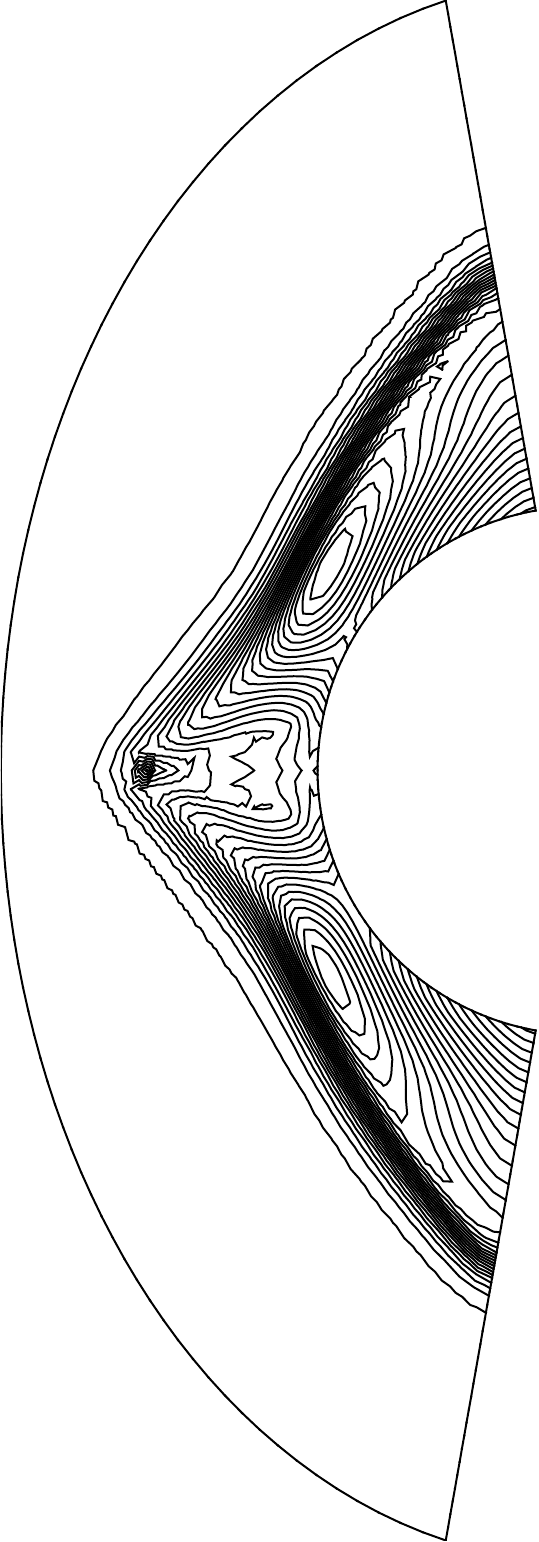} &
\includegraphics[width=0.18\textwidth]{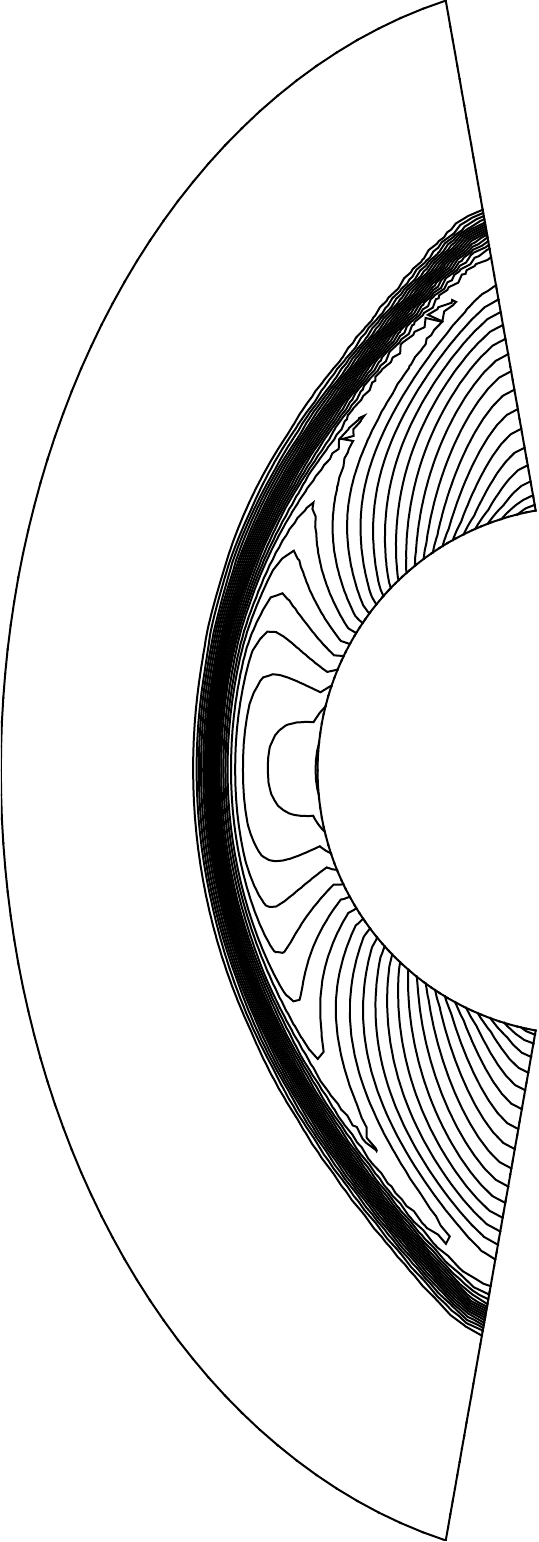} &
\includegraphics[width=0.18\textwidth]{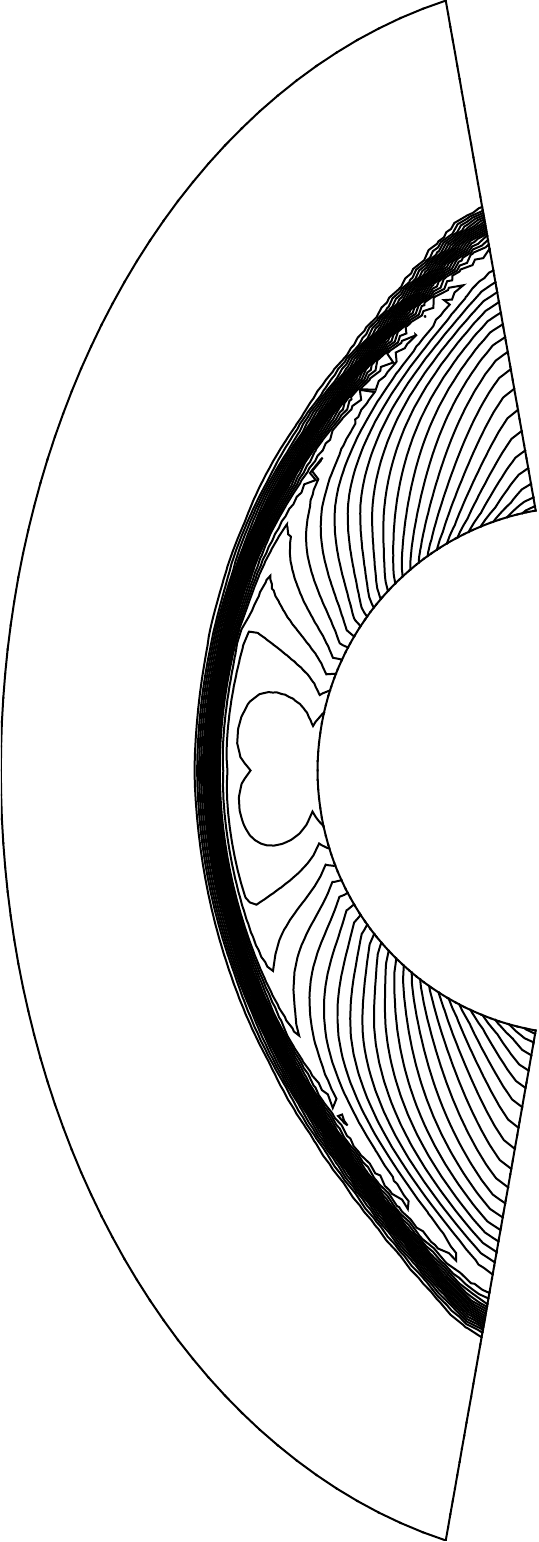} \\
(a) & (b) & (c) & (d)
\end{tabular}
\caption{Supersonic cylinder, Mach=20: Voronoi dual grid, density contours (a) KEP-EC1 (b) KEP-ES(KES) (c) KEP-ES(RUS) (d) KEP-ES(Hyb)}
\label{fig:voronoi}
\end{center}
\end{figure}

Consider the inviscid supersonic flow over a semi-cylinder; the primal triangular grid and corresponding median and voronoi dual meshes are shown in figure~(\ref{fig:cylmesh}). The voronoi cells lead to nearly structured type grids and can thus lead to carbuncle problem since the shock will be aligned with the cell faces in a more exact way than for the median dual cells. For a free-stream Mach number of 2 on the voronoi cells, the solution of the KEP-EC1 scheme is shown in figure~(\ref{fig:ec1m2voronoi}) and there is no problem of carbuncle effect. For a Mach number of 20 and using the median dual grid the KEP-EC1 scheme is able to give carbuncle free solutions with good shock resolution as shown in figure~(\ref{fig:ec1m20median}). However the same scheme gives rise to carbuncle effect on the voronoi mesh as seen in figure~(\ref{fig:voronoi}a). The carbuncle is also seen in the kinetic energy stable scheme KEP-ES(KEPS) in figure~(\ref{fig:voronoi}b). These two schemes resolve stationary contacts exactly and they also exhibit the carbuncle effect. The Rusanov and the new hybrid schemes do not give the carbuncle problem as seen in figures~(\ref{fig:voronoi}c) and (\ref{fig:voronoi}d) respectively. By looking at the shock thickness, it can also be observed that the hybrid scheme is less dissipative than the Rusanov scheme, and this is further confirmed in the boundary layer test case. We also remark that in our numerical experiments the KEP-ES scheme (i.e. without any augmentation of the eigenvalues) does not produce the carbuncle effect even on the voronoi cells but the solution is not monotone as already seen in the 1-D test cases.

\subsection{Transonic flow past NACA-0012 airfoil}
\begin{figure}
\begin{center}
\begin{tabular}{cccc}
\includegraphics[width=0.33\textwidth]{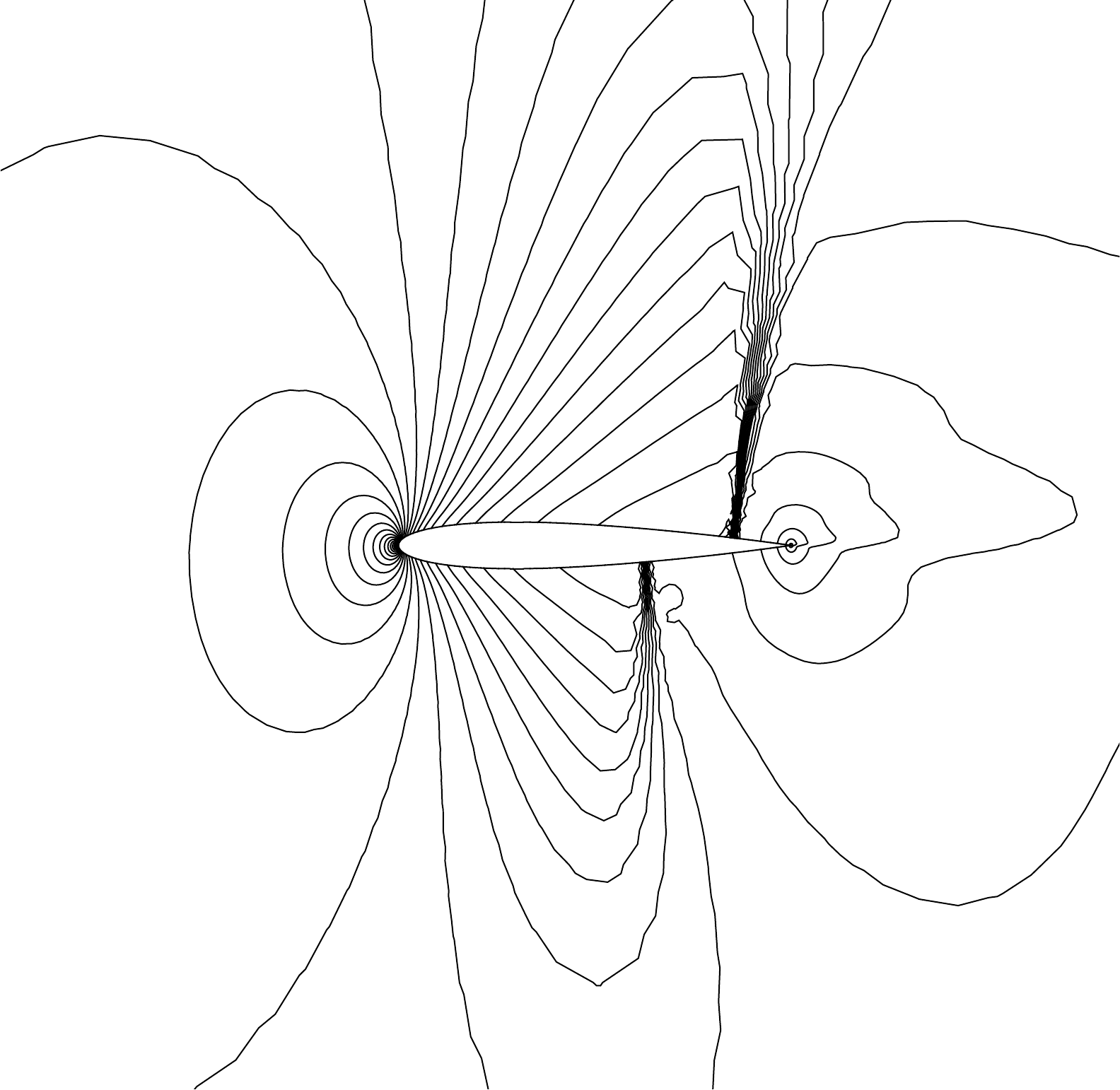} &
\includegraphics[width=0.33\textwidth]{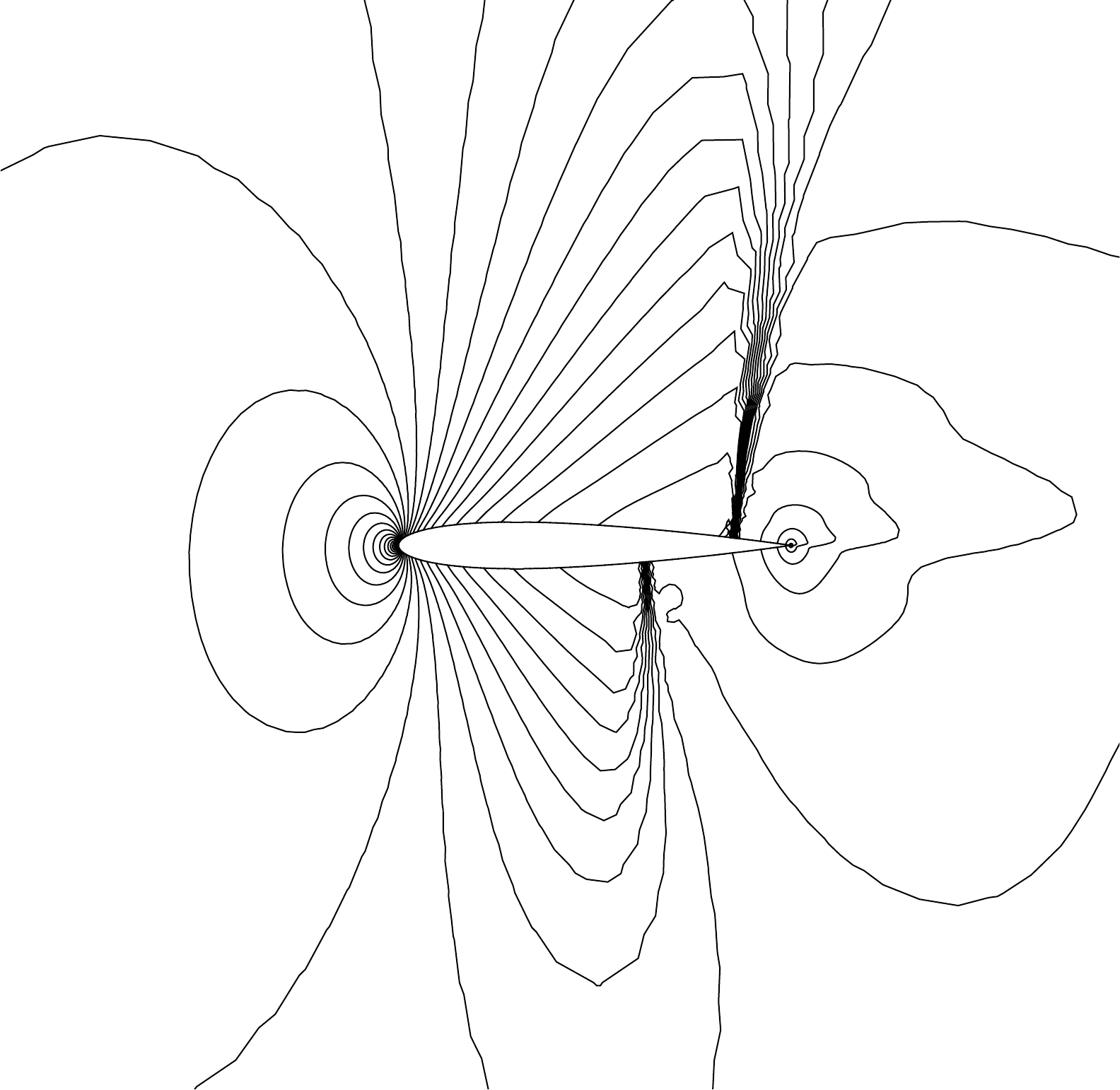} &
\includegraphics[width=0.33\textwidth]{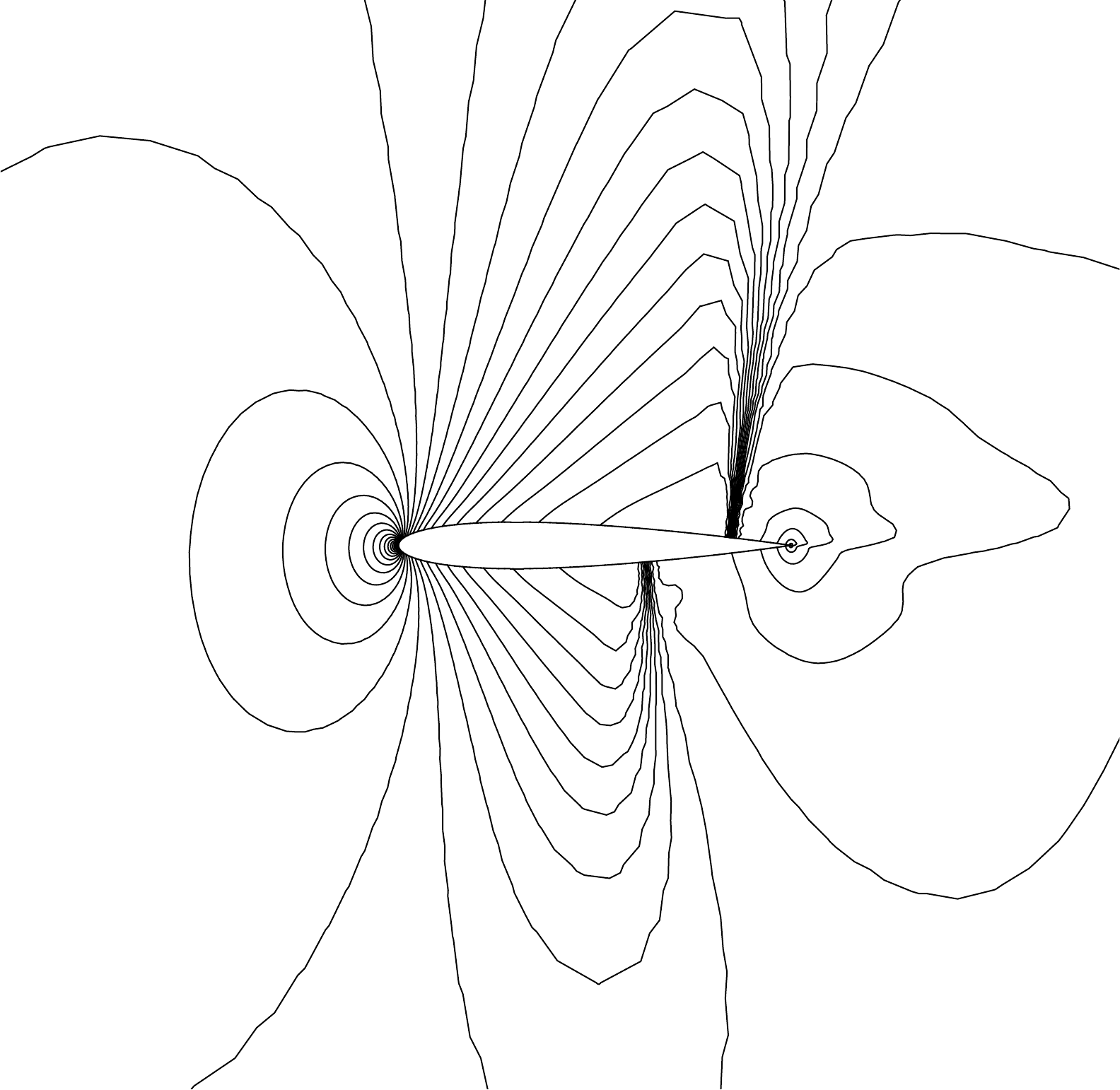} \\
(a) & (b) & (c) & (d)
\end{tabular}
\caption{Transonic flow past NACA0012 airfoil: Mach = 0.85, AOA = 2 deg., median dual grid, Mach contours (a) ROE (b) KEP-EC1 (c) KEP-ES(Hyb)}
\label{fig:naca}
\end{center}
\end{figure}

This is a standard test case for inviscid aerodynamic problems and involves a symmetric NACA-0012 airfoil placed in a freestream Mach number of 0.85 and angle of attack (AOA) of 2 degrees. The flow develops shocks both on the upper and lower airfoil surfaces. We compute this flow on a triangular grid containing 180 points on the airfoil surface and 20 points on the farfield boundary which is a circle, with a total of 6402 vertices. Higher order accuracy is achieved using MUSCL type reconstruction and van Albada limiter. The solution using the original ROE scheme, the new KEP-EC1 and KEP-ES(Hyb) schemes are shown in figures~(\ref{fig:naca}). The results with ROE and KEP-EC1 schemes are almost identical, while the KEP-ES(Hyb) scheme which has increased dissipation shows less sharp resolution of shocks. For such transonic flows, it is preferable to use the KEP-EC1 scheme which is accurate and satisfied entropy condition.
\subsection{Laminar flat-plate boundary layer}
\begin{figure}
\begin{center}
\includegraphics[width=0.9\textwidth]{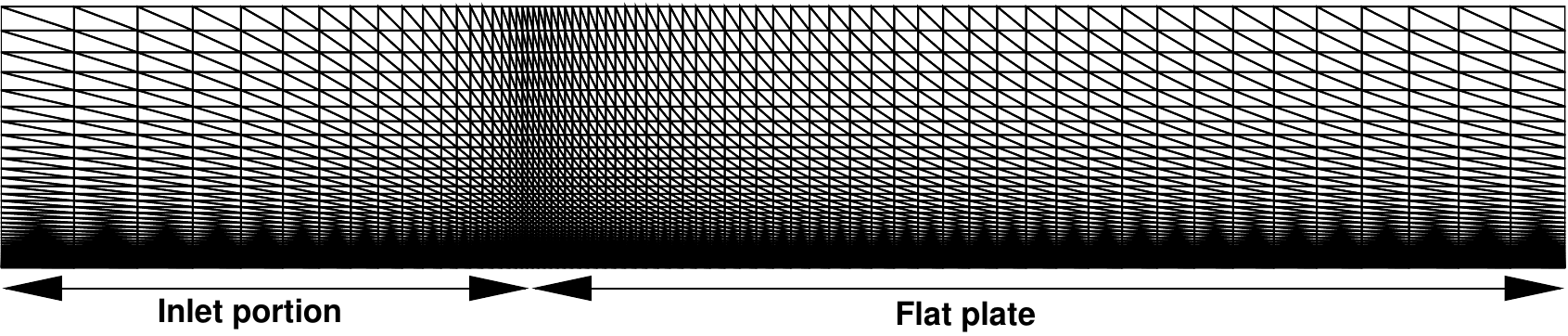}
\caption{Primal grid for laminar flat plate boundary layer problem}
\label{fig:blgrid}
\end{center}
\end{figure}
\begin{figure}
\begin{center}
\begin{tabular}{cc}
\includegraphics[width=0.48\textwidth]{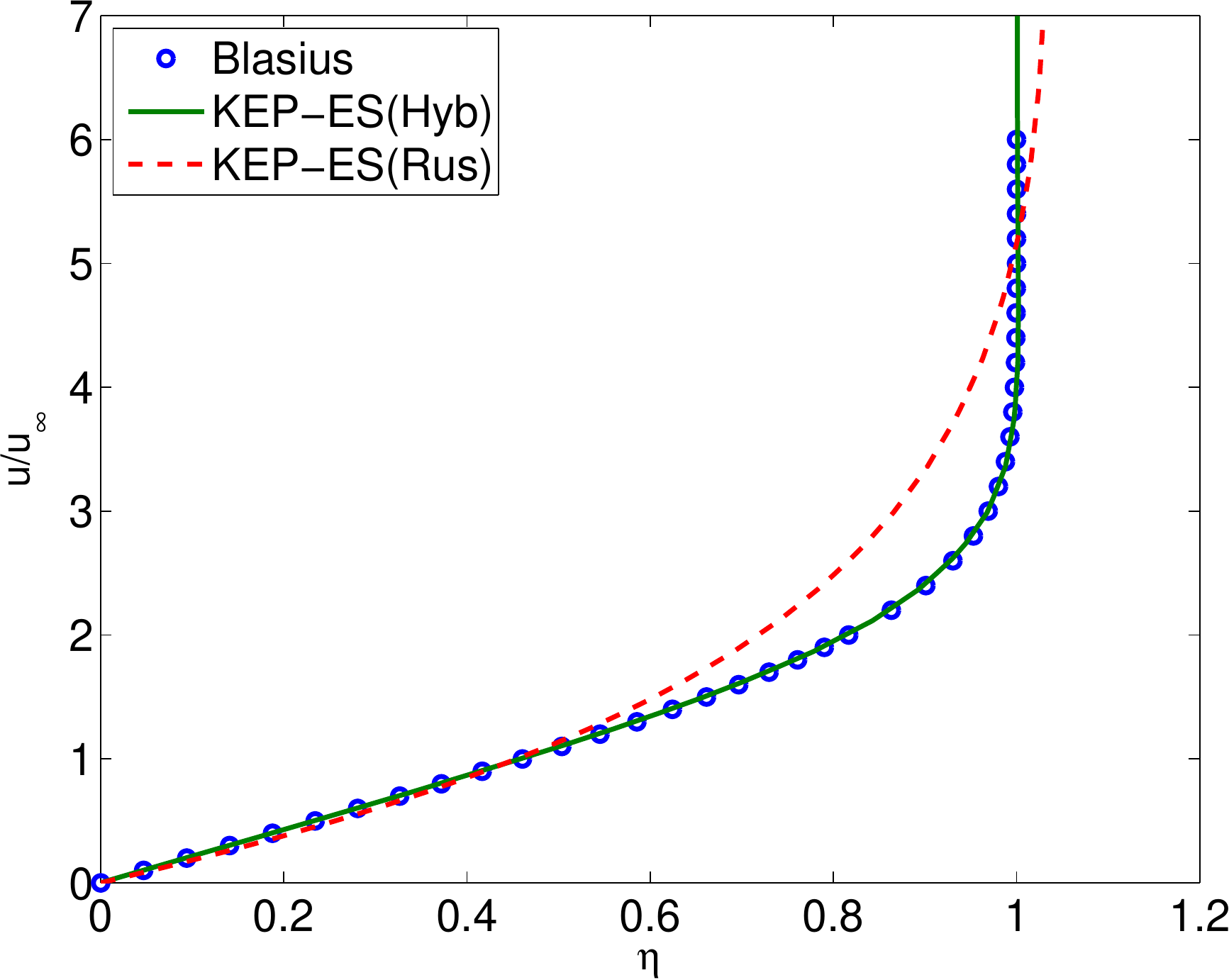} &
\includegraphics[width=0.48\textwidth]{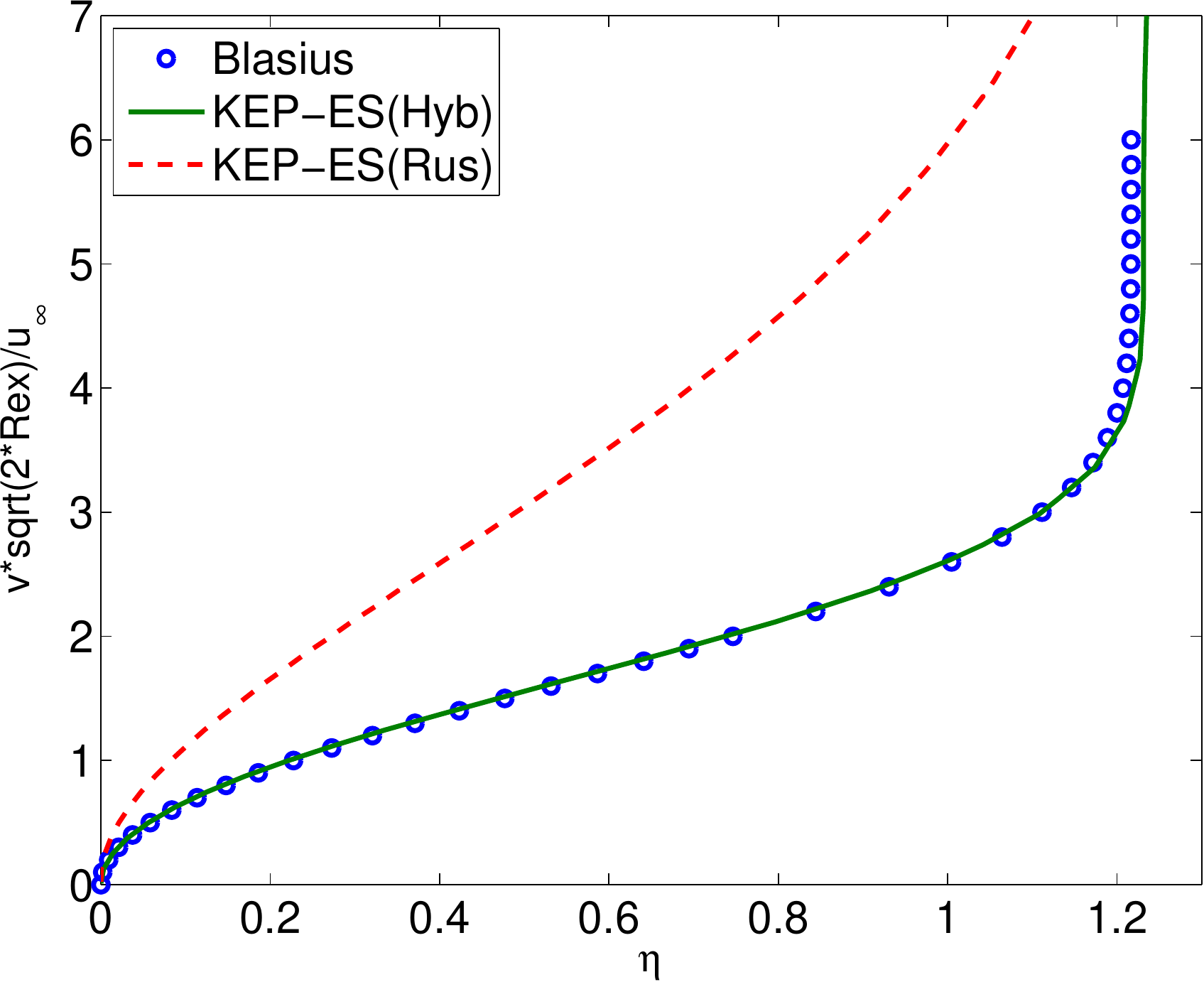} \\
(a) & (b)
\end{tabular}
\caption{Laminar flat plate boundary layer: (a) Streamwise velocity, (b) Vertical velocity}
\label{fig:lfpbl}
\end{center}
\end{figure}

This problem corresponds to viscous flow over a flat plate which leads to the   development of a boundary layer near the plate surface. The Reynolds number     corresponding to the plate length is $10^5$ while the Mach number of the        incoming flow is taken to be 0.1. The computation domain is rectangular as      shown in figure~(\ref{fig:blgrid}) which also shows the primal triangular grid used for the       computations. There is an initial inlet portion of the domain on which slip     boundary condition is imposed followed by the no-slip boundary corresponding to the flat plate. Adiabatic conditions are used on the flat plate boundary. At    the top and outlet, the free-stream pressure is specified while at the inlet    the free-stream values are used together with the numerical flux function to    compute the flux.  We compare the numerical solution of  the velocities with the Blasius semi-analytical solution in figure~(\ref{fig:lfpbl}) in the standard non-dimensional units. These results are taken on the vertical line through the   center point of the plate. We have shown the numerical results using the Rusanov and hybrid form of the dissipation terms. The results from ROE and KEP-EC1 schemes are almost similar to the KEP-ES(Hyb) scheme and are not shown here. It is clear that the new hybrid scheme is able to give accurate resolution of the boundary  layer profile while the Rusanov scheme introduces too much dissipation. The difference is more dramatic in the vertical velocity profile which is more sensitive to numerical dissipation since the vertical velocity component is much weaker than the streamwise component.
\subsection{Step in wind tunnel}
\begin{figure}
\begin{center}
\begin{tabular}{cc}
(a) & \includegraphics[width=0.8\textwidth]{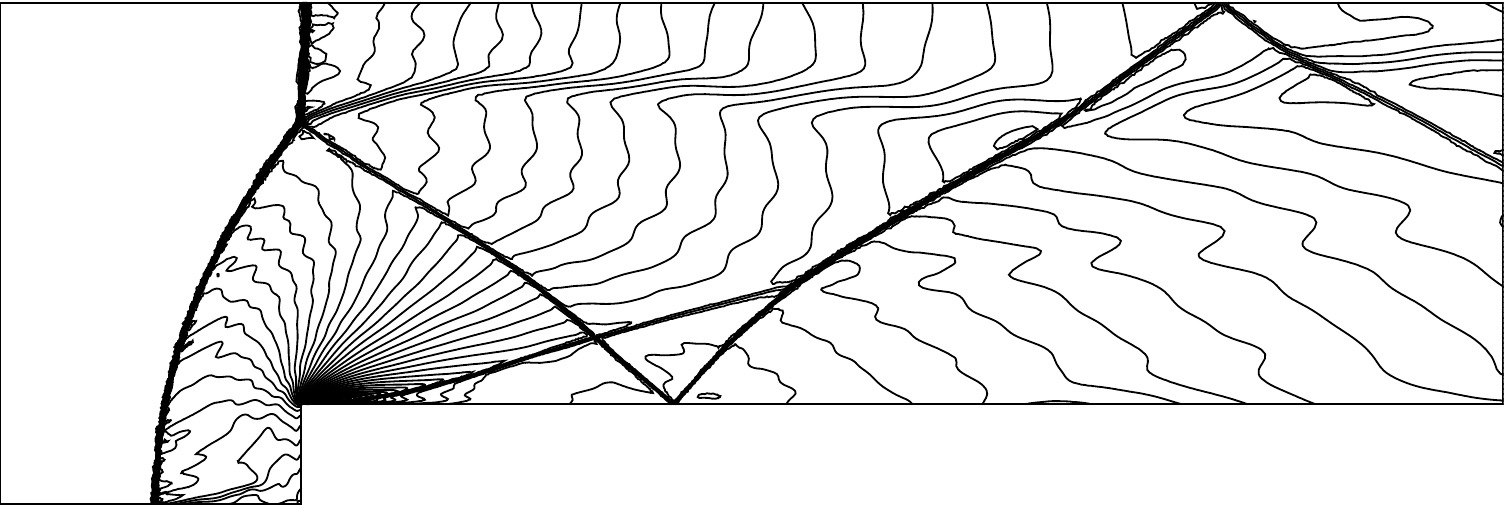} \\
(b) & \includegraphics[width=0.8\textwidth]{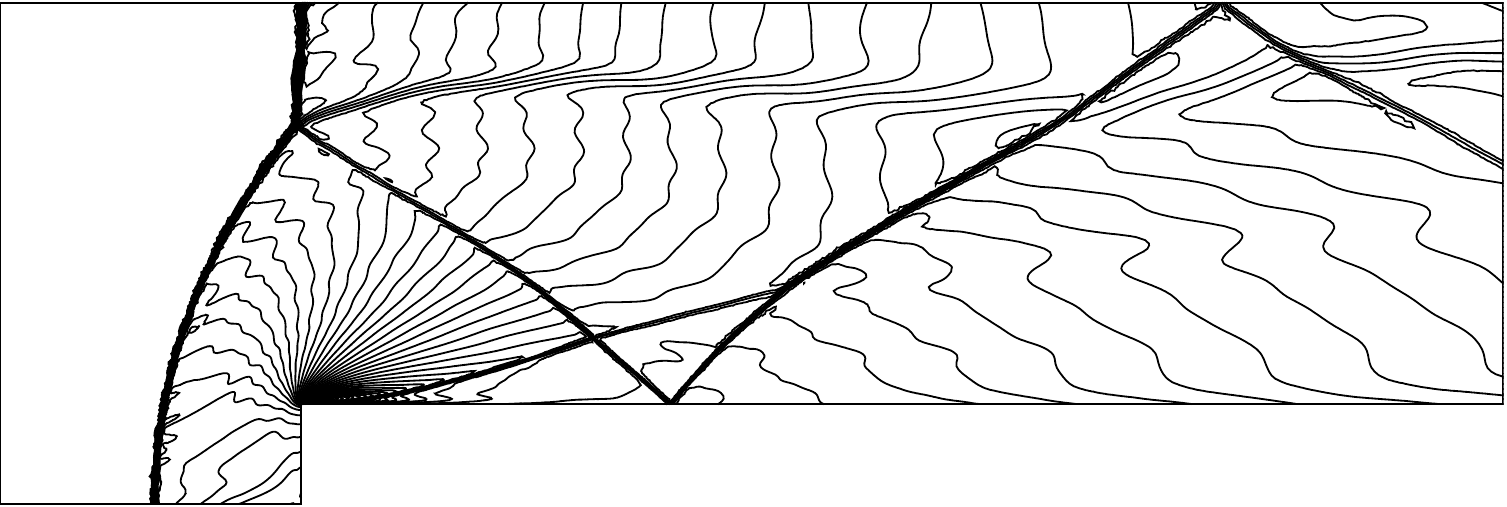}
\end{tabular}
\caption{Forward step in wind tunnel at $M=3$: Mach number contours, 50 equally spaced contours between 0 and 4.8, (a) KEP-EC1 flux, (b) KEP-ES(Hyb) flux}
\label{fig:fstep}
\end{center}
\end{figure}

This test case is described in~\cite{Woodward1984115} and involves inviscid supersonic flow past a step in a wind tunnel which is impulsively started. The initial Mach number is 3. The flow develops several shocks which undergo reflections. A shock triple point intersection leads to the formation of a slip line. The corner on the step is a singular point and many numerical schemes develop spurious entropy at the corner which then leads to the formation of a Mach stem in the downstream direction. This problem is solved using the KEP-EC1 and KEP-ES(Hyb) flux functions together with linear reconstruction of primitive variables and minmax limiter~\cite{barthencyclo}. The grid is adapted to be finer near the corner where the spacing is of $O(0.002)$ while the maximum spacing is of $O(0.01)$, the total number of grid points being 64246. The Mach number contours at time $t=4$ are shown in figure~(\ref{fig:fstep}) using the KEP-EC1 and KEP-ES(Hyb) schemes and both of them are able to resolve the main features of the flow. The slip line is resolved equally well by both schemes and the KEP-EC1(Hyb) does not add any extra dissipation at the contact waves. The reflection of the shock from the botton wall is seen to be free of the spurious Mach stem.

\section{Summary and conclusions}
We have derived a new entropy conservative flux for the Euler equations which also preserves the kinetic energy. Due to their better consistency properties, such numerical fluxes could be attractive for direct numerical simulation of Navier-Stokes equations on highly resolved meshes. These fluxes are central in nature and hence lead to second order accurate schemes. But due to their central character, for coarse meshes or for problems in which the shock is not well resolved by the mesh, some dissipation has to be added to stabilize the scheme. We have constructed scalar and matrix dissipation schemes which preserve the stability properties. In particular, dissipation of kinetic energy is shown to occur if the eigenvalues in the dissipation matrix are chosen appropriately. The matrix dissipation schemes are based on Roe-type dissipation using entropy variables; however even the first order schemes can yield oscillatory solutions even for weak shocks. The eigenvalues can be augmented as proposed by Roe which leads to monotone resolution of shocks even at hypersonic Mach numbers. Due to the entropy consistency of the central part of the numerical flux, all of the entropy consistent schemes avoid unphysical shocks at sonic points. However, like other schemes based on Riemann solvers and which resolve stationary contacts exactly, the new scheme can also suffer from 1-D shock instability and the multi-dimensional carbuncle problem at higher Mach numbers. A hybrid dissipation scheme which blends the Roe type dissipation with the Rusanov dissipation is proposed and shown to be free of carbuncle effect for the blunt body problem. The blending is achieved at the level of the eigenvalues or wave speeds which appear in the dissipation matrix. The hybrid scheme is also able to give good resolution of boundary layer flows which is an important requirement in the computation of viscous flows and is more accurate than the Rusanov type dissipation scheme. Several standard test cases that are commonly used for validation of compressible flows have been computed with the new class of schemes and shown to yield accurate solutions.

Kinetic energy and entropy conservative schemes are attractive for DNS of shock-free compressible turbulent flows. It was shown in the Sod shock tube problem solved with Navier-Stokes equations without any artificial viscosity (Section \ref{sec:nssod})  that the KEP scheme gives smaller entropy oscillations than the entropy conservative schemes, indicating that there is some inherent entropy dissipation in the KEP scheme which is not present in the entropy conservative schemes. However the KEP scheme does not completely eliminate these oscillations and there is a need to add some dissipation to control the oscillations in density and pressure. It is preferable to use the kinetic energy and entropy conservative scheme since it does not add any implicit dissipation; explicit dissipation to control density oscillations can be added whose magnitude is precisely known and hence can be controlled.

\bibliographystyle{abbrv}      
\bibliography{bibdesk}
\end{document}